\documentclass[brochure,english,12pt]{bourbaki}
\usepackage[matrix,arrow]{xy}
\usepackage{amssymb,amsfonts,amsmath,footnote}
\usepackage[francais]{babel}

\usepackage{enumerate}
\usepackage{esint}

\newcommand{\ep}{\varepsilon}

\newcommand{\Z}{\mathbb{Z}}
\newcommand{\R}{\mathbb{R}}

\newcommand{\con}{\equiv}

\newcommand{\ndiv}{\nmid}
\newcommand{\modd}[1]{\; ( \text{mod} \; #1)}
\newcommand{\bstack}[2]{#1 \atop #2}

\newcommand{\maps}{\rightarrow}
\newcommand{\intersect}{\cap}

\newcommand{\union}{\cup}

\newcommand{\al}{\alpha}
\newcommand{\be}{\beta}

\newcommand{\ga}{\gamma}
\newcommand{\del}{\delta}
\newcommand{\Del}{\Delta}

\newcommand{\om}{\omega}

\newcommand{\sig}{\sigma}
\newcommand{\lam}{\lambda}

\newcommand{\Ga}{\Gamma}

\newcommand{\Pcal}{\mathcal{P}}
\newcommand{\Bcal}{\mathcal{B}}

\newcommand{\Ical}{\mathcal{I}}

\newcommand{\Mcal}{\mathcal{M}}

\newcommand{\Scal}{\mathcal{S}}
\newcommand{\Tcal}{\mathcal{T}}

\newcommand{\onebf}{\boldsymbol1}

\newcommand{\hbf}{{\bf h}}

\newcommand{\jbf}{{\bf j}}

\newcommand{\beq}{\begin{equation}}
\newcommand{\eeq}{\end{equation}}
\newcommand{\albf}{\boldsymbol\alpha}
\newcommand{\bebf}{\boldsymbol\beta}

\newcommand{\C}{\mathbb{C}}

\date{Juin 2017}
\bbkannee{69\`eme ann\'ee, 2016-2017}
\bbknumero{1134}
\title{The Vinogradov Mean Value Theorem}
\subtitle{after Wooley, and Bourgain, Demeter and Guth}
\author{Lillian B. PIERCE}
\address{Department of Mathematics\\
Duke University\\
120 Science Drive\\
Durham, NC 27708, U.S.A. }
\email{pierce@math.duke.edu}

\begin{document}
\maketitle

\noindent{\bf INTRODUCTION}

\bigskip

In 1770, Waring wrote: \emph{Omnis integer numerus est quadratus; vel e duobus, tribus vel quatuor quadratis compositus.
Omnis integer numerus vel est cubus; vel e duobus, tribus, 4, 5, 6, 7, 8, vel novem cubis compositus: est etiam quadrato-quadratus; vel e duobus, tribus, \&c. usque ad novemdecim compositus, \& sic deinceps: consimilia etiam affirmari possunt (exceptis excipiendis) de eodem numero quantitatum earundem dimensionum} \cite[Thm. XLVII p. 349]{War70}. From this we extrapolate Waring's problem, the assertion that for each $k \geq 2$, there exists an $s=s(k)$ such that every positive integer $N$ may be expressed as 
\beq\label{N_sum}
 N=x_1^k + \cdots + x_s^k 
\eeq
with $x_1,\ldots, x_s$ non-negative integers.
Hilbert proved this assertion in 1909. In the modern interpretation, Waring's problem also refers to the study of the number $r_{s,k}(N)$ of representations of $N$ in the form (\ref{N_sum}) with $x_i \geq 1$, with the goal of finding 
the least $s=s(k)$ for which an asymptotic  of the form
\beq\label{rskN_asymp}
r_{s,k}(N) = \frac{ \Ga (1+1/k)^s}{\Ga(s/k)} \mathfrak{S}_{s,k}(N) N^{s/k-1}  + O_{s,k}(N^{s/k-1 - \del})
\eeq
holds, for some $\del=\del(s,k)>0$, for all sufficiently large $N$; here $\mathfrak{S}_{s,k}(N)$ is an arithmetic quantity to which we will return to later. In the 1920's, Hardy and Littlewood were the first to prove such an asymptotic valid  for all $k \geq 2$, with $s$ at least exponentially large relative to $k$. Their general approach, via the circle method, relies critically on estimates for exponential sums. In 1935 Vinogradov introduced a new Mean Value Method for investigating such sums, which not only greatly reduced the number of variables required to obtain the asymptotic (\ref{rskN_asymp}), but led to a new record for the zero-free region of the Riemann zeta function, which is still (in terms of its over-all shape) the best-known today.

Despite significant attention paid to sharpening the Vinogradov Mean Value Method since 1935, the cornerstone of the method, to which we will refer  as the Main Conjecture, was not resolved in full until 2015. 
In this manuscript, we explore two approaches to the Main Conjecture:
first, the work of Wooley using analytic number theory, which between 2010--2015 set significant new records very close to resolving the Main Conjecture in all cases, and  resolved it in full in the first nontrivial case;  second, the 2015 breakthrough of Bourgain, Demeter, and Guth using harmonic analysis, which resolved the Main Conjecture in full. Through this new work, connections have been revealed to areas  far beyond arithmetic questions such as Waring's problem and the Riemann zeta function, stretching to core questions in harmonic analysis, restriction theory,  geometric measure theory, incidence geometry, Strichartz inequalities, Schr\"{o}dinger operators, and beyond.

\section{The Main Conjecture in the Vinogradov Mean Value Method}

 Given integers $s,k \geq 1$, let $J_{s,k}(X)$ denote the number of integral solutions to the system of $k$ equations
\beq\label{Vin_sys_dfn}
x_1^j + \cdots + x_s^j = x_{s+1}^j + \cdots + x_{2s}^j, \qquad 1 \leq j \leq k,
\eeq
with $1 \leq x_i \leq X$ for $i=1,\ldots, 2s$. This may be interpreted as a mean value for the exponential sum 
\[ f_k(\albf;X) =  \sum_{1 \leq x \leq X} e(\al_1 x + \cdots + \al_k x^k) \]
upon observing that we may equivalently write 
\beq\label{Jsk_intro_int}
 J_{s,k}(X) = \int_{(0,1]^k} |f_k(\albf;X)|^{2s} d\albf.
 \eeq
 (Here and throughout,  we use the notation $e(t) = e^{2\pi i t}$.)
The foundational conjecture in the area of the Vinogradov Mean Value Method is as follows:\footnote{Here and throughout, we use the Vinogradov notation $A \ll_\ep B$ to denote that there exists a constant $C_\ep$ such that $|A| \leq C_\ep B$. Unless otherwise specified, any statement involving $\ep$ may be taken to hold for all arbitrarily small $\ep>0$, with associated implied constants.  In addition, the implied constant is allowed to depend on other parameters, such as $s,k$ in this section.}
\begin{conj}[The Main Conjecture]
For all integers $s,k \geq 1$, 
\beq\label{Jsk_MC0}
 J_{s,k}(X) \ll_{s,k,\ep} X^\ep ( X^s + X^{2s-\frac{1}{2}k(k+1) }),
 \eeq
for all $X \geq 1$, and every $\ep>0$.
\end{conj}
This conjecture may be refined to omit the factor $X^\ep$ if $k>2$ (see  \cite[Eqn. 7.5]{Vau97}, and  \S \ref{sec_k12}--\S \ref{sec_ep_removal} in this manuscript).
 Vinogradov's  motivation for bounding the mean value $J_{s,k}(X)$ was to extract bounds for individual sums $f_k(\albf;X)$, which (as we will summarize later) would impact many number-theoretic problems, the most famous relating to Waring's problem and the Riemann zeta function. Vinogradov's mean value perspective \cite{Vin35}, evidently inspired by an idea of Mordell \cite{Mor32}, has been influential ever since.

Historically, any nontrivial result toward the Main Conjecture has been called the Vinogradov Mean Value Theorem. Instead, we will refer to partial results as the Vinogradov Mean Value Method, and  reserve the terminology ``Vinogradov Mean Value Theorem''  for the newly proved theorems that verify the Main Conjecture in full.

The so-called critical case occurs when $s=s_k=\frac{1}{2}k(k+1)$: this is the index at which the two terms in (\ref{Jsk_MC0}) are of equivalent size. Importantly, if for a certain $k$ the bound (\ref{Jsk_MC0}) has been proved for $s_k$, then it follows immediately for all $s \geq 1$. In particular, it is elementary to verify the cases of $k=1,2$ at the critical index. For a review of these facts, and heuristics leading to the Main Conjecture, see \S \ref{sec_MC_classical}.

Due to the role of the critical index, investigations naturally divide into the case of small $s$ and large $s$. For small $s \leq k$, the diagonal solutions (those with $\{x_1, \ldots, x_s\} = \{x_{s+1}, \ldots, x_{2s}\}$ as sets) dominate, and the relation 
\[J_{s,k}(X) = s! X^s + O(X^{s-1})\]
 is an immediate consequence of the Newton-Girard identities. Hua \cite{Hua47} extended these considerations to verify the upper bound in the Main Conjecture when $s=k+1$, subsequently refined to an asymptotic in \cite{VauWoo95}.  See e.g. \cite[\S 3]{Woo14b} for further refinements and historic partial results for small $s$, before 2010.

We turn to the setting of large $s$.
 Vinogradov's  original work \cite{Vin35} was taken up by Linnik \cite{Lin43}  (who moved it to a $p$-adic setting) and polished by Karatsuba \cite{Kar73} and Stechkin \cite{Ste75}.
In total, this approach showed that for $s \geq k$, 
 \beq\label{Jsk_classical}
  J_{s,k}(X) \leq D(s,k) X^{2s - \frac{1}{2} k(k+1) + \eta_{s,k}} 
  \eeq
 for an explicit constant $D(s,k)$, and $\eta_{s,k} = \frac{1}{2}k^2 (1-1/k)^{[s/k]} \leq k^2 e^{-s/k^2}$ for $k \geq 2$. As a consequence of the decay of $\eta_{s,k}$,  one can verify  for $s \geq 3k^2 (\log k + O(\log \log k))$ the bound in the Main Conjecture, and indeed obtain an asymptotic 
 \beq\label{Jsk_asymp}
  J_{s,k}(X) \sim C(s,k) X^{2s - \frac{1}{2}k(k+1)} 
  \eeq
 for an explicit positive real constant $C(s,k)$.
 (Indeed, along the same lines of argument, the leading $3$ can be improved to a $2$; see \cite[Thm. 3.9]{ACK04}.)
 See \cite[Ch. 3]{ACK04} for a treatment of these various historic methods, or \cite[Ch. 7]{Vau97} for a modern overview.

In the 1990's, Wooley's thesis \cite{Woo92,Woo96} developed an \emph{efficient differencing} method which allowed him to extract faster decay from $\eta_{s,k}$ for $s > k^2 \log k$, and as a result he obtained the next historic leap, showing  the Main Conjecture held for $s \geq k^2 (\log k + 2 \log \log k +O(1))$.

\subsection{The work of Wooley: Efficient Congruencing}
This record remained untouched until the 2010's, when Wooley developed an \emph{efficient congruencing} method. In his initial work on this method \cite{Woo12a}, Wooley set a startling new record, proving that the Main Conjecture held for $s \geq k(k+1)$, for every $ k \geq 3$ (and in the asymptotic form (\ref{Jsk_asymp}) for $s \geq k(k+1)+1$).  This was a landmark result, since for the first time the additional logarithmic factor was removed, and thus the limitation on $s$ was only a constant multiple away from the expected truth.
 With this new method (and its ``multigrade'' version), Wooley  held the Main Conjecture under siege, making continual progress on this and related consequences in a remarkable series of papers, including \cite{Woo12b, Woo13b,Woo14b,Woo15a,Woo15b,Woo16c,Woo16a,Woo16b,Woo17b,Woo17a}, and in joint work with Ford \cite{ForWoo14}.
 
 By the end of 2015, Wooley had succeeded in proving the conjectured bound for $J_{s,k}(X)$ for $s$ pushing very close to the critical index $s_k = \frac{1}{2}k(k+1)$. For $s$ approaching the critical index from above, Wooley proved the Main Conjecture for $k \geq 3$ and $s \geq k(k-1)$ \cite{Woo14b}. For $s$ approaching the critical index from below,  Wooley  \cite{Woo17b} proved the Main Conjecture for $1 \leq s \leq D(k)$, where $D(4)=8,$ $D(5)=10,...$ and 
\beq\label{Wooley_final}
D(k) \leq \frac{1}{2}k(k+1) - \frac{1}{3} k - 8 k^{2/3},
 \eeq
 for large $k$.
This landmark result was the first ever to prove the Main Conjecture for $s$ differing from the critical index $s_k$ by only a lower order term.

Moreover, in \cite{Woo16c}, Wooley proved the $k=3$ case in full, establishing the Main Conjecture for the first nontrivial degree:

\begin{theo}[Wooley:  Vinogradov Mean Value Theorem, $k=3$]\label{thm_VMVT_cubic}
For $k=3$, for every integer $s \geq 1$, the Main Conjecture holds,
\[ J_{s,k} (X) \ll_{s,\ep} X^\ep(X^s + X^{2s - \frac{1}{2}k(k+1)}),\]
for all $X \geq 1$, and every $\ep>0$.
\end{theo}

\subsection{The work of Bourgain, Demeter and Guth: $\ell^2$ decoupling}

In December 2015, Bourgain, Demeter and Guth \cite{BDG16} announced the resolution of the final cases required for the Main Conjecture for $k \geq 4$:
\begin{theo}[Bourgain, Demeter, Guth: Vinogradov Mean Value Theorem, $k \geq 4$]\label{thm_VMVT}
For every integer $k \geq 4$ and for every integer $s \geq 1$,
\[ J_{s,k}(X) \ll_{s,k,\ep} X^\ep ( X^s + X^{2s - \frac{1}{2} k(k+1)}),\]
for all $X \geq 1$, and every $\ep>0$.
\end{theo}
By standard methods, once Theorems \ref{thm_VMVT_cubic} and \ref{thm_VMVT} are known, the $X^\ep$  may be omitted and the asymptotic (\ref{Jsk_asymp}) obtained  for all integers $k \geq 3$ and $s > \frac{1}{2}k(k+1)$; see \S \ref{sec_ep_removal}.

The resolution of the Main Conjecture, eighty years after its initiation by Vinogradov, is  a spectacular achievement with many consequences. An additional striking feature is that the Bourgain-Demeter-Guth approach is rooted in harmonic analysis. Of course, even the expression (\ref{Jsk_intro_int}) immediately indicates that the Main Conjecture is inextricably bound to ideas of Fourier analysis, as is the Hardy-Littlewood circle method, one of the core techniques of analytic number theory.
But as we will see, the Bourgain-Demeter-Guth method (and work leading up to it) takes this quite a bit further, revealing fascinating connections between the Vinogradov Mean Value Method and deep open problems that have been motivating work  in harmonic analysis over the past fifty years.

 The Bourgain-Demeter-Guth work \cite{BDG16} is part of the new area of \emph{decoupling} (and in particular $\ell^2$ decoupling). This area was initiated in work of Wolff \cite{Wol00}, who introduced the study of an $\ell^p$ decoupling inequality  for the cone, motivated by the local smoothing conjecture for the wave equation, see e.g. the survey \cite{Sog95}. See also  early work on decoupling, then called ``Wolff's inequality,'' by {\L}aba, Pramanik, and Seeger  \cite{LabWol02,LabPra06,PraSee07}.
 
In the present context, decoupling was deeply developed by Bourgain \cite{Bou13} and then Bourgain and Demeter \cite{BouDem13a,BouDem14x,BouDem14xb,BouDem15a,BouDem15x,BouDem15,BouDem16a,BouDem16b,BouDem17a}, and is now a growing area of research, including for example \cite{BouWat15x,Bou14x,Bou16x,Bou17a,BDG16x,DGS16x,DGL16x,FSWW16,DGG17x,Guo17x}.
As a whole, the decoupling method has ramifications far broader than resolving the Main Conjecture in the Vinogradov Mean Value Method, several of which we will mention in \S \ref{sec_mot_cons} and \S \ref{sec_decoupling_other}.

While decoupling has deep ties to many aspects of harmonic analysis, it is worth noting (as per Wooley \cite{Woo17b}) that the decoupling method shares similarities with Vinogradov's initial framework for mean value investigations in the 1930's, which used small real intervals. Indeed, it seems reasonable to speculate, as Wooley does, that efficient congruencing (in its most recent ``nested'' formulation) and $\ell^2$ decoupling will ultimately be understood as $p$-adic and Archimedean perspectives of one unified method (see \S \ref{sec_parallels} for a few clear parallels).

\subsection{Overview of the paper}

In \S \ref{sec_mot_cons} we survey several of the key arithmetic problems that motivated the Vinogradov Mean Value Method, and state the new records imposed by the recent work. In \S \ref{sec_MC_classical} we gather classical observations about Vinogradov systems and the Main Conjecture, before turning in \S \ref{sec_EC} to an exploration of Wooley's efficient congruencing ideas, focusing on the cases of $k=2$ and $k=3$ for the purposes of demonstration. 

We then turn to the setting of decoupling. In order to motivate the form that decoupling theorems take, in \S \ref{sec_decoupling_intro} we familiarize the notion of decoupling with a few simple examples, and also place decoupling questions within the broader context of historic questions in harmonic analysis. In \S \ref{sec_decoupling_VMVT} we state the Bourgain-Demeter-Guth $\ell^2$ decoupling theorem for the moment curve, and demonstrate how it implies the Main Conjecture, followed by a brief mention of decoupling in more general settings. In \S \ref{sec_Kakeya} we introduce multilinear methods and Kakeya problems, and demonstrate in the simple case of $k=2$ how these play a critical role in proving the $\ell^2$ decoupling theorem. Finally in \S \ref{sec_anatomy} we survey the anatomy of the Bourgain-Demeter-Guth proof. Given the boundary-crossing nature of the new work, we have taken particular care to make both arithmetic and analytic motivations and methods clear to a wide audience.

\section{The Vinogradov Mean Value Method: Motivations and Consequences}\label{sec_mot_cons}
In this section we gather together a few classical arithmetic problems that both motivated the Vinogradov Mean Value Method, and are now decisively affected by the resolution of the Main Conjecture and the underlying methods: Waring's problem, bounds for Weyl sums, and the Riemann zeta function.

\subsection{Waring's problem and the Hardy-Littlewood circle method}
In Waring's problem, the least $s$ such that \emph{every} $N$ has a representation as a sum of $s$ non-negative $k$-th powers is traditionally called $g(k)$. After Hilbert proved the existence of such a $g(k)$ for every $k$ in 1909 \cite{Hil09},  the study of reducing $g(k)$ to its optimal lowest value was taken up by Hardy and Littlewood in the \emph{Partitio Numerorum} series. But $g(k)$ is controlled by a few pathological numbers; indeed,  $g(k) \geq 2^k + \lfloor (3/2)^k \rfloor -2$, simply because the integer $n=2^k \lfloor (3/2)^k \rfloor -1$ can only be written as a sum of $\lfloor (3/2)^k \rfloor -1$ terms $2^k$ and   $2^k-1$ terms $1^k$.\footnote{This example was apparently already known in the 1770's to J. A. Euler, son of L. Euler; see \cite[Item 36 p. 203-204]{Eul62}.}  This is conjecturally the lower bound as well (see e.g. \cite[Ch. 1]{Vau97} for progress),
but the Hardy and Littlewood circle method \cite{HarLit20} made it more satisfactory to study $G(k)$, the least $s$ such that every \emph{sufficiently large} $N$ may be represented by (\ref{N_sum}), with $x_i \geq 1$.

We will focus here on the related quantity $\tilde{G}(k)$, traditionally denoting the least $s$ for which the asymptotic (\ref{rskN_asymp}) for $r_{s,k}(N)$ is known to hold, for all sufficiently large $N$. Since the work of Hardy and Littlewood in the 1920's, the quest  to reduce $\tilde{G}(k)$ as far as possible has fueled many technical innovations in the circle method, the most powerful analytic method for counting integral solutions to Diophantine equations.

Certainly, we must have $\tilde{G}(k) \geq k+1$, since purely geometric considerations show that when $s=k$, the number of  integer tuples $0<x_1 \leq \cdots \leq x_k$ with  $x_1^k + \cdots + x_k^k \leq X$ is asymptotically $\ga X$ as $X \maps \infty$ for a value of $\ga$ strictly smaller than 1 (see e.g. \cite[Ch. 9]{Dav05}). The case $k=2$ of sums of squares is somewhat special since it is closely connected with theta functions; sums of at least $5$ squares were treated by Hardy and Littlewood \cite{HarLit20}, four squares by Kloosterman \cite{Klo49}, and three squares  by Duke \cite{Duk88} and Iwaniec \cite{Iwa87} (via modular forms of half-integral weight; see \cite[Ch. 20]{IK} for an overview of the quadratic case).
 For $k \geq 3$, it can be hoped that the circle method will ultimately succeed for $s \geq k+1$, and certainly one may see from heuristics  that at least $s \geq 2k+1$ is a reasonable goal for the method (see e.g. \cite{HB07b}).

In broad strokes, the idea of the circle method is to express 
\[ r_{s,k}(N) = \int_0^1 g_k(\al;X)^s e(-N\al) d\al ,\]
where we choose $X \approx N^{1/k}$ and define
\[g_k(\al;X)=\sum_{1 \leq x \leq X} e(\al x^k).\]
One then dissects the integral over $(0,1]$ into a portion called the major arcs $\mathfrak{M}$, comprised of small disjoint intervals centered at rational numbers $a/q$ with $q$ sufficiently small (determined relative to $N,k$), and the remainder, called the minor arcs $\mathfrak{m}$. The major arcs contribute the main term in the asymptotic (\ref{rskN_asymp}) for $r_{s,k}(N)$, with singular series
\[ \mathfrak{S}_{s,k}(N) = \sum_{q=1}^\infty \sum_{\bstack{a=1}{(a,q)=1}}^q ( q^{-1} \sum_{m=1}^q e(am^k/q))^s e(-Na/q).\]
 For $s \geq \max\{5,k+2\}$, $\mathfrak{S}_{s,k}(N)$ converges absolutely to a non-negative real number $\ll 1$, and for $N$ satisfying reasonable congruence conditions, $\mathfrak{S}_{s,k}(N) \gg 1$, so that the main term in (\ref{rskN_asymp}) is non-zero, and genuinely larger than the error term, for $N$ sufficiently large. (See \cite[\S 4.3, 4.5, 4.6]{Vau97} for precise conclusions.)

Thus the main challenge in proving the asymptotic for $r_{s,k}(N)$ is showing that the minor arcs contribute exclusively to the smaller remainder term; we will summarize a few key points in the historical development with a focus on large $k$. (Detailed sources on the circle method include \cite{Vau97}, \cite[Ch. 20]{IK}, \cite{Dav05}, and for Waring's problem, \cite{VauWoo02}, \cite{Woo12b}.)
 Hardy and Littlewood ultimately showed in \cite{HarLit22} that $\tilde{G}(k) \leq (k-2)2^{k-1}+5$, using the bound for $f_k(\albf;X)$ (and hence for $g_k(\al;X)$) due to Weyl, of the form
\beq\label{Weyl_bound}
 |f_k(\albf;X)| \ll_{k,\ep} X^{1+\ep} (q^{-1} + X^{-1}  + qX^{-k})^{\frac{1}{2^{k-1}}},
 \eeq
which holds whenever the leading coefficient $\al_k$ has rational approximation $|\al_k - a/q| \leq q^{-2}$ with $(a,q)=1$. 
(This improves on the trivial bound $X$ when, say, $X \leq q \leq X^{k-1}$, which is roughly the situation of the minor arcs.)

Hua \cite{Hua38} improved Hardy and Littlewood's strategy and obtained $\tilde{G}(k) \leq 2^k +1$ by using what is now called Hua's inequality, namely that for each $1 \leq \ell \leq k$,
\beq\label{Hua_ineq}
 \int_{0}^1 |g_k(\al;X)|^{2^\ell} d\al \ll_{\ell,k,\ep} X^{2^\ell - \ell + \ep}.
 \eeq
Next, Vinogradov's Mean Value Method came on the scene, with the observation that
\[ \int_0^1 |g_k(\al;X)|^{2s} d\al \]
is the number of solutions of the single equation
\[ x_1^k + \cdots + x_s^k = x_{s+1}^k  + \cdots + x_{2s}^k \]
with $1 \leq x_i \leq X$, so that 
\[ \int_0^1 |g_k(\al;X)|^{2s} d\al  = \sum_{|h_1| \leq sX} \cdots \sum_{|h_{k-1}| \leq sX^{k-1}} \int_{(0,1]^{k}} |f_k(\albf;X)|^{2s} e(-h_1 \al_1 - \cdots -h_{k-1} \al_{k-1} )d\albf.\]
Hence by the triangle inequality, 
\[ \int_0^1 |g_k(\al;X)|^{2s} d\al \ll_{s,k} X^{\frac{1}{2} k(k-1) } J_{s,k}(X).\]
From this, Vinogradov \cite{Vin35} deduced that $\tilde{G}(k) \leq (C+o(1))k^2 \log k$ for $C=10$, ultimately refined to Hua's statement \cite{Hua49} that $C=4$ suffices, a record that held until the leading constant was refined to $2$ by Wooley \cite{Woo92} and to 1 by Ford \cite{For95}.

Wooley's breakthrough with efficient congruencing removed the logarithmic factor, and then over many papers he successively improved $\tilde{G}(k)$ for all $k \geq 5$. By the time decoupling entered the picture, Wooley \cite[Thm. 12.2]{Woo17b} had established, for large $k$, the record $\tilde{G}(k) < (1.540789... +o(1)) k^2 +O(k^{5/3})$.

\subsubsection{Improvements on the Weyl bound}
Now, with the Vinogradov Mean Value Theorem in hand, more can be said. There are two types of improvements that can affect the circle method: ``pointwise''  bounds for  individual Weyl sums, and  mean value estimates. 
To improve the classical Weyl bound (\ref{Weyl_bound}), the goal is to show that if $k \geq 3$ and for some $2 \leq j \leq k$ one has $|\al_j - a/q| \leq q^{-2}$ with $(a,q)=1$ then 
\beq\label{Weyl_bound'}
 |f_k(\albf;X)| \ll_{k,\ep} X^{1+\ep} (q^{-1} + X^{-1}  + qX^{-j})^{\sig(k)},
\eeq
with $\sig(k)$ bigger than $2^{-(k-1)}$.
As Vinogradov knew, a bound for an individual sum can be extracted from a mean value estimate, since  under the aforementioned hypotheses,
\[
|f_k(\albf;X)| \ll ( X^{\frac{1}{2}k(k-1)} J_{s,k-1}(2X) (q^{-1} + X^{-1} + qX^{-j}))^{\frac{1}{2s}} \log (2X).
\]
(See e.g. \cite[Thm. 5.2]{Vau97}, or a heuristic explanation in \cite[\S 1]{Woo14b}.)
 Improvements on Weyl's bound were already known from \cite{HB88} and \cite{RobSar00}; after successive improvements with efficient congruencing, Wooley \cite[Thm. 7.3]{Woo14b} showed one can take $\sig(k)^{-1} =  2(k-1)(k-2)$. Now, as an application of the Vinogradov Mean Value Theorem, Bourgain has confirmed that one can take $\sig(k)^{-1} = k(k-1)$  \cite[Thm. 5]{Bou16x}. (This is as expected, see e.g. \cite[Ch. 4 Thm. 4]{Mon94} and the subsequent remark.)
Anticipating the resolution of the Main Conjecture, Wooley \cite[Thm. 4.1]{Woo12b} had already deduced that with this improvement on the Weyl bound in hand, 
\beq\label{Woo_Gk}
 \tilde{G}(k) \leq k^2 + 1 - \left[ \frac{\log k }{\log 2 } \right], \qquad  \text{for all $k \geq 3$.}
 \eeq

\subsubsection{Improvements on Hua's inequality}\label{sec_Hua}
Regarding the second type of input to the circle method, namely mean value estimates, Bourgain \cite[Thm. 10]{Bou16x} has now improved Hua's inequality (\ref{Hua_ineq}), showing that for each $1 \leq \ell \leq k$,
\beq\label{Hua_ineq_Bou}
\int_0^1 |g_k(\al;X)|^{\ell(\ell+1)}d\al \ll X^{\ell^2+\ep} .
\eeq
This is not a consequence  of the Vinogradov Mean Value Theorem itself; Bourgain instead deduces it from an $\ell^2$ decoupling result (for $L^{\ell(\ell+1)}$) for a curve in $\R^\ell$:
\beq\label{Gamma_Hua}
 \Gamma_{\ell,k} = \{ (t, t^2, \ldots, t^{\ell-1},t^k), t \in [1,2] \}.
 \eeq
(See \S \ref{sec_an_heur} for a remark on this decoupling result.)
With (\ref{Hua_ineq_Bou}) in place of certain estimates in \cite{Woo12b}, Bourgain \cite[Thm. 11]{Bou16x} has established the current record in Waring's problem, which improves on (\ref{Woo_Gk}) for all $k \geq 12$;  for large $k$ it can be written as
\[ \tilde{G}(k) \leq k^2 -k +O(\sqrt{k}).\]

\subsection{The Riemann zeta function}
As is well-known, our ability to prove good error terms in the asymptotic $\pi(x) \sim x/\log x$ for the prime counting function depends on our understanding of the Riemann zeta function $\zeta(s)$ within the critical strip $0 \leq \Re(s) \leq 1$.
This in turn rests heavily on estimating exponential sums, since for $s=\sig+it$ near the line $\sig=1$ and  with $t \geq 3$,
\[ \zeta( s) = \sum_{n \leq t} n^{-\sig - it} + O(t^{-\sig}).\]
After partial summation, the key is to understand sums of the form 
\[ \sum_{N < n \leq M} n^{it} = \sum_{N<n \leq M} e^{i t \log n},\]
which via a Taylor expansion can be made to resemble an exponential sum of a real-valued polynomial.
Vinogradov's Mean Value Method, in the form of the classical result (\ref{Jsk_classical}), 
resulted in the so-called Vinogradov-Korobov zero-free region: $\zeta(\sig +it) \neq 0$ for
\[ \sig \geq 1 - \frac{C}{(\log t)^{2/3} (\log \log t)^{1/3}}, \qquad t \geq 3,\]
for a certain absolute constant $C>0$ (see e.g. \cite[Ch. 6]{Ivi85}).
 This shape of result remains the best-known today, with Ford \cite{For02} providing an explicit value of $C$. 
 
 Interestingly, the current methods of resolving the Main Conjecture  do not improve the zero-free region for $\zeta(s)$, since any such improvement would depend delicately upon the implicit constant in Theorem \ref{thm_VMVT}, with its dependence on $k$ and $\ep$ dictating the upper bound implied for the growth of $\zeta(s)$ near the line $\Re(s)=1$.
However, Heath-Brown has applied the Vinogradov Mean Value Theorem to deduce a new $k$-th derivative estimate for Weyl sums \cite{HB16x}, an improvement on the classical van der Corput $k$-th derivative estimate. Heath-Brown's  work provides new upper bounds for $\zeta(s)$ throughout the critical strip, zero-density estimates, and certain moment estimates for $\sig$ sufficiently close to 1. Here we specify only the first of his results.
The so-called convexity bound for $|t| \geq 1$ and $0 \leq \sig \leq 1$ shows that
\beq\label{zeta_convexity}
 \zeta(\sig+it) \ll_\ep |t|^{\frac{1}{2}(1-\sig) + \ep}, 
 \eeq
for all $\ep>0$.
One consequence of Heath-Brown's new work is that uniformly in $|t| \geq 1$ and $0 \leq \sig \leq 1$, 
  \[ \zeta(\sig + it) \ll_\ep |t|^{\frac{1}{2}(1-\sig)^{3/2}+\ep},\]
  for all $\ep>0$.
  This improves the coefficient of $(1-\sig)^{3/2}$ from $4.45$ in a previous record of Ford \cite{For02}; for further results and refinements, see \cite{HB16x}.

Other instances of decoupling have impacted our knowledge of the Riemann zeta function in two further ways, not directly via the Vinogradov Mean Value Theorem.  On the critical line, the convexity estimate (\ref{zeta_convexity}) shows that $\zeta(1/2+it)\ll_\ep |t|^{1/4+\ep}$,  whereas the Lindel\"{o}f Hypothesis (implied by the Riemann Hypothesis) posits that $\zeta(1/2+it)\ll_\ep |t|^\ep$,  as $|t| \maps \infty$, for every $\ep>0$.  One method to improve on the convexity estimate is the Bombieri-Iwaniec method for bounding exponential sums, which, as codified by Huxley \cite{Hux96}, includes a so-called first and second spacing problem. Bourgain used a specific instance of a decoupling theorem for a nondegenerate curve in $\R^4$ (also deduced in \cite{BouDem16a} from a decoupling theorem for nondegenerate surfaces in $\R^4$) to make improvements on the first spacing problem, with two consequences.
 First, Bourgain \cite{Bou17a} proved   the current best-known upper bound on the critical line,
 \[ |\zeta(\frac{1}{2} + it) | \ll |t|^{\frac{13}{84}+\ep} = |t|^{ 0.1548...+\ep},\]
improving on the previous record exponent $32/205 = 0.1561...$ of Huxley \cite{Hux05}.
 Furthermore, Bourgain and Watt \cite{BouWat15x} have set new records pertaining to second moments of the zeta function.

\subsection{Further consequences}

Considerations related to the Vinogradov Mean Value Method have  consequences for many other arithmetic problems, such as Tarry's problem and the distribution of polynomial sequences modulo 1 (e.g. \cite{Woo12a}), and solutions to congruences in short intervals \cite{CCGHSZ14}.
It is also natural to consider multi-dimensional versions of the Vinogradov system (\ref{Vin_sys_dfn}), and via multi-dimensional efficient congruencing, Parsell, Prendiville and Wooley proved nearly optimal results related to the analogous Main Conjecture in broad generality \cite{PPW13}. The relevant Main Conjecture for Vinogradov systems of degree $k$ in $n$ dimensions has now been resolved for $k=2=n$ by Bourgain and Demeter \cite{BouDem15x} and for $k=3,n=2$ by Bourgain, Demeter, and Guo \cite{BDG16x} (see also \cite{BouDem16b}), via a sharp decoupling result for a two-dimensional surface in $\R^9$; this has also been generalized in \cite{Guo17x}. Both in the setting of the original Vinogradov system and its higher-dimensional variants, the Vinogradov Mean Value Theorem has implications for Burgess-type bounds for short mixed character sums \cite{HBP15}, \cite{Pie16}.

Wooley's efficient congruencing methods can be applied to Vinogradov-type problems in number field and function field settings (as already remarked in \cite{Woo12a}), and further results are anticipated in such settings, see e.g. \cite{KLZ14}. 
Decoupling methods presently include the flexibility to consider non-integral moments and real (not necessarily integral) solutions to Diophantine equations and inequalities (see \S \ref{sec_dis_Vin}). In addition to the arithmetic implications of decoupling for certain carefully chosen curves and surfaces, there are implications for analytically-motivated problems, to which we will return in \S \ref{sec_decoupling_other}. 
The strengths of both the efficient congruencing and decoupling methods will no doubt be further developed.

\section{Classical considerations for the Main Conjecture}\label{sec_MC_classical}
In this section, we record the heuristics that motivate the Main Conjecture, discuss the critical index, the trivial cases $k=1,2$, the derivation of an asymptotic for $J_{s,k}(X)$ from an upper bound, and the translation-dilation invariance of the underlying system of equations. 

\subsection{Heuristics motivating the Main Conjecture}\label{sec_MC_Heuristics}
When we consider solutions $1 \leq x_i \leq X$ to the Vinogradov system
\beq\label{Vin_sys_dfn2}
x_1^j + \cdots + x_s^j = x_{s+1}^j + \cdots + x_{2s}^j, \qquad 1 \leq j \leq k,
\eeq
we see immediately that a trivial upper bound is $J_{s,k}(X) \leq X^{2s}$. On the other hand, the diagonal solutions provide a lower bound:
\beq\label{Jsk_lowerbd1}
 J_{s,k}(X) \geq s! X^s  + O_s(X^{s-1}). 
 \eeq
An alternative lower bound may be obtained from the observation that for $1 \leq x_i \leq X$,
\[
|(x_1^j -x_{s+1}^j) + \cdots + (x_s^j - x_{2s}^j )| \leq sX^j , \qquad 1 \leq j \leq k
\]
so that
\[ X^{2s} \ll \sum_{|h_1| \leq sX} \cdots \sum_{|h_k| \leq sX^k} \int_{(0,1]^k}  |f_k(\albf;X)|^{2s} e(\albf \cdot \hbf) d\albf .
\]
By the triangle inequality, the right-hand side is at most $\ll X \cdot X^{2} \cdots X^k J_{s,k}(X)$, so that we may extract the lower bound 
\beq\label{Jsk_lowerbd2}
 J_{s,k}(X) \gg_{s,k} X^{2s - \frac{1}{2} k(k+1)}.
 \eeq
 Alternatively, we may think of this as the contribution to $J_{s,k}(X)$ in (\ref{Jsk_intro_int}) from $\albf$ near the origin, 
e.g. with $|\al_j| \leq \frac{1}{8k} X^{-j}$ for $1 \leq j \leq k$. The set of such $\albf$ has measure $(1/4k)^kX^{-\frac{1}{2}k(k+1)}$, and for each such $\albf$, 
\[|f_k(\albf;X)| \geq |\Re ( f_k(\albf;X))| \geq  X \cos(\pi /4),\]
leading to a contribution to (\ref{Jsk_intro_int}) of at least 
$(1/ \sqrt{2})(1/4k)^k X^{2s-\frac{1}{2}k(k+1)}$.
The Main Conjecture posits that (\ref{Jsk_lowerbd1}) and (\ref{Jsk_lowerbd2}) represent the true order of $J_{s,k}(X)$.

\subsection{The critical case}\label{sec_VMVT_critical_case}
If one has obtained for a certain $k$ the estimate at the critical index $s_k = \frac{1}{2}k(k+1)$,
\beq\label{critical_exponent_Jsk}
\int_{(0,1]^k} |f_k(\albf;X)|^{k(k+1)} d\albf \ll_{k,\ep} X^{\frac{k(k+1)}{2} +\ep}, 
\eeq
one may immediately deduce the Main Conjecture holds for all $s \geq 1$, for this $k$. For supposing $s > s_k$, we may write
\[
J_{s,k}(X)	\leq \sup_{\albf \in (0,1]^k} |f_k(\albf;X)|^{2s - k(k+1)}  \int_{(0,1]^k}  |f_k(\albf;X)|^{k(k+1)} d\albf,
\]
	whereupon applying the trivial bound $ |f_k(\albf;X)| \leq X$ in the first term and the critical case (\ref{critical_exponent_Jsk}) to the integral, 
one obtains the bound $X^{2s-k(k+1)/2+ \ep}$.
Considering $s <s_k$, one may apply H\"{o}lder's inequality  with $q=s_k/s>1$ and its conjugate exponent $q'$:
\[
\int_{(0,1]^k}  |f_k(\albf;X)|^{2s} d\albf
	\leq ( \int_{(0,1]^k} 1^{q'})^{1/q'} \left( \int_{(0,1]^k} |f_k(\albf;X)|^{2sq} d\albf \right)^{1/q}.
\]
Again applying the critical case (\ref{critical_exponent_Jsk}) to the last factor, 
this is $\ll_{k,\ep} (X^{s_k+\ep})^{1/q} \ll_{k,\ep} X^{s+\ep}$, as desired.

\subsection{Trivial cases $k=1$ and $k=2$}\label{sec_k12}
The Main Conjecture holds trivially for $k=1$ and $s_k=1.$
The case of $k=2$ and $s_2=3$ may also be handled in an elementary fashion. After rearranging, the claim is that there exist at most $X^{3+\ep}$ solutions to the system
\begin{eqnarray*}
x_1 + x_2 - y_3 &=& y_1 + y_2  - x_3 \\
x_1^2 + x_2^2 -y_3^2 & = & y_1^2 + y_2^2- x_3^2 .
\end{eqnarray*}
By subtracting the second equation from the square of the first equation and applying the identity $(a+b-c)^2  - (a^2+b^2  -c^2) =2(a-c)(b-c)$, we  obtain 
\beq\label{k2_elementary_id}
 (x_1 - y_3)(x_2 - y_3)  = (y_1 - x_3)(y_2 -x_3).
 \eeq
In cases where both sides of this identity are nonzero, we may fix $x_1,x_2,y_3$  from $O(X^3)$ possible choices, and for each such triple there are $\ll d(X^2) \ll X^\ep$ choices for the factors on the right-hand side (where $d(n)$ denotes the number of divisors of $n$). Then via the linear equation, written as
\[ x_1 + x_2 - y_3 = (y_1 - x_3) + (y_2 - x_3) + x_3 ,\]
we see that once we have chosen $(y_1 - x_3)$ and $(y_2-x_3)$, the choice for $x_3$ is unique. Then $y_1,y_2$ are also uniquely determined. In cases where the identity (\ref{k2_elementary_id}) vanishes, assuming without loss of generality $x_1=y_3$ and $y_1=x_3$, we may choose $x_1,x_2,x_3$ freely and then $y_3$ and $y_1$ are uniquely determined, while $y_2$ is uniquely determined  by the linear equation, showing that $J_{3,2}(X) \ll X^{3+\ep}$. 
More sophisticated considerations  show that $J_{3,2}(X) \sim \frac{18}{\pi^2} X^3 \log X$, with explicit lower order terms (see \cite{Rog86}, \cite{BloBru10});  in particular this verifies that the $X^\ep$ term cannot be omitted from the Main Conjecture for $k=2$.

\subsection{The asymptotic for $J_{s,k}(X)$}\label{sec_ep_removal}
Recall the statement following Theorem \ref{thm_VMVT} that once the upper bound (with a factor $X^\ep$) in the Main Conjecture  has been verified for $s=\frac{1}{2}k(k+1)$, an asymptotic (without a factor $X^\ep$) can be obtained for $J_{s,k}(X)$, as long as $k \geq 3$ and the number of variables satisfies $s > \frac{1}{2} k(k+1)$. Indeed historically each improved upper bound toward the Main Conjecture has enabled the proof of an asymptotic for $J_{s,k}(X)$, possibly in slightly more variables; see e.g. \cite[\S 7.3]{Vau97} for the deduction of the asymptotic (\ref{Jsk_asymp}) in the Vinogradov range, or  Wooley's proof of an asymptotic for $s \geq k(k+1)+1$ once he proved the Main Conjecture  for $s \geq k(k+1)$ \cite[\S 9]{Woo12a}.

The strategy for deducing an asymptotic is related to the Prouhet-Escott-Tarry and Hilbert-Kamke problems (see e.g. \cite[\S X.3]{Hua65},\cite[Ch. 3]{ACK04},\cite{Woo96}), and is a nice application of the circle method.
Roughly speaking, the idea is to dissect the integral expression (\ref{Jsk_intro_int}) for $J_{s,k}(X)$ as
\beq\label{J_int}
 J_{s,k}(X)  = \int_{(0,1]^k}  |f(\albf)|^{2s} d\albf = \int_{\mathfrak{M}}  |f(\albf)|^{2s} d\albf+ \int_{\mathfrak{m}}  |f(\albf)|^{2s} d\albf
 \eeq
with $f(\albf)  = f_k(\albf;X) = \sum_{1 \leq x \leq X} e(\al_1 x + \cdots + \al_k x^k)$, according to 
 an appropriate choice of major arcs $\mathfrak{M}$ and minor arcs $\mathfrak{m}$.
Once the Main Conjecture is known for the critical index $s_k$, the major arcs can be chosen to be very small, since the only property required of the minor arcs is that there exists some $\tau = \tau(k)>0$ such that 
 \beq\label{Wooley1}
 \sup_{\albf \in \mathfrak{m}} |f(\albf)| \ll_{k,\ep} X^{1 - \tau + \ep} ,
 \eeq
 for every $\ep>0$.
 With any such $\tau$ in hand, supposing $s > \frac{1}{2}k(k+1)$, we may apply the known result of the  Main Conjecture (including the $X^\ep$ factor)
  for $s_k = \frac{1}{2}k(k+1)$ to see that
 \[ \int_{\mathfrak{m}} |f(\albf)|^{2s} d\albf \ll (\sup_{\albf \in \mathfrak{m}} |f(\albf)|)^{2s -k(k+1)} \int_{\mathfrak{m}} |f(\albf)|^{k(k+1)} d\albf
	\ll_{k,\ep} (X^{1 - \tau +\ep})^{2s - k(k+1)} X^{\frac{1}{2}k(k+1)+\ep},
\]
for every $\ep>0$.
Since $\tau$ is a fixed positive number, and $s>\frac{1}{2}k(k+1)$, this is now  $\ll_{s,k} X^{2s - \frac{1}{2}k(k+1) - \tau'}$ for some $\tau'  = \tau'(s,k)>0$.

Defining the major and minor arcs as in \cite[Eqn. (7.1), Lemma 7.1]{Woo17a}, Wooley obtains $\tau= (4k^2)^{-1}$ on the minor arcs, and shows that on the major arcs, as soon as $k \geq 3$ and $2s> \frac{1}{2}k(k+1)+2$ (a far weaker condition than we assume), 
\[
 \int_{\mathfrak{M}} |f(\albf)|^{2s} d\albf  \ll X^{2s - \frac{1}{2}k(k+1)},
\]
thus completing the removal of the $X^\ep$ factor from both major and minor arcs as soon as $s > \frac{1}{2}k(k+1)$ (as noted in \cite[\S 5]{BDG16}). 

To obtain the asymptotic formula, one can adapt e.g. the argument of \cite[\S 9]{Woo12a} combined with \cite[\S 7.3]{Vau97} to
show that on the major arcs as chosen above (or chosen to be even smaller, given how little a savings is needed on the minor arcs), for $k \geq 3$ and integers $s > \frac{1}{2}k(k+1)$,
\beq\label{Wooley2}
 \int_{\mathfrak{M}} |f(\albf)|^{2s} d\albf  = \mathfrak{S}(s,k)\mathfrak{I}(s,k)  X^{2s - \frac{1}{2}k(k+1) } + O(X^{2s - \frac{1}{2}k(k+1) - \tau''}),
 \eeq
 for some $\tau'' = \tau''(s,k)>0$,
 with singular series
\[
\mathfrak{S}(s,k) = \sum_{q=1}^\infty \sum_{a_1=1}^q \cdots \sum_{\bstack{a_k=1}{(a_1,\ldots, a_k,q)=1}}^q
	| q^{-1} \sum_{r=1}^q e((a_1r + \cdots + a_k r^k)/q)|^{2s} \]
and singular integral
\[
\mathfrak{I}(s,k) = \int_{\R^k} | \int_0^1 e(\be_1 \ga + \cdots + \be_k \ga^k ) d\ga |^{2s} d\bebf .
\]
Here the singular series converges absolutely in the more flexible range $2s > \frac{1}{2}k(k+1) + 2$, and the singular integral if $2s> \frac{1}{2}k(k+1) + 1$ \cite[Thm. 3.7]{ACK04}.
 As Vaughan has pointed out \cite[\S 7.3]{Vau97}, the positivity of the singular integral and singular series for $s> \frac{1}{2}k(k+1)$ is guaranteed by the known lower bound (\ref{Jsk_lowerbd2}) for $J_{s,k}(X)$. Similarly, the lower bound (\ref{Jsk_lowerbd1}) indicates that the asymptotic (\ref{Jsk_asymp}) would not be expected to hold for $s < \frac{1}{2}k(k+1)$.
(In the setting of (non-translation-dilation invariant) discrete restriction problems, such ``$\ep$-removal'' strategies have also been successful in work of  Bourgain \cite{Bou89b,Bou93}, Wooley \cite{Woo17a}, and elsewhere.)

\subsection{Translation-dilation invariance}\label{sec_trans_dil}

The Vinogradov system is said to be translation-dilation invariant (or affine invariant): given real numbers $\lam,\xi$ with $\lam  \neq 0$, then
$x_1,\ldots, x_{2s}$  are solutions to the system 
\beq\label{Vin_sys_xy}
\sum_{i=1}^s x_i^j = \sum_{i=1}^s x_{s+i}^j  \qquad 1 \leq j \leq k
\eeq
if and only if they are solutions to the system
\[ \sum_{i=1}^s (\lam x_i + \xi)^j = \sum_{i=1}^s (\lam x_{s+i} + \xi)^j  \qquad 1 \leq j \leq k.\]
The invariance under $\lam$ is perhaps the easiest to see, as each individual equation is homogeneous of degree $j$.
 Now with $\lam=1$, to understand the invariance under $\xi \neq 0$, an upward induction in $j$ combined with the binomial theorem shows that for each $j$,
\[ \sum_{i=1}^s ( (x_i -\xi +\xi)^j - (x_{s+i} - \xi+\xi)^j) = \sum_{\ell=1}^j \binom{j}{\ell} \xi^{j-\ell} \sum_{i=1}^s ( (x_i-\xi)^\ell - (x_{s+i}-\xi)^\ell),\]
so that for $\xi \neq 0$ and a particular choice of the $x_i$, the left-hand side vanishes if and only if (after the upward induction) the $j$-th term on the right-hand side vanishes.

A consequence of translation-dilation invariance is the following crucial observation: counting solutions $x_1, \ldots, x_{2s}$ to (\ref{Vin_sys_xy}) such that $x_i \leq X$ and $x_i \con \xi \modd{q}$ for some fixed residue class $\xi$ and $q \leq X$ is essentially equivalent to counting solutions $z_1, \ldots, z_{2s}$ to (\ref{Vin_sys_xy}) with no congruence condition, but in the smaller range $z_i \leq X/q$. This is made formal e.g. in Lemma \ref{k2_EC_lemma1} below.

\section{Wooley's efficient congruencing method}\label{sec_EC}
We now turn to Wooley's powerful efficient congruencing approach to the Vinogradov Mean Value Method. After defining some initial notation, we sketch the classical approach from the era of Vinogradov and Linnik. We then indicate via the simple exemplar  case of $k=2$ how Wooley's approach explosively magnifies any purported extraneous solutions to the Vinogradov system that might cause $J_{s,k}(X)$ to exceed the conjectured bound, thus arriving at a contradiction, unless $J_{s,k}(X)$ does indeed satisfy the bound. 
Finally, we sketch the method in the nontrivial case of $k=3$. 

\subsection{Initial set-up}
We first prepare some notation, which we will employ both for the classical and the efficient congruencing arguments. We define
\[ f(\albf) = \sum_{1 \leq x \leq X} e(\al_1 x + \cdots + \al_k x^k) \]
and the new notation
\[ f_a(\albf;\xi) = \sum_{\bstack{1 \leq x \leq X}{x \con \xi \modd{p^a}}} e(\al_1 x + \cdots + \al_k x^k) ,\]
where $p>k$ is a prime that will be fixed once and for all, $\xi$ is any integer, and $a \geq 1$ is an integer. We will use the latter function to count solutions to the Vinogradov system 
\beq\label{Vin_sys}
 x_1^j + \ldots + x_s^j = x_{s+1}^j + \cdots +x_{2s}^j, \qquad 1 \leq j \leq k
 \eeq
with $x_i$ satisfying certain congruence conditions. Precisely, if we want to
count solutions to the Vinogradov system in which 
 \[ x_i \con \xi \modd{p^a} \qquad \text{for $1 \leq i \leq m$ and $s+1 \leq i \leq s+m$}\]
 and 
  \[ x_i \con \eta \modd{p^b} \qquad \text{for $m+1 \leq i \leq s$ and $m+s+1 \leq i \leq 2s$},\]
we may also express this count as the bilinear integral
  \beq\label{EC_bilinear}
   I_m(X;\xi,\eta;a,b) := \int_{(0,1]^k} |f_a(\albf;\xi)|^{2m} |f_b(\albf;\eta)|^{2(s-m)} d\albf.
   \eeq
We define $I_m(X;\xi,\eta;a,b)$  for each $0 \leq m \leq s-1$; in the event that $m=0$, this expression is in fact independent of $\xi$ and $a$.
 We will later take a maximum over the possible residues $\xi$ and $\eta$ modulo $p$, and so we define 
 \[ I_0(X;a,b) = \max_{\eta \modd{p}} I_0(X;\xi,\eta; a,b) \]
 (in which we note that $\xi,a$ are irrelevant),
 and for $1 \leq m \leq s-1$ set 
 \beq\label{Im_dfn}
  I_m(X;a,b) = \max_{\xi \not\con \eta \modd{p}} I_m(X;\xi,\eta; a,b). 
  \eeq
 
\subsection{The classical $p$-adic method}

We will not prove Vinogradov's result (\ref{Jsk_classical}) rigorously, but instead summarize the two points in the classical argument that are most relevant for the comparison to efficient congruencing. 
We assume that $p \leq X$ is a fixed prime, to be chosen momentarily.
First, we claim that for any $k' \geq 1,$ 
\beq\label{Jsk_relation1}
 J_{k+k',k}(X) \leq p^{2k'} I_k(X;0,1).
\eeq
(Here we understand that when applying (\ref{Im_dfn}) in (\ref{Jsk_relation1}) we used $m=k$, $s=k+k'$.) This means we have passed to counting solutions to the Vinogradov system that satisfy a congruence restriction modulo $p$.
Second, if $p^k > X$, we claim that for any $k'\geq 1$, the main contribution to $I_k(X;0,1)$ is of size at most 
\beq\label{Jsk_relation2}
 k! p^{\frac{1}{2}k(k-1)} X^k J_{k',k}(X/p).
 \eeq
 With this, we have passed back to counting solutions to the Vinogradov system without a congruence restriction, but now in fewer variables and in a smaller range.

We will indicate the reasoning behind these relations in a moment, but first let us see how they contribute to Vinogradov's original bound (\ref{Jsk_classical}). We  aim to prove that  if we define $\lam_{s,k}^*$ to be the infimum of those exponents $\lam_{s,k}$ for which $J_{s,k}(X) \ll_{s,k,\ep} X^{\lam_{s,k}+\ep}$ for all $X \geq 1$ and all $\ep>0$,  and define $\eta_{s,k}$ by 
\beq\label{eta_dfn}
 \lam_{s,k}^* = 2s - \frac{1}{2}k(k+1) + \eta_{s,k},
 \eeq
then the remainder $\eta_{s,k}$ is small, at least if $s$ is quite large. (The Main Conjecture  for $J_{s,k}(X)$ claims that $\eta_{s,k}=0$.)

The strategy is to bound $\eta_{s,k}$ recursively, by controlling $\eta_{s+k,k}$ in terms of $\eta_{s,k}$. We first apply the relation (\ref{Jsk_relation1}) with $k'=s$, followed by the claim that the main contribution is from (\ref{Jsk_relation2}), so that  (roughly, ignoring other contributions),
\[ J_{s+k,k}(X) \ll p^{2s} I_k(X;0,1) \ll_k p^{2s + \frac{1}{2} k(k-1)} X^k J_{s,k}(X/p)
	\ll_{s,k,\ep} p^{2s+\frac{1}{2}k(k-1)- \lam_{s,k}^*} X^{k+ \lam_{s,k}^*+\ep}.\]
It now makes sense to choose $p$, and the hypothesis underlying (\ref{Jsk_relation2}) motivates us to choose $p>k$ to be a prime 
\[ X^{1/k} < p \leq 2X^{1/k},\]
which exists by Bertrand's postulate (or the prime number theorem), if say $X \geq k^k$.
Then we see from our upper bound for $J_{s+k,k}(X)$ and the definition (\ref{eta_dfn}) that 
we have shown $\lam_{s+k,k}^* \leq 2(s+k) - \frac{1}{2}k(k+1) + \eta_{s+k,k}$ with $\eta_{s+k,k} \leq \eta_{s,k}(1-1/k)$. This, combined with the recollection that $J_{k,k} \leq k! X^{k}$, allows induction to prove that 
\[\eta_{s,k} \leq \frac{1}{2}k^2 (1-1/k)^{[s/k]},\]
 the outcome of Vinogradov's original program.

\subsubsection{Underlying philosophy of the key relations}
The relation (\ref{Jsk_relation1}) claims that solutions to the Vinogradov system  may be bounded by counting solutions with a congruence condition. We see this as follows: for any $\ell \geq 1$, by H\"{o}lder's inequality,
\[ 
|f(\albf)|^{2\ell} = \left| \sum_{\xi=1}^p f_1(\albf;\xi) \right|^{2\ell} \leq p^{2\ell-1} \sum_{\xi=1}^p |f_1(\albf;\xi)|^{2\ell}.
\]
Now applying this with $\ell =k'$ we have 
\beq\label{Jkk'}
 J_{k+k',k}(X) \leq p^{2k'} \max_{1 \leq \xi \leq p} \int_{(0,1]^k} |f(\albf)|^{2k} |f_1(\albf;\xi)|^{2k'} d\albf,
 \eeq
and the claim follows.

The extraction of the main contribution (\ref{Jsk_relation2}) is the heart of the proof, and most important for the later comparison to efficient congruencing. For any fixed $\xi$, the integral on the right-hand side of (\ref{Jkk'}) is counting solutions $(x_1,\ldots, x_{2s})$ to the Vinogradov system with $x_i$ having no congruence restrictions for $1 \leq i \leq k$ and $s+1 \leq i \leq s+k$, and the restriction $x_i \con \xi \modd{p}$ for all the other $2k'$ variables.
Re-naming variables, this is counting solutions to the system 
\[ 
\sum_{i=1}^k (x_i^j  - y_i^j)  = \sum_{\ell=1}^{k'} ( ( pu_\ell + \xi)^j - (pv_\ell + \xi)^j), \qquad 1 \leq j \leq k,
\]
where $1 \leq x_i,y_i \leq X$ and $(1-\xi)/p \leq u_\ell, v_\ell \leq (X-\xi)/p$.
By translation-dilation invariance, this is equivalent to counting solutions to 
\beq\label{xyuv}
\sum_{i=1}^k ((x_i - \xi)^j  - (y_i - \xi)^j)  = p^j \sum_{\ell=1}^{k'} ( u_\ell^j - v_\ell^j), \qquad 1 \leq j \leq k.
\eeq
This reveals a new system of not equations but \emph{congruences},
\beq\label{xyp} 
\sum_{i=1}^k (x_i - \xi)^j \con \sum_{i=1}^k (y_i - \xi)^j \modd{p^j}, \qquad 1 \leq j \leq k.
\eeq
Now we will make the assumption that for the main contribution we may focus only on those $(x_1,\ldots,x_k)$ that are, in Wooley's terminology, \emph{well-conditioned}, that is, they are distinct modulo $p$. (In this setting, solutions that are not well-conditioned satisfy other restrictions that make them sufficiently small in number to be controlled; see e.g. \cite[Ch. 4]{Mon94} for a complete proof.)
If we fix any integral tuple $(n_1,\ldots,n_k)$, then a classical result called Linnik's Lemma \cite[Ch. 4 Lemma 3]{Mon94} shows there are at most $k! p^{\frac{1}{2}k(k-1)}$ well-conditioned solutions $(x_1,\ldots, x_k)$ with $1 \leq x_i \leq p^k $ to the system 
\beq\label{p_classical} 
\sum_{i=1}^k (x_i - \xi)^j \con n_j \modd{p^j}, \qquad 1 \leq j \leq k.
\eeq
In particular, after choosing $(y_1,\ldots, y_k)$ freely (with $ \leq X^k$ possibilities), we obtain at most $k! p^{\frac{1}{2}k(k-1)}$ well-conditioned solutions $(x_1,\ldots, x_k)$ modulo $p^k$ to (\ref{xyp}). (We have ``lifted'' knowledge modulo $p^j$ for various $j \leq k$ to knowledge modulo $p^k$, which is stronger.) Now we use our crucial hypothesis that $p^k >X$ so that the residue class of each $x_i$ modulo $p^k$ actually uniquely identifies it as an integer. 
Now with all the $x_i,y_i$ fixed as integers, one can argue from (\ref{xyuv}), again using translation-dilation invariance, that the number of solutions $u_\ell,v_\ell$ to this system is accounted for by $J_{k',k}(X/p)$, providing the final factor in the main contribution (\ref{Jsk_relation2}).

Here we have followed the presentation of Wooley's beautiful first efficient congruencing paper \cite{Woo12a}, which then goes on to give a clear heuristic comparison of the classical method to efficient differencing and then efficient congruencing. We take a different approach here, and instead focus on demonstrating efficient congruencing in action in the simplest case, that of the degree $k=2$ Vinogradov system.

\subsection{Efficient congruencing for the degree $k=2$ case}\label{sec_EC_k2}
To demonstrate the basic principles of efficient congruencing, we will prove the Main Conjecture for degree $k=2$ with critical index $s_k=3$:
\begin{theo}[Main Conjecture for $k=2$]\label{thm_MC_k2_EC}
\[ J_{3,2}(X) \ll_\ep X^{3+\ep} \]
for all $X \geq 1$ and every $\ep>0$.
\end{theo}
(Of course, we have already seen this is true by elementary means.)

 \subsubsection{The building block lemmas for $k=2$}
For simplicity we temporarily let $J(X)$ denote $J_{3,2}(X)$, and we suppose that $p \geq 3$ is a fixed prime.
First we show that counting solutions with congruence restrictions can be controlled by counting solutions without congruence restrictions, in a smaller range:
\begin{lemm}\label{k2_EC_lemma1}
For any $a,b \geq 1$, if $p^b \leq X$ then 
\[ I_0(X;a,b) \leq J(2X/p^b).\]
\end{lemm}
This is a result of translation-dilation invariance.
Next we show that in the other direction, $J(X)$ can be controlled by counting solutions with certain congruence conditions:
\begin{lemm}\label{k2_EC_lemma2}
If $p\leq X$ then
\[J(X) \ll p J(2X/p) + p^6 I_1(X;1,1).\]
\end{lemm}
Clearly these two lemmas may be compared to the building blocks in the classical method, (\ref{Jsk_relation2}) and  (\ref{Jsk_relation1}), respectively.
We now diverge from the classical method, instead capitalizing on a congruence restriction modulo $p^b$ to build a congruence restriction modulo the higher power $p^{2b}$:
\begin{lemm}\label{k2_EC_lemma3}
If $1 \leq a \leq 2b$  then 
\[ I_1(X;a,b) \leq p^{2b-a} I_1(X;2b,b).\]
\end{lemm}

Finally, we show that on the right-hand side of Lemma \ref{k2_EC_lemma3} the roles of $b$ and $2b$ can be reversed:
\begin{lemm}\label{k2_EC_lemma4}
If $1 \leq a \leq 2b$ and $p^b \leq X$ then 
\[ I_1(X;a,b) \leq p^{2b-a} I_1(X;b,2b)^{1/2} J(2X/p^b)^{1/2}.\]
\end{lemm}
This lemma is the key step to building a recursion relation that will allow us to pass to congruence restrictions modulo ever higher powers of $p$: first $p^b$, then $p^{2b}$, $p^{4b}$, and so on until $p^{2^nb}$, for any $n$ we choose. We may compare this to the classical argument, which in the case of $k=2$ would encounter only congruences modulo  $p^b$ with $b=1,2$ in the system (\ref{p_classical}).

\subsubsection{The iteration for $k=2$}
Before proving the lemmas, we show how to assemble them into a recursion argument that proves Theorem \ref{thm_MC_k2_EC}. Trivial upper and lower bounds for $J(X)$ in this setting are $X^3 \ll J(X) \ll X^6$. Thus we can define 
\beq\label{Delta_dfn_k2_EC}
 \Delta = \inf \{ \del \in \R : J(X) \ll X^{3+\del} \; \text{for all $X \geq 1$} \},
 \eeq
in which case we know that $\Del \in [0,3]$ and that $J(X) \ll_\ep X^{3+\Del+\ep}$ for all $\ep>0$. To prove Theorem \ref{thm_MC_k2_EC}, we will show that $\Del=0$.

First, we claim that for all $b \geq 1$, for all $n \geq 0$, 
\beq\label{k2_EC_induction}
 I_1(X;b,2b) \ll_{b,n,\ep} X^{3+\Del+\ep} p^{-n\Del b},
 \eeq
for every $\ep>0$, as long as $p^{2^nb} \leq X$.
 For each fixed $b \geq 1$ we prove this by induction on $n$ as follows.  The case $n=0$ is trivially true since $I_1(X;b,2b)$ counts solutions with congruence restrictions, so of course
 \[I_1(X;b,2b) \leq J(X) \ll_\ep X^{3+\Del+\ep}.\]
 Now assuming (\ref{k2_EC_induction}) for $n$, we apply Lemma \ref{k2_EC_lemma4} followed by the induction hypothesis (assuming $p^{2^{n+1}b} \leq X$) to write 
 \begin{eqnarray*}
  I_1(X;b,2b)& \leq & 		p^{4b-b} I_1(X;2b,4b)^{1/2}J(2X/p^{2b})^{1/2}	\\
  	& \ll_{b,n,\ep} & p^{3b} (X^{3+\Del+\ep} p^{-n\Del (2b)})^{1/2} ( (2X/p^{2b})^{3 +\Del+\ep})^{1/2} \\
	& \ll_{b,n,\ep} & X^{3+\Del+\ep}  p^{-(n+1)\Del b} ,
  \end{eqnarray*}
which completes the induction.

Now we describe the recursion with which we bound $J(X)$: by applying Lemma \ref{k2_EC_lemma2}, then  Lemma \ref{k2_EC_lemma4} with $a=1,b=1$, then (\ref{k2_EC_induction}) with $b=1$, we see that 
\begin{eqnarray*}
J(X) & \ll & p J(2X/p) + p^6 I_1(X;1,1) \\
	& \leq & pJ(2X/p) + p^7 I_1(X;1,2)^{1/2} J(2X/p)^{1/2} \\
	& \ll_{n,\ep} & pJ(2X/p) + p^7 \{X^{3+\Del+\ep} p^{-n\Del} \}^{1/2} \{ (2X/p)^{3 + \Del + \ep} \}^{1/2} \\
	& \ll_{n,\ep} & pJ(2X/p) + X^{3+\Del+\ep} p^{\frac{11}{2} - \frac{(n+1)\Del}{2}},
\end{eqnarray*}
for all $n \geq 1$ such that $p^{2^{n}} \leq X$. 
We now fix a prime $p \geq 3$ with 
\beq\label{p_large}
 \frac{1}{2} X^{1/2^n}  \leq p \leq X^{1/2^n},
 \eeq
which exists by Bertrand's postulate (or the prime number theorem) as long as $X \geq 10^{2^n}$.
Then we have shown that 
\[ J(X)
 \ll X^{3 + \Del + \ep} p^{-2-\Del} + X^{3+\Del+\ep} p^{\frac{11}{2} - \frac{(n+1)\Del}{2}} ,\]
for every $n\geq1$ and $\ep>0$, as $X \maps \infty$.
Now, assuming that $\Del>0$, we may choose $n$ such that $(n+1)\Del \geq 12$, so that $J(X) \ll_\ep X^{3+\Del+\ep} p^{-1/2},$ for every $\ep>0$. This contradicts the definition of $\Del$, due to (\ref{p_large}), so we must have that $\Del=0$, as desired.

Philosophically, what powered this iterative engine? We showed that, given a suspiciously large number of integral solutions to the Vinogradov system, they could be concentrated in such a short $p$-adic interval (corresponding to restrictions modulo increasingly high powers of $p$) that we obtained a contradiction to a known bound.

\subsubsection{Proof of the lemmas for $k=2$}
For completeness, we indicate the simple proofs of the four lemmas, and then turn to a comparison to the mechanism of efficient congruencing in the nontrivial case $k=3$.

To prove Lemma \ref{k2_EC_lemma1}, we note that by the definition of $I_0(X;a,b)$, there exists some $\eta \modd{p}$ for which $I_0(X;a,b)$ counts solutions of (\ref{Vin_sys}) such that $x_i \con \eta \modd{p^b}$ for all $1 \leq i \leq 2s$, that is, each $x_i$ may be written as $x_i = \eta + p^b y_i$ where $0 \leq y_i \leq X/p^b$. We now set $z_i = y_i + 1$, so that $1 \leq z_i \leq 1+X/p^b \leq 2X/p^b$, due to the assumption that $p^b \leq X$. Then by translation-dilation invariance, the set of $y_i$ and the set of $z_i$ are each solutions to (\ref{Vin_sys}) as well, so that $I_0(X;a,b) \leq J(2X/p^b)$, as claimed.

To prove Lemma \ref{k2_EC_lemma2}, we first count those solutions $x_i$ to (\ref{Vin_sys}) such that $x_1 \con \cdots \con x_6 \modd{p}$. The number of such solutions is 
\[ \sum_{\eta \modd{p}} I_0(X;\cdot,\eta; \cdot,1) \leq p I_0(X;\cdot,1) \leq p J(2X/p) \]
by Lemma \ref{k2_EC_lemma1}; this contributes the first term in Lemma \ref{k2_EC_lemma2}.
We now consider the remaining solutions, for which we know that for some $i\neq j$ we have $x_i \not\con x_j \modd{p}$. In particular, the contribution from these solutions is at most 
\[ \binom{6}{2}\sum_{\xi \not\con \eta \modd{p}} \int_{(0,1]^2} |f_1(\albf,\xi)| |f_1(\albf;\eta)| |f(\albf)|^4 d\albf.\]
There are at most $p(p-1)$ terms in the sum; for each fixed $\xi,\eta$, we apply H\"{o}lder's inequality to see that 
the integral is bounded above by
\[ ( \int |f_1(\albf;\xi)|^2 |f_1(\albf;\eta)|^4 d\albf )^{1/6}
		 ( \int |f_1(\albf;\xi)|^4 |f_1(\albf;\eta)|^2 d\albf )^{1/6}
		 (\int |f(\albf)|^6 d\albf)^{2/3}.
		\]
This is at most $I_1(X;1,1)^{1/6}I_1(X;1,1)^{1/6}J(X)^{2/3}$, and in total we conclude that 
\[ J(X) \ll pJ(2X/p) + p^2I_1(X;1,1)^{1/3}J(X)^{2/3}.\]
If the first term dominates, then the bound in Lemma \ref{k2_EC_lemma2} holds; if otherwise the second term dominates, then
\[ J(X) \ll p^2I_1(X;1,1)^{1/3}J(X)^{2/3},\]
from which we deduce that $J(X) \ll p^6I_1(X;1,1)$, leading to the second term in Lemma \ref{k2_EC_lemma2}.

To prove Lemma \ref{k2_EC_lemma3}, recall that for fixed $\xi,\eta$, 
$I_1(X;\xi,\eta;a,b)$ counts solutions to (\ref{Vin_sys}) such that 
$x_i = \xi + p^a y_i$ for $i=1,4$ and $x_i = \eta + p^b y_i$ for $i=2,3,5,6$.
Defining $\nu = \xi - \eta$, we see by translation invariance that the variables $z_i$ defined by 
\[ z_i = \nu + p^a y_i, \;\; i=1,4; \qquad z_i = p^by_i, \;\; i = 2,3,5,6 \]
are also solutions to (\ref{Vin_sys}). 
Now upon examining the quadratic equation in the system (\ref{Vin_sys}) in terms of the variables $z_i$, we learn that
\[ (\nu+p^ay_1)^2 \con (\nu + p^a y_4)^2 \modd{p^{2b}}.\]
Recall that in $I_1(X;1,1)$ it is specified that $\xi \not\con \eta \modd{p}$, so $p \ndiv \nu$; this allows us to deduce 
$\nu+p^ay_1 \con \nu + p^a y_4 \modd{p^{2b}}$ and hence $y_1 \con y_4 \modd{p^{2b-a}}$.
Once we fix a choice for $y_4$ modulo $p^{2b-a}$ (of which there are $p^{2b-a}$ choices) it thus fixes $y_1$ modulo $p^{2b-a}$, and we consequently see that $x_1$ and $x_4$ are fixed modulo $p^{2b}$, say $x_1 \con x_4 \con \xi' \modd{p^{2b}}$.
We have shown that for any $\xi \not\con \eta \modd{p}$, 
\[ I_1(X;\xi,\eta,a,b) \leq p^{2b-a}I_1(X;2b,b),\]
which suffices for Lemma \ref{k2_EC_lemma3}.

To obtain Lemma \ref{k2_EC_lemma4} we start from Lemma \ref{k2_EC_lemma3} and observe that for any $\xi,\eta$, by Cauchy-Schwarz,
\begin{eqnarray*}
 I_1(X;\xi,\eta,2b,b) 
 	&=& \int_{(0,1]^2} |f_{2b}(\albf;\xi)|^2|f_b(\albf;\eta)|^4 d\albf \\
	& \leq & ( \int_{(0,1]^2} |f_{2b}(\albf;\xi)|^4|f_b(\albf;\eta)|^2 d\albf )^{1/2}( \int_{(0,1]^2} |f_b(\albf;\eta)|^6 d\albf )^{1/2}.
	 \end{eqnarray*}
The first term is $I_1(X;\eta,\xi;b,2b)^{1/2}$ while the second term is $I_0(X;\cdot, \eta;\cdot,b)^{1/2}$, which by Lemma \ref{k2_EC_lemma1} is bounded above by $J(2X/p^b)^{1/2}$, as long as $p^b \leq X$. Taking the supremum over $\xi,\eta$, we have proved the claim of Lemma \ref{k2_EC_lemma4}.

\subsection{Efficient congruencing for the degree $k=3$ case}\label{sec_k3_EC}
Theorem \ref{thm_MC_k2_EC}, the Main Conjecture for $k=2$, was a useful model for demonstrating a basic principle of efficient congruencing, but was a trivial result. On the other hand, the Main Conjecture for degree $k=3$ is not trivial at all, and was only finally obtained by Wooley in \cite{Woo16c}. However, a simplified presentation of Wooley's argument in the cubic case, due to Heath-Brown \cite{HB15x}, gives a slim 10-page proof of:
\begin{theo}[Main Conjecture for $k=3$]\label{thm_MC_k3_EC}
\[ J_{6,3}(X) \ll_\ep X^{6+\ep} \]
for all $X \geq 1$ and every $\ep>0$.
\end{theo}
In fact  our description for $k=2$ was written exactly parallel to Heath-Brown's presentation, so that we may now make a few quick comparisons to the (nontrivial) case of $k=3$. Let $J(X)$ denote the critical case $J_{6,3}(X)$, define all quantities relative to $s=6$, and assume $p \geq 5$ is a fixed prime. Analogous to the $k=2$ case, we start with a lemma that shows solutions with congruence restrictions can be counted by solutions without congruence restrictions, in a shorter range:
\begin{lemm}\label{k3_EC_lemma1}
For any $a,b \geq 1$, if $p^b \leq X$ then 
\[ I_0(X;a,b) \leq J(2X/p^b).\]
\end{lemm}
We have a lemma that introduces a congruence restriction:
\begin{lemm}\label{k3_EC_lemma2}
If $p\leq X$ then
\[J(X) \ll p J(2X/p) + p^{12} I_2(X;1,1).\]
\end{lemm}
Note that now the new quantity $I_2(X;a,b)$ plays a role, shifting the weights of the exponents in the bilinear expression (\ref{EC_bilinear}).
We have two further lemmas that show how to pass from restrictions modulo a lower power of $p$ to a higher power of $p$:
\begin{lemm}\label{k3_EC_lemma5}
If $1 \leq a \leq 3b$  then 
\[ I_1(X;a,b) \leq p^{3b-a} I_1(X;3b,b).\]
\end{lemm}
\begin{lemm}\label{k3_EC_lemma6}
If $1 \leq a \leq b$, 
\[ I_2(X;a,b) \leq 2b p^{4(b-a)} I_2(X;2b-a,b).\]
\end{lemm}
And finally there are two lemmas that show how to pass between $I_1$ and $I_2$ and reverse the weights of the powers of $p$ in the bilinear integral (\ref{EC_bilinear}):
\begin{lemm}\label{k3_EC_lemma3}
\[I_2(X;a,b) \leq I_2(X;b,a)^{1/3}I_1(X;a,b)^{2/3}\]
\end{lemm}
\begin{lemm}\label{k3_EC_lemma4}
If $p^b \leq X$ then 
\[ I_1(X;a,b) \leq I_2(X;b,a)^{1/4} J(2X/p^b)^{3/4}.\]
\end{lemm}

With these lemmas in hand, upon setting 
\[ \Del = \inf \{ \del \in \R : J(X) \ll X^{6+\del} \; \text{for all $X \geq 1$}\},\]
 Heath-Brown establishes a recursion analogous to (\ref{k2_EC_induction}), 
in the shape
 \beq\label{HB_recursion}
 I_2(X;a,b) \ll_{a,b,n,\ep} X^{6+\Del + \ep}p^{-5(b-a)} p^{-\frac{n\Del}{6}(3b-a)},
 \eeq
for all $1 \leq a \leq b$ and every $n \geq 0$, provided that $p^{3^nb} \leq X$.
(For the philosophy leading to this recursion relation, see \cite[\S 4]{HB15x}.)
From this, after choosing $p \geq 5$ with $\frac{1}{2}X^{1/3^n} \leq p \leq X^{1/3^n}$ (which exists if $X \geq 10^{3^n}$), 
one deduces that 
\[ J(X) \ll pJ(2X/p)  + p^{12}I_2(X;1,1) \ll p(X/p)^{6+\Del + \ep}  + X^{6+\Del + \ep} p^{12 - n\Del/3},\]
so that if $\Del$ were strictly positive, by choosing $n$ sufficiently large that $n\Del/3>13$ and using the size of $p$, this bound would violate the definition of $\Del$. 

Superficially this argument appears as simple as our model case $k=2$, but  nontrivial work is hidden in the structure of the recursion (\ref{HB_recursion}), and in the proof of Lemma \ref{k3_EC_lemma6}; this can be reduced to counting solutions to a certain pair of congruences in four variables \cite[Lemma 8]{HB15x}, and these solutions fall into two cases: nonsingular and singular, with a solution being nonsingular if two associated gradient vectors are not proportional modulo $p$ (which we may  regard as a kind of transversality). In the case of Lemma \ref{k3_EC_lemma6}, the number of singular solutions turns out to be the same order of magnitude as the number of nonsingular solutions, and a streamlined argument  succeeds. 

But as Heath-Brown points out, for  Vinogradov systems of degree $k \geq 4$, the number of singular solutions can exceed the number of nonsingular solutions in arguments of this style, and one of Wooley's impressive technical feats in the general case is to carry out a ``conditioning'' step, which removes singular solutions from consideration (thus facilitating later uses of Hensel's Lemma, which ``lifts'' solutions of congruences modulo lower powers of $p$ to solutions of congruences modulo higher powers of $p$). This conditioning process is especially difficult, since even in the general case, the single auxiliary prime $p$ is fixed once and for all at the beginning of the argument. Once this has been selected, a basic philosophy at the heart of the efficient congruencing method is to use an intricate iteration to successively extract congruence restrictions modulo $p^{kb}$, then $p^{k^2b}$, then $p^{k^3b}$, then $p^{k^4b}$, and so on, in order to concentrate putative integral solutions of the Vinogradov system into successively shorter $p$-adic neighborhoods.

\section{Setting the stage for decoupling}\label{sec_decoupling_intro}
We now turn to our second principal focus: the approach to the Main Conjecture via $\ell^2$ decoupling. To motivate the form that decoupling theorems take, and to introduce several fundamental concepts that underpin the recent proofs of decoupling theorems,  we first recall several key concepts: orthogonality of functions, Littlewood-Paley theory, square functions, and restriction problems.

\subsection{Orthogonality}

Suppose that we consider functions in the Hilbert space $L^2(\R^n)$ with inner product $\langle f,g \rangle = \int_{\R^n} f(x) \overline{g(x)} dx$ (in fact any Hilbert space will do). 
If we have a countable collection of orthogonal functions $\{f_j\}$ in this space, so that $\int f_i \overline{f_j}=0$ for $i \neq j$, we see that for any finite collection $I$ of indices, 
\beq\label{Hilbert_L2_ex1} \| \sum_{j\in I} f_j \|_{L^2}  =(  \sum_{j \in I} \|f_j\|_{L^2}^2 )^{1/2}.
\eeq
Indeed,
\[ \| \sum_{j \in I} f_j \|_{L^2} = (\int (\sum_{i \in I} f_i) (\overline{\sum_{j \in I} f_j}) )^{1/2}
	 = ( \int \sum_{i,j \in I} f_i \overline{f_j} )^{1/2}
	 = ( \sum_{j \in I} \int |f_j|^2 )^{1/2}.\]
On the other hand, if we make no assumption about orthogonality, we could trivially apply the triangle inequality and obtain 
 \[ \| \sum_{j \in I} f_j \|_{L^2} \leq \sum_{j \in I} \|f_j \|_{L^2}.  \]
 To make this more comparable to our previous observation, we can further apply the Cauchy-Schwarz inequality to see that without any assumption on orthogonality, 
  \beq\label{Hilbert_L2_ex2}
   \| \sum_{j \in I} f_j \|_{L^2} \leq |I|^{1/2} (\sum_{j \in I} \|f_j \|^2_{L^2})^{1/2} .
   \eeq
This is certainly  much weaker than our conclusion (\ref{Hilbert_L2_ex1}), which we say exhibits square-root cancellation, as a result of the orthogonality assumption. 
  
  What are sources of orthogonality? As an example, within $L^2([0,1])$ with inner product $\langle f, g \rangle = \int_{[0,1]} f(x) \overline{g(x)}dx$, if we set $f_j(x) = e( jx) = e^{2\pi i j x}$ for $j \in \Z$ then we can see that the collection $\{f_j\}$ is orthogonal. Or, on $L^2 (\R^n)$, recall that the Fourier transform is a unitary operator, with  Plancherel theorem 
  \beq\label{polarized_Plancherel}
  \int_{\R^n} f(x) \overline{g(x)}dx = \int_{\R^n} \hat{f}(x) \overline{\hat{g}(x)} dx, \qquad f,g \in L^2(\R^n),
  \eeq
yielding the important  special case $\|f \|^2_{L^2}  = \| \hat{f} \|^2_{L^2}$.
The identity (\ref{polarized_Plancherel}) shows that orthogonality arises when the Fourier transforms $\hat{f}$ and $\hat{g}$ behave sufficiently differently that their inner product is zero. A particularly nice case occurs when the supports of $\hat{f}$ and $\hat{g}$ are disjoint sets of $\R^n$, and then we see immediately from (\ref{polarized_Plancherel}) that $f$ and $g$ must be orthogonal.\footnote{Recall that the support of a function is the closure of the set upon which it takes nonzero values.}
 
\subsection{Littlewood-Paley theory and square functions}
 This observation leads to the basic principles underlying a useful technique called Littlewood-Paley theory, dating back to \cite{LitPal31,LitPal36}, which decomposes a given function $f$ by dissecting $\hat{f}$ into pieces.
 Let us consider the case of a function $f$ on $\R$. For each integer $j$ let $\Del_j(\xi)$ denote the indicator function for the one-dimensional annulus $\{2^j \leq |\xi| < 2^{j+1}\}$, and define an operator $P_j$ that acts on a function $f$ by 
  \[ (P_j f) (x) = ( \Del_j(\xi) \hat{f})\check{\;} (x);\]
 this modifies $f$ by ``projecting'' (or ``restricting'') $\hat{f}$ onto the chosen annulus.
  In particular, if $f \in L^2(\R)$ then we note that the functions $f_j = P_j f$ have disjoint Fourier supports and hence are orthogonal functions, so that (\ref{Hilbert_L2_ex1}) holds. 
 Here we note that because of the compatibility of the $L^2$ and $\ell^2$ norms, (\ref{Hilbert_L2_ex1}) can be re-written as 
 \beq\label{Hilbert_L2_ex4}
   \| \sum_{j \in I} f_j \|_{L^2}  = (  \sum_{j \in I} \|f_j\|_{L^2}^2 )^{1/2}= \| ( \sum_j |f_j|^2)^{1/2} \|_{L^2} .
   \eeq
Thus, our argument above may be summarized by the identity
 \beq\label{LP_oneD_L2}
  \|  ( \sum_j |P_j f|^2)^{1/2} \|_{L^2}  = \|f\|_{L^2}.
  \eeq
 The principal insight of Littlewood-Paley theory is that a relationship analogous to (\ref{LP_oneD_L2}) still holds if the $L^2(\R)$ norms are replaced by $L^p (\R)$ norms: for $1<p<\infty$, there exist constants $A_p, B_p$  such that for every $f \in L^p(\R)$,
 \beq\label{LP_oneD_Lp}
A_p \|f\|_{L^p(\R)} \leq  \| ( \sum_j |P_j f|^2)^{1/2} \|_{L^p(\R)} \leq B_p \|f\|_{L^p(\R)}.
 \eeq
 A key advantage of this approach is that even if one is ultimately interested in a rather general function $f$, the Littlewood-Paley decomposition allows one to focus piece by piece on the special type of function $f_j = P_jf$ that has its Fourier transform ``localized'' to a particular region (such as an annulus); we will gain some notion of why this is advantageous in \S \ref{sec_uncertainty}. 
 
Importantly, in higher dimensions, something more subtle is required than the construction above: on $\R^n$ with $n \geq 2$, if a  sharp cut-off function (such as an indicator function) is used to to dissect $\hat{f}$ into pieces supported on $n$-dimensional annuli $\{ 2^j \leq |\xi| < 2^{j+1} \}$, a celebrated result of C. Fefferman \cite{Fef71} shows that an $L^p$ relationship of the kind (\ref{LP_oneD_Lp}) fails to hold for any $p \neq 2$:
 \begin{theo}[Ball Multiplier Theorem]\label{thm_ball_multiplier}
 Define an operator $T$ acting on $f \in L^2(\R^n) \intersect L^p(\R^n)$ by $(Tf)\hat{\:}(\xi) = \onebf_{B}(\xi) \hat{f}(\xi)$, where $\onebf_B$ is the characteristic function of a ball $B$. Then for every  dimension $n \geq 2$, the operator $T$ is a bounded operator on $L^2(\R^n)$ and an unbounded operator on $L^p(\R^n)$ for every $ p \neq 2$.
 \end{theo}
This provided an unexpected answer to a long-standing question on the convergence of truncated Fourier integrals (see e.g. \cite[p. 298]{Tao01}); as we will  later see, the methods Fefferman used are intertwined with several phenomena closely related to decoupling.
 
This also brings  us to the notion of \emph{almost orthogonality}. 
In the context of annuli in $\R^n$, if instead of a sharp cut-off one uses an appropriate choice of a $C^\infty$ bump function $\Del_j$ that is only essentially concentrated on the annulus $\{ 2^j \leq |\xi| < 2^{j+1} \}$ and decays rapidly outside of it,  then even though  the analogous functions $f_j=P_j f$ do not strictly speaking have Fourier transforms with disjoint supports, they are sufficiently dissociated that the $f_j$ are said to be \emph{almost orthogonal}, and one can again prove a relation of the form (\ref{LP_oneD_Lp}) for $1<p<\infty$ (see e.g. \cite[Ch. 6]{Gra14a}). All told, there is a vast array of Littlewood-Paley methods, which are key tools in the theory of singular integral operators, maximal operators, and multiplier operators; for an overview see for example \cite[Ch. IV]{SingInt}, \cite[Ch. 8]{Duo01}.\footnote{The subtlety of constructing truncations on the Fourier side is revealed by the fact that in some situations, sharp cut-offs do work well, such as: in $\R^1$ dissecting $\hat{f}$ according to arbitrary disjoint intervals \cite{RdF85},  in $\R^2$ dissecting $\hat{f}$ according to disjoint sectors in the plane \cite{NSW78}; in $\R^n$ for any $n$, dissecting $\hat{f}$ according to disjoint dyadic boxes, with sides parallel to the axes \cite[Ch IV\S5]{SingInt}. In each of these cases the $L^2$ result (\ref{LP_oneD_L2}) holds immediately by orthogonality, and the analogue to (\ref{LP_oneD_Lp}) is the key result.}

For our purposes, Littlewood-Paley theory also introduces a relevant object called a \emph{square function}: given a collection of operators $T_j$ and an appropriate function $f$,  the square function is the operator
 \[ f \mapsto ( \sum_{j} |T_j f|^2 )^{1/2} \]
and a fundamental question is whether this operator is bounded on certain $L^p$ spaces. 
 As in  (\ref{Hilbert_L2_ex4}), the case $p=2$ is special, since
 \[ \|  ( \sum_{j} |T_j f|^2 )^{1/2} \|_{L^2} = (\int \sum_j |T_jf(x)|^2 dx )^{1/2}= ( \sum_j \|T_j f\|_{L^2}^2)^{1/2}.\]
For other $L^p$ spaces, one may aim to prove a Littlewood-Paley inequality analogous to the right-hand inequality in (\ref{LP_oneD_Lp}). The left-hand estimate in (\ref{LP_oneD_Lp}) is called  a reverse Littlewood-Paley inequality or a \emph{square function estimate}, which we can state with $f_j = T_j f$ as the claim:
 \beq\label{square_function_Lp_estimate}
  \| \sum_j f_j \|_{L^p} \leq A \|  ( \sum_{j} |f_j|^2 )^{1/2}  \|_{L^p}. 
 \eeq
Square functions appear in many guises throughout harmonic analysis, and in particular play a fundamental role in methods applying Littlewood-Paley theory or related ideas. (For a historical overview, see \cite{Stein82}; the reader in search of a particularly transparent example of the power of a square function may read the proof of the $L^2$ boundedness of the maximal function on the parabola, \cite[Ch. XI \S 1.2]{SteinHA}.)
 
 \subsection{A first notion of $\ell^2$ decoupling}
For the moment, we will say (somewhat informally) that a collection of functions $\{f_j\}_{j \in I}$ in $L^p$ satisfies an \emph{$\ell^2$ decoupling inequality for $L^p$} if 
 \beq\label{dfn_generic_decoupling}
  \| \sum_{j \in I} f_j \|_{L^p} \ll_{p,\ep} |I|^\ep ( \sum_{j \in I} \|f_j\|_{L^p}^2 )^{1/2},
  \eeq
 for every $\ep>0$, with an implied constant independent of $I$ and the functions $f_j$.
 The terminology specifies $\ell^2$  because of the appearance of the discrete $\ell^2$ norm $\| \{a_j\}\|_{\ell^2} = (\sum_j |a_j|^2)^{1/2}$, and one could more generally study $\ell^r$ decoupling for $2 \leq r <\infty$.
 We will formally state in \S \ref{sec_decoupling_VMVT} the $\ell^2$ decoupling inequality developed by Bourgain, Demeter and Guth, but first we give some motivating examples.

We may compare (\ref{dfn_generic_decoupling}) to our trivial statement (\ref{Hilbert_L2_ex2}), which applies to  functions in any Banach space, and which we now write more generally as 
 \beq\label{Hilbert_L2_ex3}
   \| \sum_{j \in I} f_j \|_{L^p} \leq |I|^{1/2} (\sum_{j \in I} \|f_j \|^2_{L^p})^{1/2} 
   \eeq
 for any collection $f_j \in L^p$, and any $1 \leq p \leq \infty$. We see then that an $\ell^2$ decoupling inequality for $L^p$ saves over (\ref{Hilbert_L2_ex3}) by nearly a factor $|I|^{1/2}$; that is, it exhibits square-root cancellation. We now examine four examples of $\ell^2$ decoupling, as popularized in talks of Demeter and e.g. Tao \cite{Tao15ablog}.

 \subsubsection{Example 1: orthogonal and almost orthogonal functions}\label{sec_almost_orthogonal_example}
Any collection of orthogonal functions in $L^2$ satisfies an $\ell^2$ decoupling inequality for $L^2$ since they satisfy the stronger identity (\ref{Hilbert_L2_ex1}), which we might think of as ``perfect decoupling.'' 
Recalling that functions with disjoint Fourier supports are orthogonal, suppose we instead consider a weaker condition, that a family $\{f_j\}$ is a collection of functions such that the supports of the Fourier transforms $\{ \hat{f_j}\}$  have finite (or bounded) overlap: that is to say, there exists a uniform constant $C_0$ such that for any $\xi$, at most $C_0$ of the functions $\hat{f_j}(\xi)$ are nonzero at $\xi$. Such a collection is  \emph{almost orthogonal}, and we may see immediately that an $\ell^2$ decoupling inequality for $L^2$ still holds under this weaker assumption. Indeed, letting $\omega_j$ denote the support of $\hat{f_j}$, by Plancherel's theorem followed by the Cauchy-Schwarz inequality,
\[ \| \sum_j f_j \|_{L^2} = \| \sum_j \hat{f_j} \|_{L^2} \leq \|  ( \sum_j |\hat{f_j}|^2)^{1/2} ( \sum_j |\onebf_{\omega_j}|^2)^{1/2}  \|_{L^2}. \]
The bounded overlap assumption implies the uniform bound $\sum_j |\onebf_{\omega_j}(\xi)|^2 \leq C_0$ as a function of $\xi$, 
so that 
\[ \| \sum_j f_j \|_{L^2}  \leq  C_0^{1/2} \|  ( \sum_j |\hat{f_j}|^2)^{1/2} \|_{L^2} = C_0^{1/2}( \sum_j \|\hat{f_j} \|^2_{L^2} )^{1/2},\]
upon recalling the way that $\ell^2$ and $L^2$ norms interact via the identity (\ref{Hilbert_L2_ex4}). Finally, applying Plancherel again to each term shows that (\ref{dfn_generic_decoupling}) holds for $L^2$. This argument has quite a bit of flexibility, and given a collection of $N$ functions $f_j$, we could have allowed up to $N^\ep$ to have overlapping Fourier support at any point, and still obtained (\ref{dfn_generic_decoupling}) for $L^2$.

\subsubsection{Example 2: comparison to the square function estimate}\label{sec_decoupling_square_function}
A decoupling inequality in the form (\ref{dfn_generic_decoupling}) is clearly much stronger than (\ref{Hilbert_L2_ex3}). How would it compare to a square function estimate of the form (\ref{square_function_Lp_estimate}), in which the roles of the $L^p$ and $\ell^2$ norms are reversed? 
In fact, for $p > 2$ the square function estimate is stronger, because we can verify directly that for any collection of functions $f_j$,
\[ \| ( \sum_j |f_j|^2)^{1/2} \|_{L^p} \leq ( \sum_j \|f_j\|_{L^p}^2)^{1/2}. \]
To see this, we re-write this claim in the equivalent form
\[ \left( \int \left| \sum_j |f_j|^2  \right|^{p/2} \right)^{2/p} \leq \sum_j \left( \int  \left| \; |f_j|^2  \right|^{p/2} \right)^{2/p}.\]
This is an identity if $p=2$, and holds for $p > 2$ by Minkowski's inequality for $L^q$ norms with $q=p/2>1$, applied to the functions $\tilde{f}_j = |f_j|^2$. In general, as we will see in further examples, a square function estimate is extremely powerful. But given the current intractability of certain square function estimates of great interest in harmonic analysis, we are motivated to study weaker replacements, such as $\ell^2$ decoupling. (Indeed, Wolff \cite{Wol00} was motivated by a conjecture that would follow from an unproved square function estimate, and he made progress by formulating a new $\ell^p$ decoupling question.)

\subsubsection{Example 3: bi-orthogonality}\label{sec_bi_orthogonality}
Consider for $j \in \{1,\ldots, N\}$ the function $f_j(x) = e(j^2 x)$. These functions satisfy an $\ell^2$ decoupling theorem in $L^4([0,1])$, namely,
\beq\label{decoupling_ex1}
 \| \sum_{j=1}^N e(j^2 x) \|_{L^4([0,1])} \ll N^\ep ( \sum_{j=1}^N \| e(j^2 x) \|_{L^4([0,1])}^2)^{1/2}.
 \eeq
Here we clearly see the motivation for the terminology ``decoupling:'' on the left-hand side, the function $\sum_{j=1}^N e(j^2 x)$ comprises many different frequencies, while on the right-hand side, the frequencies have been separated, or decoupled, and  each frequency $e(j^2x)$ is considered independently.

To verify (\ref{decoupling_ex1}), after raising both sides to the 4th power and applying orthogonality to compute the integral on the left, it is equivalent to show that
\[ \# \{ 1 \leq x_1,\ldots,x_4 \leq N: x_1^2 + x_2^2  = x_3^2 + x_4^2 \} \ll N^{2+\ep}.\]
This is true: after choosing any two variables freely, say $x_1,x_2$, yielding a fixed value $c$ for $x_1^2 + x_2^2$, then for any $\ep>0$ there are at most $c^\ep \ll N^{2\ep}$ choices for $x_3, x_4$ such that $x_3^2 + x_4^2 = c$.\footnote{In general, for any integer $n$ there are at most $d(n)$ representations of $n$ as a sum of two squares, where $d(n)$ is the divisor function. This may be seen by noting that $n=x^2+y^2$ is equivalent to stating that $n$ is the norm of a Gaussian integer $x+iy$ in $\mathbb{Q}(i)$, and at most $d(n)$ Gaussian integers may have norm $n$; see \cite[Lemma 4.2]{Nar80}. One then finally applies the fact that $d(n)\ll_\ep n^\ep$ for every $\ep>0$, as in \cite[Thm. 315]{HarWri75}.
} 
Moreover, a more general phenomenon is underfoot: we can regard the above $L^4$ result as an $L^2$ bound for the functions $f_i f_j$. For any bi-orthogonal collection of functions, so that the products $\{f_if_j\}_{i,j}$ are pairwise orthogonal for $(i,j) \neq (i',j')$, $\ell^2$ decoupling holds for $L^4$, since
\[ \| \sum_{j} f_j \|^4_{L^4} = \| ( \sum_{j} f_j)^2 \|^2_{L^2} = \| \sum_{i,j} f_{i} f_j \|^2_{L^2} =  \sum_{i,j,i',j'}  \int f_if_j \overline{f_{i'}f_{j'}} = \sum_{i,j} \int |f_if_j|^2 .
\]
Applying Cauchy-Schwarz, we see that
\[ \| \sum_{j} f_j \|^4_{L^4}\leq \sum_{i,j} (\int |f_i|^4 )^{1/2} (\int |f_j|^4 )^{1/2}   = ( \sum_j \|f_j\|_{L^4}^2)^2,
\] 
which, after taking $4$-th roots, reveals $\ell^2$ decoupling for $L^4$. 
As in Example 2, this argument still proves decoupling under the weaker assumption that each product $f_if_j$ is orthogonal to all but $O(N^\ep)$ of the other products. (Similarly, any family of tri-orthogonal functions satisfies $\ell^2$ decoupling for $L^6$, and so on.)

\subsubsection{Example 4: lacunarity}
Consider for $j \in \{1,\ldots, N\}$ the function $f_j(x) = a_j e(2^j x)$ with a coefficient $a_j \in \C$. These  functions satisfy an $\ell^2$ decoupling theorem  in $L^p([0,1])$ for each even integer $p=2k$.
In fact the following expressions are comparable:
\[ \| \sum_{j=1}^N a_j e(2^j x) \|_{L^p([0,1])} \approx ( \sum_{j=1}^N \| a_j e(2^j x) \|_{L^p([0,1])}^2)^{1/2}.\]
Raising the claim  to the $p$-th power and using  orthogonality, we see that the left-hand side becomes 
\beq\label{decoupling_ex2}
\sum_{1 \leq j_1,\ldots, j_{2k} \leq N} a_{j_1} \cdots a_{j_k}\overline{a_{j_{k+1}}} \cdots \overline{a_{j_{2k}}} \del_{\mathbf{j}}
\eeq
where $\del_{\mathbf{j}}=1$ if $2^{j_1} + \cdots + 2^{j_{k}} = 2^{j_{k+1}} + \cdots+ 2^{j_{2k}}$ and zero otherwise.
By the uniqueness of representations in base 2, this identity occurs precisely when $\{j_1,\ldots, j_k\} = \{j_{k+1},\ldots, j_{2k}\}$ as sets. Thus  (\ref{decoupling_ex2}) is comparable to $(\sum_{j} |a_j|^2)^k$, and upon taking $2k$-th roots, this proves the claim for $p=2k$. 

\subsubsection{Example 5: Vinogradov Mean Value Theorem}
From Examples 3 and 4 it is likely already clear that the Main Conjecture for $J_{s,k}(N)$ may be framed as a decoupling estimate for $L^p([0,1]^k)$ at the critical exponent  $p = p_k= 2s_k = k(k+1)$, putatively of the form
\beq\label{decoupling_VMVT_model0}
 \| \sum_{j=1}^N  e(jx_1 + j^2x_2 + \cdots + j^kx_k)\|_{L^p([0,1]^k)} \ll N^\ep ( \sum_{j=1}^N \|e(jx_1 + j^2x_2 + \cdots + j^kx_k)\|_{L^p([0,1]^k)}^2)^{1/2}.
\eeq
To make this patently familiar, we raise both sides to the $p$-th power and observe this is equivalent to the claim
\beq\label{decoupling_VMVT_model}
\int_{(0,1]^k} |\sum_{j=1}^N  e(jx_1 + j^2x_2 + \cdots + j^kx_k)|^{k(k+1)} dx_1\cdots dx_k \ll N^\ep N^{\frac{k(k+1)}{2}},
\eeq
which would show the Main Conjecture in the critical case.

\subsubsection{Limiting cases}
For $p<2$, an $\ell^2$ decoupling inequality for $L^p$ cannot in general be expected to hold, even for a family for which $\ell^2$ decoupling for $L^2$ holds.
For example, supposing that $\{f_j\}$ are disjointly supported functions (and hence orthogonal in $L^2$), then
\[ \| \sum_{j \in I} f_j \|_{L^p}
	 = ( \int |\sum_{j \in I} f_j|^p dx )^{1/p}
	  = ( \int \sum_{j \in I} |f_j|^p dx)^{1/p}
	   = ( \sum_{j \in I} \|f_j\|_{L^p}^p)^{1/p}.\]
Then if an $\ell^2$ decoupling inequality held for $L^p$ for this family, it would assert that 
\beq\label{supposed_p_2}
 ( \sum_{j \in I} \|f_j\|_{L^p}^p)^{1/p} \ll_\ep |I|^\ep ( \sum_{j \in I} \|f_j\|_{L^p}^2)^{1/2},
 \eeq
 for every $\ep>0$,
and such an inequality need not hold in general.\footnote{This is reflective of the behavior of discrete $\ell^p$ spaces, namely $\ell^p \subsetneq \ell^q$ if $0<p <q \leq \infty.$}
At the other extreme, an $\ell^2$ decoupling inequality for $L^p$ is also not expected to hold in general for $p= \infty$. For example take $f_j(x)  =  a_j \Re( e( jx))$ for $j \in \{1,\ldots,N\}$, $a_j \in \R_{\geq 0}$.  Then these functions are orthogonal on $L^2([0,1])$ and all attain their maximum modulus at the point $x=0$.
 As a result, $\|\sum_j f_j \|_{L^\infty} = \sum_j \|f_j \|_{L^\infty}$, and then $\ell^2$ decoupling for $L^\infty$ would assert that 
 \[
\sum_{j \in I} \|f_j\|_{L^\infty} \ll_\ep |I|^\ep ( \sum_{j \in I} \|f_j\|_{L^\infty}^2)^{1/2},
 \]
 for every $\ep>0$, which need not hold in general. Thus in any setting of interest, one aims to prove an $\ell^2$ decoupling inequality for $L^p$ for some range $2 \leq p \leq p_0$, where the critical exponent $p_0$ depends on the setting.

\subsection{Extension and restriction operators}
So far we have considered examples of decoupling only in the setting of finite exponential sums. To introduce decoupling inequalities in the form considered in the work of Bourgain, Demeter and Guth, we require one more motivating principle, that of restriction problems, and the corresponding restriction and extension operators, which are intricately related to properties of Fourier transforms. 

The Fourier transform operator is in many ways ideally suited to working with functions in $L^2(\R^n)$, as it is a (surjective)  unitary operator on this Hilbert space. But there is a way in which an inconvenience arises, due to the general principle that a ``function'' $f$ in an $L^p$ space is actually an equivalence class of functions, so that $f$ may take completely arbitrary values on any set of measure zero. 
Consider instead that for a given function $f \in L^1(\R^n)$, its Fourier transform 
\[ \hat{f}(\xi) = \int_{\R^n} f(x) e^{-2\pi i x\cdot \xi}dx\]
is in fact a bounded, continuous function, which moreover (by the Riemann-Lebesgue lemma) vanishes at infinity (see e.g. \cite{SteinWeiss}). In this way, $L^1$ is more convenient than $L^2$, since given an arbitrary function $f \in L^2(\R^n)$, its Fourier transform  $\hat{f}$ can be any  $L^2$ function, so that there is in general no meaningful way to  consider the values of $\hat{f}$ on any set of measure zero.  More generally, in  intermediate spaces between $L^1$ and $L^2$, the Hausdorff-Young inequality states that if $f \in L^p(\R^n)$ with $1 \leq p \leq 2$ then $ \hat{f} \in L^{p'}(\R^n)$ with conjugate exponent $p'$ defined by $1/p + 1/p'=1$, so that \emph{a priori} we may only think of $\hat{f}$ as being defined almost everywhere, that is, on sets of positive measure. 

A remarkable insight of Stein in the 1970's is that if within $\R^n$ ($n \geq 2$) one considers not an arbitrary set of measure zero, but a smooth submanifold $S$ of appropriate curvature, then there is an exponent $p_0$ (depending on $S$) so that for every function $f \in L^p(\R^n)$ with $1 \leq p < p_0$, the Fourier transform $\hat{f}$ can be meaningfully restricted to $S$. With this, Stein \cite{Stein79} initiated a new field, that of \emph{restriction estimates}.

To be more precise, following the notation of \cite[Ch. VIII \S4]{SteinHA}, letting $d\sig$ be the induced Lebesgue measure on a compact smooth submanifold $S \subset \R^n$, we say that $(L^p,L^q)$ restriction holds for $S$ if
\beq\label{restriction_basic_ineq}
  \| \left. \hat{f} \right|_S \|_{L^q(S,d\sig)} = ( \int_{S} |\hat{f}(\xi)|^q d\sig (\xi))^{1/q} \leq A_{p,q}(S) \|f\|_{L^p(\R^n)} 
 \eeq
holds for every Schwartz function $f$ on $\R^n$.\footnote{We recall the Schwartz functions $\Scal(\R^n)$ are functions that are infinitely differentiable and have rapid decay, that is, the functions (along with all of their derivatives) decay as rapidly as $O((1+|x|)^{-M})$ for every $M \geq 1$, as $|x| \maps \infty$.  These very nicely behaved functions are particularly useful because the space $\Scal (\R^n)$ is dense in $L^p(\R^n)$ for $1 \leq p < \infty$ and is acted upon by the Fourier transform as a linear isomorphism.}
Once the inequality (\ref{restriction_basic_ineq}) is obtained for Schwartz functions,  one then uses the fact that they are dense in $L^p$ to see that for any $f \in L^p(\R^n)$, $\left. \hat{f} \right|_S$ is well-defined as an $L^q(S,d\sig)$ function; that is to say, in particular, that $\hat{f}$ is defined (almost everywhere, with respect to $d\sig$) on $S$. 

For which surfaces $S$, and which exponents $p$ and $q$ does such an inequality (\ref{restriction_basic_ineq}) hold? 
Our previous observations (Fourier transforms of $L^1$ functions are extremely nice on measure zero sets; Fourier transforms of $L^2$ functions are not) show that $(L^1,L^q)$ restriction holds trivially for all $1 \leq q \leq \infty$, while $(L^2,L^q)$ restriction fails for all $1 \leq q \leq \infty$.
The general expectation is that as long as $S$ has sufficient \emph{curvature}, and $1 \leq p<2$ is sufficiently close to 1, then an inequality of the form (\ref{restriction_basic_ineq}) should hold for an appropriate range of $q$ depending on $p$ and $S$.
 (The curvature of $S$ is critical, as it ensures that $S$ is not contained in any hyperplane; there exist functions $f$ which belong to an $L^p$ space but have Fourier transform unbounded along a hyperplane, see e.g. \cite{Tao04}.)
  Many questions of this type remain mysterious still today, in particular relating to how far $p$ can be pushed above 1 and still obey an $(L^p,L^q)$ restriction result.

There is an adjoint form of the restriction problem, the \emph{extension} problem, which we will use to state decoupling theorems. Above we defined the  \emph{restriction operator} $R_S$ acting on $f$ (initially  a Schwartz function on $\R^n$) by 
\[ R_S f(\xi) = \int_{\R^n}  f(x) e^{-2\pi i x \cdot \xi} dx , \]
for $\xi \in S$, so that $R_Sf  = \left. \hat{f} \right|_S$.
We now define an adjoint \emph{extension operator} $E_S$ acting on an integrable function $g$ on $S$ by 
\beq\label{extension_operator_dfn}
 E_S g(x) = \int_S g(\xi) e^{2\pi i x \cdot \xi} d\sig(\xi) ;
 \eeq
we see this is the inverse Fourier transform $(gd \sig)\check{\;} (x)$ along $S$. 
The dual inequality to (\ref{restriction_basic_ineq}) is then a putative extension estimate 
\beq\label{extension_basic_ineq}
\| E_S  g \|_{L^{p'}(\R^n)} \leq A_{p,q}(S) \|g\|_{L^{q'}(S,d\sig)},
\eeq
where $1/p+1/p'=1$ and $1/q+1/q'=1$. 
The sense in which (\ref{restriction_basic_ineq}) and (\ref{extension_basic_ineq}) are dual arises from the Plancherel identity, 
\[ 
\int_{\R^n} (gd\sig)\check{\;}(x) \overline{f(x)} dx = \int_{S} g(\xi) \overline{\hat{f}(\xi)} d\sig(\xi).
\]
As a result, proving the restriction inequality (\ref{restriction_basic_ineq}) for $(L^p,L^q)$ is equivalent to proving the extension inequality (\ref{extension_basic_ineq}) for $(L^{q'},L^{p'})$. The main conjecture in this area, called the restriction conjecture, may be therefore be framed in terms of either operator; in preparation for decoupling estimates, we use the adjoint extension form here:

\begin{conj}[Restriction conjecture, adjoint form]\label{conj_restriction_paraboloid}
Let $S$ be a compact $C^2$ hypersurface in $\R^n$ with nonvanishing Gaussian curvature at every point and surface measure $d \sig$. Then for $p' > \frac{2n}{n-1}$ and $q \leq  (\frac{n-1}{n+1})p'$,
\[ \| (g d\sig)\check{\;} \|_{L^{p'}(\R^n)} \leq A_{p,q}(S) \|g\|_{L^{q'}(S,d\sig)} \]
holds for all $g \in L^{q'}(S,d\sig)$.
\end{conj}
This applies, for example, to a compact piece of the paraboloid $\{ (\xi,|\xi|^2):  \xi \in  [-1,1]^{n-1}\}$, or a ``cap'' on the unit sphere $S^{n-1}$.
This conjecture highlights the quantity $\frac{2n}{n-1}$, which we will call the \emph{restriction exponent}.
Conjecture \ref{conj_restriction_paraboloid} is known in full for the case $n=2$; see \cite{Fef70} and also \cite{Zyg71}, or \cite[Ch. IX \S 5.1 and p. 432]{SteinHA} for a full exposition.
For the moment, for $n \geq 3$ we will only mention the early breakthrough work of  Tomas \cite{Tom75} and Stein (see e.g. \cite[Ch. IX \S 2.1]{SteinHA} for a more general statement):
 \begin{theo}[Tomas-Stein Restriction Theorem]\label{thm_SteinTomas}
  Let $S$ be a compact $C^2$ hypersurface in $\R^n$ with nonvanishing Gaussian curvature at every point and  surface measure $d\sig$. Then for $p' \geq \frac{2(n+1)}{n-1}$ and every  $g \in L^2(S,d\sig),$
  \[ \| (gd\sig)\check{\;} \|_{L^{p'}(\R^n)} \leq A \|g\|_{L^{2}(S,d\sig)}.\]
 \end{theo}
  The value $\frac{2(n+1)}{n-1}$ is called the \emph{Tomas-Stein exponent}, and we note that it is larger than the conjectured restriction exponent in Conjecture \ref{conj_restriction_paraboloid}. 
  
  Restriction estimates are an active field of current interest within harmonic analysis; for an evolving story of progress in dimensions $n \geq 3$ see e.g. \cite[Figure 1]{Tao04}. \footnote{Very recently a new world record was set for the (linear) restriction conjecture in the case of truncated cones, using polynomial partitioning \cite{OuWan17x}.}
See e.g. \cite{Tao04}, \cite{SteinHA} for  further considerations, such as restriction to the cone or submanifolds of lower dimensions, examples leading to necessary conditions, and connections to other areas such as PDE's.
Our motivation for introducing restriction estimates is to set the stage for the decoupling inequalities we will turn to now.

\section{Decoupling for the moment curve and how it proves Vinogradov}\label{sec_decoupling_VMVT}
\subsection{Statement of sharp $\ell^2$ decoupling for the moment curve}
We are ready to state the decoupling inequality that implies the Main Conjecture in the Vinogradov Mean Value Method: $\ell^2$ decoupling for $L^{n(n+1)}$ for the moment curve in $\R^n$. For $n \geq 2$ and any interval $J \subseteq [0,1]$, define 
\[ \Gamma_J = \{ (t,t^2, t^3, \cdots, t^n): t\in J \}.\]
Given an integrable function $g : [0,1] \maps \C$ and an interval $J \subseteq [0,1]$, define the operator
\beq\label{extn_Ga_op}
 E_J g(x_1,\ldots, x_n) = \int_J g(t) e(tx_1 + t^2 x_2 + \cdots + t^n x_n) dt.
 \eeq
We may recognize this now as the extension operator $(g d\sig)\check{\;}$ for the curve $\Gamma_J \subset \R^n$.

Recalling the notion of an $\ell^2$ decoupling inequality from (\ref{dfn_generic_decoupling}), we will study the decoupling of $E_{[0,1]}g$ into pieces $E_Jg$ as $J$ ranges over a partition of $[0,1]$ into (disjoint) sub-intervals of length $\del$, for  small $\del>0$. We will denote a sum over such a partition by $\sum_{\bstack{J \subset [0,1]}{|J|=\del}}$.\footnote{We will always assume this is well-defined, e.g. when $\del^{-1}$ is an integer, in which case $J$ ranges over the partition $[j\del, (j+1)\del]$ with $0 \leq j \leq \del^{-1}-1$.}

For technical reasons, we will measure the $L^p$ norm of any function $f : \R^n \maps \C$ according to a positive weight function $v: \R^n \maps \R_{>0}$, via the weighted norm
\beq\label{decoupling_weighted_norm}
 \|f\|_{L^p(v)} = ( \int_{\R^n} |f(x)|^p v(x) dx ) ^{1/p}.
 \eeq
In particular, we will apply this with weights of the following form: for a ball $B$ in $\R^n$ centered at a point $x_0$ and with radius $R$, we define the weight
\beq\label{decoupling_weight_dfn}
 w_B(x) = \frac{1}{ ( 1 + \frac{|x-x_0|}{R})^{E}},
 \eeq
 where $E$ is a sufficiently large positive integer (say $E \geq 100n$).\footnote{We think of this weight as a smoother version of the indicator function $\onebf_B$. The exponent $E$ is simply chosen large enough that we may apply e.g. H\"{o}lder's inequality as many times as the argument requires, and still yield a weight, which we will also call $w_B$, that is tailored to the ball $B$ and decays sufficiently to be integrable in $\R^n$.
The rigorous argument to prove decoupling theorems must in fact pass via an inductive hypothesis that assumes every intervening proposition is known relative to weights $w_{B,E}$ defined as in (\ref{decoupling_weight_dfn}), for every $E \geq 100n$. While this is important for the rigorous induction, it would affect our discussion below mainly in terms of notation, and we simply continue to write $w_B$ in each instance, without specifying the large value $E$.}

The landmark result of Bourgain, Demeter and Guth \cite{BDG16} is:
\begin{theo}[sharp $\ell^2$ decoupling for the moment curve]\label{thm_decoupling_moment_curve}
Let $n \geq 2$. For every $\ep>0$ there exists a constant $C_\ep = C(\ep,n)$ such that for every $0<\del \leq 1$,  for each ball $B \subset \R^n$ with radius $\del^{-n}$, and for every integrable function  $g : [0,1] \maps \C$, 
\beq\label{thm_decoupling_moment_curve_ineq}
 \| E_{[0,1]} g \|_{L^{n(n+1)}(w_B)} \leq C_\ep \del^{-\ep} ( \sum_{\bstack{J\subset [0,1]}{|J| = \del}} \|E_J g\|^2_{L^{n(n+1)}(w_B)} )^{1/2}.
 \eeq
\end{theo}
Importantly, the constant $C_\ep$ is independent of $\del$, the ball $B$, and the function $g$. To establish terminology, we think of this as the claim that we can \emph{decouple in $L^{n(n+1)}(\R^n)$ down to the scale} $\del$ (on the Fourier side), by detecting cancellation on  balls of radius $\del^{-n}$ (on the spatial side). There is a tension here between the scale down to which we can decouple, and the size of the ball  detected by the weight $w_B$: the larger the ball, the easier it is to detect cancellation on the spatial side. (One may deduce as a corollary to Theorem \ref{thm_decoupling_moment_curve} that the result continues to hold with $B$ replaced by any ball of radius larger than $\del^{-n}$; see Corollary \ref{thm_decoupling_moment_curve_cor}.) The relationship evident in Theorem \ref{thm_decoupling_moment_curve}, which also plays a role throughout the proof, is that  generally speaking, the smaller the scale to which we wish to decouple, the larger the ball must be.

\subsubsection{Special cases}
Incidentally, we could have stated Theorem \ref{thm_decoupling_moment_curve} to include dimension $n=1$. In dimension $n=1$ (in contrast to $n \geq 2$), the result may be proved very simply---in this case it is a claim about $\ell^2$ decoupling for $L^2(\R)$, claiming decoupling down to the scale $\del$ by detecting cancellation on balls of radius $\del^{-1}$, that is, of radius equal to the \emph{inverse} of the decoupling scale. Importantly, we will  see (Proposition \ref{prop_ell2_L2_linear}) that this strong $\ell^2$ decoupling result for $L^2$, which decouples down to scale $\del$ relative to balls of radius $\del^{-1}$ (as opposed to $\del^{-n}$) holds not only on $\R$ but on $\R^n$ for all $n \geq 1$; this depends on the special relationship between $\ell^2$ and $L^2$, and uses almost orthogonality.

Note also that the case of Theorem \ref{thm_decoupling_moment_curve}  with $n=2$ is special, since then the moment curve $\{ (t,t^2) \}$ is the 1-dimensional parabola in $\R^2$.
In this case the result of Theorem \ref{thm_decoupling_moment_curve} was already known due to sharp $\ell^2$ decoupling inequalities proved by Bourgain and Demeter for paraboloids $\{ (t,|t|^2)\}$ in $\R^n$ for all $n \geq 2$ \cite{BouDem15}.
The new achievement of \cite{BDG16} is proving Theorem \ref{thm_decoupling_moment_curve} for $n \geq 3$. 

\subsubsection{Previous results}
While we will focus here nearly exclusively on the paper \cite{BDG16} of Bourgain, Demeter and Guth, the sequence of breakthroughs on decoupling which led to this fantastic result  originated in earlier work of Bourgain \cite{Bou13}, followed by the Bourgain-Demeter proof of the $\ell^2$ decoupling conjecture for paraboloids in \cite{BouDem15}, and continuing with many other works (to which we briefly return in  \S \ref{sec_decoupling_other}). 
In particular, we note that \cite{BDG16} builds on an earlier approach of \cite[Thm. 1.4]{BouDem14x} which proved $\ell^2$ decoupling for $L^p(\R^n)$ for nondegenerate curves in $\R^n$, $n \geq 2$ in a more limited range $2 \leq p \leq 4n-1$ instead of up to the best possible critical exponent $p_n = n(n+1)$ (see \S \ref{sec_an_heur}). This consequently proved only weak results toward the Main Conjecture in the setting of the Vinogradov Mean Value Method, but certainly demonstrated the applicability of the decoupling method to this arithmetic question.

\subsection{From decoupling to discrete decoupling}	 
	Before we discuss the proof of Theorem \ref{thm_decoupling_moment_curve}, we show how it implies  Theorem \ref{thm_VMVT}, proving the Main Conjecture in the Vinogradov Mean Value Method for the remaining cases of degree at least $4$.
	The first step is to pass from Theorem  \ref{thm_decoupling_moment_curve} to  a discrete version, following \cite[Thm. 4.1]{BDG16}:
	
\begin{theo}[discrete $\ell^2$ decoupling for the moment curve]\label{thm_decoupling_moment_curve_discrete}
Let $n \geq 2$ and $p \geq 2$ be fixed. For every $\ep>0$ there exists a constant $C_\ep = C(\ep,n,p)$ such that the following holds: for every $N \geq 1$,  and for each choice of a fixed set of points $\{t_1,\ldots, t_N\}$ with  $t_i \in (\frac{i-1}{N}, \frac{i}{N}]$,  for each ball $B_R$ of radius $R \geq N^n$ in $\R^n$, and every set of coefficients $\{a_i\}_{1\leq i \leq N}$ with $a_i \in \C$,
\begin{multline*}
 ( \frac{1}{|B_R|} \int_{\R^n} | \sum_{i=1}^N a_i e(t_i x_1 + t_i^2x_2 + \cdots + t_i^n x_n)|^p w_{B_R}(x) dx_1 \cdots dx_n)^{1/p} \\
 \leq C_\ep N^\ep\{1 + N^{\frac{1}{2}(1 - \frac{n(n+1)}{p}) }\}( \sum_{i=1}^N |a_i|^2)^{1/2}.
\end{multline*}
\end{theo}
Note that the constant $C_\ep$ is independent of $N, R$, the points $t_i$, and the coefficients $a_i$.
It is an important feature of this theorem that although it concerns exponential sums, no (arithmetic) properties are assumed of the frequencies $t_i$, except that they fall into distinct intervals $((i-1)/N,i/N]$, that is, are $1/N$-separated, and so cannot cluster together.

To prove Theorem \ref{thm_decoupling_moment_curve_discrete}, it suffices to prove the inequality 
  in the critical case $p_n=n(n+1)$, in which case the right-hand side of the decoupling inequality is $2C_\ep N^\ep (\sum_i |a_i|^2)^{1/2}$. (The sufficiency follows from an argument entirely analogous to \S \ref{sec_VMVT_critical_case}; note that for reasons of tradition, the \emph{critical index} in the arithmetic setting is $s_n=\frac{1}{2}n(n+1)$ while the \emph{critical exponent} in the real-variable setting is considered to be $p_n=2s_n=n(n+1)$, but the meaning is analogous.)

Next, we use the following immediate consequence of Theorem \ref{thm_decoupling_moment_curve}, which extends a decoupling inequality known for balls of a fixed radius to balls of all larger radii  (e.g.  the proof of \cite[Thm. 4.1]{BDG16}):
\begin{coro}\label{thm_decoupling_moment_curve_cor}
The conclusion of Theorem \ref{thm_decoupling_moment_curve} holds for any ball of radius $R \gg \del^{-n}$ with an implied constant independent of $R$.
\end{coro}

Arguing formally, we may apply Corollary \ref{thm_decoupling_moment_curve_cor} to the function $g :[0,1]\maps \C$ chosen to be the finite linear combination $\sum_{i=1}^N a_i \del_{t=t_i}$ of Dirac delta functions, and with the choice $\del = 1/N$. In this case,
\[ E_{[0,1]}g(x_1,\ldots, x_n) = \sum_{i=1}^N a_i e( t_ix_1+\cdots +t_i^n x_n).\]
On the other hand, for each interval $J = (\frac{i-1}{N},\frac{i}{N}] \subset [0,1]$, 
\[ E_{J}g(x_1,\ldots, x_n) =a_i e( t_ix_1+\cdots +t_i^n x_n) \]
for the unique $t_i \in J$, 
so that for $p=n(n+1)$,
\beq\label{E_J_upper_bound}
\| E_J g \|_{L^p(w_{B_R})} = ( \int_{\R^n} |a_i|^{p} w_{B_R}(x) dx)^{1/p} \ll_p |B_R|^{1/p} |a_i| .
\eeq
 Thus applying Corollary \ref{thm_decoupling_moment_curve_cor} with the ball $B_R$, we may deduce that
\[ 
 ( \frac{1}{|B_R|} \int_{\R^n} | \sum_{i=1}^N a_i e(t_i x_1 + t_i^2x_2 + \cdots + t_i^n x_n)|^{p} w_{B_R}(x) dx_1 \cdots dx_n)^{1/p} 
 \leq C_\ep N^\ep (\sum_i|a_i|^2)^{1/2},
\]
for some appropriate $C_\ep$,
as claimed. More rigorously, in place of delta functions one may apply Corollary  \ref{thm_decoupling_moment_curve_cor} to the integrable function
\[ g_\mu(t) = \sum_{i=1}^N a_i K_\mu (t_i-t),\]
where  the family  $\{K_\mu\}_{\mu>0}$ is a so-called approximation to the identity (e.g. $K_\mu(t) = \frac{1}{2\mu} \onebf_{(-\mu,\mu)}(t) $ or a smooth approximation of this). Then the above result continues to hold via an appropriate limiting argument, as $\mu \maps 0$.

\subsection{From discrete decoupling to Vinogradov}\label{sec_dis_Vin}
We now focus on how to deduce Theorem \ref{thm_VMVT} from the discrete decoupling inequality, thus proving the Main Conjecture in Vinogradov's Mean Value Method. In fact,  while Theorem \ref{thm_VMVT} considers the classical scenario of counting integral solutions $(x_1,\ldots, x_{2s})$ to the Vinogradov system (\ref{Vin_sys_dfn}), we may now extend our consideration to counting solutions $(z_1,\ldots, z_{2s})$ that take values in a fixed collection of real numbers that are instead merely 1-separated, that is, a  fixed collection of real numbers $\Mcal=\{\xi_i\}_{i \geq 1}$ such that for every $i \geq 1$, $i-1 < \xi_i \leq i$. 
With respect to such a set, we denote by $J_{s,k,\Mcal}(N)$ the number of solutions $(z_1,\ldots, z_{2s})$ to the inequalities (analogous to the Vinogradov system)
\beq\label{VMVT_decoupling_inequalities}
 |x_1^j + \cdots + x_s^j -( x_{s+1}^j + \cdots + x_{2s}^j) |\leq \frac{1}{2}\frac{1}{N^{k-j}}, 	\quad 1 \leq j \leq k ,
 \eeq
with all the $z_i \in \Mcal \intersect (0,N]$. (Note that by construction $\Mcal \intersect (0,N]$ contains precisely $N$ points, for $N \geq 1$.)

\begin{theo}[Vinogradov Mean Value Theorem, generalized version]\label{thm_VMVT_generalized}
For all integers $s, k \geq 1$ and every fixed collection $\Mcal$ of 1-separated real numbers  as above,
\[ 
J_{s,k,\Mcal}(N) \ll_{s,k,\ep}  N^\ep \{ N^s + N^{2s - \frac{1}{2} k(k+1)} \},\]
for every $N \geq 1$, and all $\ep>0$, with an implied constant  independent of the set $\Mcal$.
\end{theo}
 Theorem \ref{thm_VMVT} follows immediately from this more general version: if $\Mcal$ is the set of positive integers, for each $N \geq 1$ the inequalities (\ref{VMVT_decoupling_inequalities}) specify that the left-hand side is an integer with absolute value $< 1$, in which case it must be the integer zero.

By the argument in \S \ref{sec_VMVT_critical_case} it suffices, for each $k \geq 1$, to prove Theorem \ref{thm_VMVT_generalized} in the critical case $s=s_k = \frac{1}{2}k(k+1)$.
We first note that the case $k=1$ may be proved trivially, since it is counting solutions to $|x_1 -x_2| \leq \frac{1}{2}$ with $x_1,x_2  \in \Mcal \intersect (0,N]$; given $x_1 \in \Mcal,$ there are at most $2$ choices of $x_2 \in \Mcal$ that satisfy this inequality, from which we deduce that $J_{1,1,\Mcal}(N) \ll N$, as desired. Thus we restrict our attention to $k \geq 2$, and  the moment curve $(t,t^2,\ldots, t^k) \subset \R^k$.

The first step is to replace the generic weight function $w_B$ appearing in Theorem \ref{thm_decoupling_moment_curve_discrete} by an advantageously chosen weight. We  fix a Schwartz class function $\phi : \R^k \maps [0,\infty)$ that satisfies $\phi(x) >0$ for all $x \in \R^k$ and also has positive Fourier transform $\hat{\phi}(\xi)>0$ for all $\xi \in \R^k$, with the further requirement that $\hat{\phi}(\xi) \geq 1$ for $|\xi| \leq 1$. (That such a function  exists may be observed by recalling that the Gaussian $e^{-\pi |x|^2}$ has Fourier transform $e^{-\pi |\xi|^2}$, so that appropriate scalings of these functions provide options for $\phi$.)
We now rescale $\phi$ by defining for any $M>0$, $\phi_M(x) = \phi(x/M)$; we note for later reference that by the standard rules for rescaling Fourier transforms,
\beq\label{rescaling_principle_0}
(\phi_M)\hat{\;}(x) = M^k \hat{\phi} (Mx).
\eeq
 Because any Schwartz function $\phi$ has  rapid decay, given $N\geq 1$, there exists a ball $B_R$ of radius $N^k \ll R \ll N^k$ such that $\phi_{N^k}(x)$ is majorized by $w_{B_R}$, so that we may apply Theorem \ref{thm_decoupling_moment_curve_discrete}, with $n=k$, $p=2s=2s_k$ and all the coefficients $a_i=1$ to conclude that for any choice of points $t_i \in (\frac{i-1}{N},\frac{i}{N}]$,
\beq\label{VMVT_decoupling_app_ineq}
N^{-k^2} \int_{\R^k} | \sum_{i=1}^N e(t_i x_1 + t_i^2x_2 + \cdots + t_i^k x_k)|^{2s} \phi_{N^k}(x) dx_1 \cdots dx_k\\
\ll_{k,\ep} N^{s+\ep}.
\eeq
In particular, we may choose $t_i = \xi_i/N$, where $\xi_i \in \Mcal$, in which case we see that by a change of variables,
\begin{multline*}
N^{-k^2}  \int_{\R^k} | \sum_{i=1}^N e((\frac{\xi_i}{N}) x_1 + (\frac{\xi_i}{N})^2x_2 + \cdots + (\frac{\xi_i}{N})^k x_k)|^{2s} \phi_{N^k}(x_1,\ldots,x_k) dx_1 \cdots dx_k
\\= N^{\frac{1}{2}k(k+1)-k^2}\int_{\R^k} | \sum_{i=1}^N e(\xi_i x_1 + \xi_i^2 x_2+ \cdots + \xi_i^k x_k)|^{2s} \phi_{N^k}(Nx_1,\ldots,N^kx_k) dx_1 \cdots dx_k.
 \end{multline*}
 Here it is convenient to note that 
 \[ \phi_{N^k} (Nx_1, N^2x_2,\ldots,N^kx_k) =\phi( \frac{x_1}{N^{k-1}},\frac{x_2}{N^{k-2}}, \ldots, x_k).\]
Now we expand the $2s$-moment  to see that the integral expression above is precisely
\beq\label{VMVT_deduction_identity}
N^{\frac{1}{2}k(k+1)-k^2} \sum_{\bstack{\xi = (\xi_1,\ldots, \xi_{2s}) }{\xi_i \in \Mcal \intersect (0,N]}} \int_{\R^k} e(\theta_1(\xi) x_1 + \theta_2(\xi) x_2+ \cdots + \theta_k(\xi) x_k) \phi( \frac{x_1}{N^{k-1}}, \ldots, x_k) dx_1 \cdots dx_k,
\eeq
where for each $1 \leq j \leq k$ we use $\theta_j(\xi)$ to denote
\[ \theta_j(\xi) = \xi_1^j + \cdots+ \xi_s^j  - (\xi_{s+1}^j + \cdots+ \xi_{2s}^j) .\]
Observe that for each fixed $\xi$, the integral  in (\ref{VMVT_deduction_identity}) is the Fourier transform of the function $\phi( \frac{x_1}{N^{k-1}}, \ldots, x_k),$ evaluated at the point $(\theta_1(\xi), \ldots, \theta_k(\xi))$. 
According to the standard rescaling rule for Fourier transforms (a generalization of (\ref{rescaling_principle_0})), 
\[  \phi ( \frac{\cdot}{N^{k-1}}, \frac{\cdot}{N^{k-2}},\ldots, \frac{\cdot}{1} ) \hat{\;}(\eta_1,\ldots, \eta_k) =  
N^{k-1} \cdots N^{k-2} \cdots 1 \cdot \hat{\phi} (N^{k-1}\eta_1,N^{k-2}\eta_2,\ldots, \eta_k).
\]
We may re-write the leading factor as $N^{k^2 - \frac{1}{2}k(k+1)}$ and plug this identity for the Fourier transform into (\ref{VMVT_deduction_identity}), to see that the full expression in (\ref{VMVT_deduction_identity}) becomes
\beq\label{VMVT_deduction_lower_bound}
 \sum_{\bstack{\xi = (\xi_1,\ldots, \xi_{2s}) }{\xi_i \in \Mcal \intersect (0,N]}} 
 \hat{\phi} (N^{k-1}\theta_1(\xi),N^{k-2}\theta_2(\xi),\ldots, \theta_k(\xi)).
\eeq
Now we use the critical aspect of our choice for $\phi$: recall that $\hat{\phi}$ is always positive, and moreover that $\hat{\phi}(\eta_1,\ldots, \eta_k) \geq 1$ for $|\eta| \leq 1$ (and hence certainly for $\eta$ with $|\eta_j| \leq 1/2$ for each $1 \leq j \leq k$). Thus in particular, every term in (\ref{VMVT_deduction_lower_bound}) is positive, and we may bound the sum in  (\ref{VMVT_deduction_lower_bound}) from below by 
the number of tuples $\xi = (\xi_1, \ldots, \xi_{2s})$ of points in $\Mcal \intersect (0,N]$ 
such that 
\[ N^{k-j}|\theta_j(\xi)| \leq 1/2 \qquad  \text{for each $1 \leq j \leq k.$}\]
This is equivalent to specifying that $\xi$ is a solution to the system of inequalities (\ref{VMVT_decoupling_inequalities}). 
We combine this observation with the key decoupling inequality (\ref{VMVT_decoupling_app_ineq}) to conclude that 
\[  J_{s,k,\Mcal}(N) \leq \sum_{\bstack{\xi = (\xi_1,\ldots, \xi_{2s}) }{\xi_i \in \Mcal \intersect (0,N]}} 
 \hat{\phi} (N^{k-1}\theta_1(\xi),N^{k-2}\theta_2(\xi),\ldots, \theta_k(\xi)) 
 	\ll_{k,\ep} N^{s+\ep} ,
	\]
	in the case $s=s_k$. This proves Theorem \ref{thm_VMVT_generalized} and hence the Main Conjecture in Vinogradov's Mean Value Method.

We will shortly turn to the proof of decoupling for the moment curve, which will occupy our focus for the remainder of the manuscript, but we pause to comment on the critical exponent $p_n = n(n+1)$ and to state in more generality the type of decoupling theorem that can now be obtained in several important cases.

\subsection{Limiting cases and critical exponents}\label{sec_an_heur}
We have already seen that in general we would not expect an $\ell^2$ decoupling result for $L^p$ to hold for $p<2$. In the setting of Theorem \ref{thm_decoupling_moment_curve} we can also observe, at least heuristically, that the statement of the theorem is best possible (perhaps up to the factor $\del^{-\ep}$) in the sense that we would not expect an $\ell^2$ decoupling result to hold for $L^p(\R^{n})$ if $p > n(n+1)$. This can be seen, roughly speaking, from an example like the following (see Tao \cite{Tao15ablog}): if an interval $J$ has length $\del \leq 1$, then the arc $\Gamma_J = \{ (t,t^2,\ldots, t^n) : t \in J\}$ fits inside a box of proportions $\del \times \del^2 \times \cdots \times \del^n$. Suppose we construct a particular function $g$ as a sum $\sum_J g_J$ of smooth bump functions with $g_J$ supported in $J$, where $J$ runs over a dissection of $[0,1]$ with $|J|=\del$. Then computations related to the uncertainty principle for Fourier transforms (see \S \ref{sec_uncertainty} for related notions) imply that for each $J$, $|E_Jg| = |E_Jg_J|$ is essentially  constant $\approx \del$ on a ``dual'' box of proportions $\del^{-1} \times \del^{-2} \times \cdots \times \del^{-n}$, and decays outside of the box; the dual box fits inside a ball of radius $\del^{-n}$, and has volume $\del^{-n(n+1)/2}$.
Thus heuristically for each such $J$, 
\[
 \|E_J g\|_{L^{n(n+1)}(w_B)} = \|E_J g_J\|_{L^{n(n+1)}(w_B)}
	\approx \del \cdot (\del^{-n(n+1)/2})^{\frac{1}{n(n+1)}} \approx \del^{1/2}.
\]
	Summing over the $O(\del^{-1})$ many intervals $J$, in this example the right-hand side in Theorem \ref{thm_decoupling_moment_curve} is comparable to $\del^{-\ep}(\sum_J (\del^{1/2})^2)^{1/2}  = O(\del^{-\ep})$ for arbitrarily small $\ep>0$.
Regarding the left-hand side of Theorem \ref{thm_decoupling_moment_curve}, similar considerations indicate that $E_{[0,1]}g$ is comparable to 1 on the unit ball centered at the origin in $\R^n$, so that its $L^{n(n+1)}(w_B)$ norm is $\gg 1$. Now putatively if we considered larger $L^p$ norms, say $p= n(n+1) +\kappa_0$ for some fixed $\kappa_0>0$, we would instead see a strictly positive power of $\del$ on the right-hand side of Theorem \ref{thm_decoupling_moment_curve} (while the left-hand side stayed unchanged), so that as $\del \maps 0$ a decoupling inequality could not hold.
From this example, we also gain the intuition that we would not in general expect to be able to prove $\ell^2$ decoupling for the moment curve down to the scale $\del$, at the critical exponent $n(n+1)$, if we used balls of radius significantly smaller than $\del^{-n}$: this is the smallest radius that can capture the full proportions $\del^{-1} \times \del^{-2} \times \cdots \times \del^{-n}$ of the dual box.\footnote{There are however open questions about proving $\ell^2$ decoupling for $L^p$ for $p< n(n+1)$, or $\ell^p$ decoupling for $L^p$ for certain $p$, for spatial balls of smaller radii.}

Now we recall from \S  \ref{sec_Hua}  Bourgain's $\ell^2$ decoupling result for $L^{\ell(\ell+1)}(\R^\ell)$ for the curve $\Gamma_{\ell,k}$ defined in (\ref{Gamma_Hua}), from which he deduced an improvement of Hua's inequality and $\tilde{G}(k)$ in Waring's problem.
Here we note that for a fixed integer $s \geq 1$, the $2s$-moment associated to $\Gamma_{\ell,k}$, namely
\beq\label{Jsk'_int}
 \int_{(0,1]^{\ell}} | \sum_{1 \leq x \leq X} e(\al_1 x + \al_2 x^2 + \cdots + \al_{\ell-1}x^{\ell-1} + \al_\ell x^k)|^{2s} d\albf 
 \eeq
counts solutions $1 \leq x_i \leq X$ to the system of $\ell$ equations
\beq\label{sys_lk} 
x_1^j + \cdots + x_s^j = x_{s+1}^j + \cdots + x_{2s}^j, \quad \text{$j=1,2,\ldots,\ell-1,$ and $k$.}
\eeq 
This system has at least $\gg_s X^s$ diagonal solutions, and on the other hand (roughly) a contribution of size $ \gg_\ell X^{2s - \frac{1}{2}\ell(\ell-1) - k}$ to (\ref{Jsk'_int}) from $\albf$ near the origin, say with $|\al_j| \leq \frac{1}{8 \ell}X^{-j}$ for $1 \leq j \leq \ell-1$ and $|\al_\ell| \leq \frac{1}{8\ell }X^{-k}$. Thus heuristically, we would expect the critical index $s$ for the counting problem associated to this system of Diophantine equations to occur when these contributions are of equal order, namely $s=s_{\ell,k}= \frac{1}{2}\ell(\ell-1) + k$. This suggests that in the context of an $\ell^2$ decoupling result for $L^p$, the largest $p$ for which the decoupling result might be expected to hold would occur at $p=2s_{\ell,k} = \ell(\ell-1) + 2k$. 
 For example, for $\ell=k$, this would be $L^{\ell(\ell+1)}$, and indeed $\Gamma_{\ell,k}$ is the moment curve and (\ref{sys_lk}) is the Vinogradov system of degree $k$, so that Theorem \ref{thm_decoupling_moment_curve} applies directly. The new cases are when $1 \leq \ell  \leq k-1$ and the curve $\Gamma_{\ell,k}$ is not affine invariant; in these cases Bourgain's result of $\ell^2$ decoupling for $L^{\ell(\ell+1)}$  is for $L^p$ with $p$ strictly smaller than the putative critical exponent $\ell(\ell-1)+2k$. Thus, in this sense, the decoupling result currently known for $\Gamma_{\ell,k}$ is quantitatively weaker than that obtained for the affine invariant moment curve (even though it already has significant arithmetic applications).\footnote{This is related to open questions for curves that are not affine invariant, such as $(t,t^3) \subset \R^2$. Here, reasoning as above, we see the critical index for the system of equations  associated to the $2s$-moment occurs at $s=4$, so that the largest $p$ for which one might hope to prove an $\ell^2$ decoupling result for $L^p$ is $p=2s=8$. Demeter (personal communication) notes that there is a counterexample to $\ell^2$ decoupling for $L^p$ for $6<p \leq 8$ for this curve, and it is open what the best possible $\ell^r$ decoupling results for $L^p$, for appropriate $r,p$, should be in this case. See also Wooley's work \cite{Woo15a}, relevant to $p=8,9$.}

\subsection{Decoupling inequalities in other settings}\label{sec_decoupling_other}
We turn briefly to decoupling in a more general setting.
Let $\Mcal$ be a compact smooth manifold in $\R^n$ with an associated measure $\sig$ on $\Mcal.$ Suppose that $\Mcal$ has been partitioned (or covered with finite overlap) by ``caps'' $\tau$ of size $\delta$ (where ``size'' is measured in a sense appropriate to $\Mcal$). 
Given an integrable function $g : \Mcal \maps \C$, for each cap $\tau$ let $g_\tau$ denote the restriction $g \onebf_\tau$ of $g$ to $\tau$.
In particular, $g = \sum_{\tau} g_\tau$. 
In this setting, an $\ell^2$ decoupling result for $L^p$ (in terms of extension operators on $\Mcal$) states that there exists a critical index $p_c>2$ and some $\kappa \geq 2$ (with both $p_c$ and $\kappa$ depending on $\Mcal$) such that 
\beq\label{decoupling_general_form}
 \| (g d\sig)\check{\;} \|_{L^p(B)} \ll_\ep \del^{-\ep} ( \sum_{\tau} \| (g_\tau d\sig)\check{\;} \|^2_{L^p(B)})^{1/2} 
 \eeq
for each ball $B \subset \R^n$ with radius $\del^{-\kappa}$ and for each $2 \leq p \leq p_c$. 
Note that (\ref{decoupling_general_form}), without the $\del^{-\ep}$ loss, immediately holds for $p=2$, since the functions $g_\tau$ are (nearly) orthogonal.

\subsubsection{Hypersurfaces with positive definite second fundamental form}
The setting in which Bourgain and Demeter proved the first landmark sharp $\ell^2$ decoupling result is that of a compact $C^2$ hypersurface $S \subset \R^n$ with appropriate curvature, say $S = \{ (\xi, \ga(\xi)): \xi \in [0,1]^{n-1}\}$ for a certain defining function $\ga : \R^{n-1} \maps \R$. 
For any cube $Q \subset [0,1]^{n-1}$,  define the  operator acting on integrable functions $g: [0,1]^{n-1} \maps \C$ by
\[ E_Q g(x_1,\ldots, x_n)  = \int_Q g(\xi_1,\ldots, \xi_{n-1}) e(\xi_1 x_1 + \cdots + \xi_{n-1}x_{n-1} + \ga(\xi) x_n) d\xi.\]
We recognize this as the extension operator mapping $g$ to $(g d\sig ) \check{\;}$ for the manifold $S$, analogous to our general definition (\ref{extension_operator_dfn}).
To each cube $B \subseteq \R^n$ of center $x_0$ and side-length $\ell(B)=R$, we associate the weight $w_B$ as defined in (\ref{decoupling_weight_dfn}).
We also recall the weighted norm (\ref{decoupling_weighted_norm}). (Furthermore, to make everything well-defined, we work with dyadic cubes, i.e. with $\del \in 2^{-2\mathbb{N}}$  we may partition $[0,1]$ into intervals of length $\del^{1/2}$.)

The breakthrough result of Bourgain and Demeter \cite{BouDem15}, building on \cite{Bou13}, proves the following:
\begin{theo}[$\ell^2$ decoupling for hypersurfaces, extension form]\label{thm_decoupling_paraboloid}
Let $S$ be a compact $C^2$ hypersurface in $\R^n$ with positive definite second fundamental form.  For $2 \leq p \leq \frac{2(n+1)}{n-1}$, for every cube $B \subset \R^n$ of side-length $\ell(B)=\del^{-1}$, for every integrable function $g :[0,1]^{n-1} \maps \C$,
\[ \|E_{[0,1]^{n-1}} g\|_{L^p(w_{B})} \ll_{p,n,\ep} \del^{-\ep} ( \sum_{\bstack{Q \subset [0,1]^{n-1}}{\ell(Q) = \del^{1/2}}} \|E_Q g\|_{L^p(w_{B})}^2)^{1/2}.
\]
\end{theo}
The strength of this result lies in the fact that it extends all the way up to $p \leq 2(n+1)/(n-1)$, the Tomas-Stein exponent, while the initial progress in \cite{Bou13} treated $p \leq 2n/(n-1)$, the restriction exponent.\footnote{An analogous square function estimate is only conjectured to hold for $2 \leq p \leq 2n/(n-1)$, although as mentioned in \S \ref{sec_decoupling_square_function}, a square function estimate would be stronger in that range.}
Examples of  appropriate $S$ include  a compact piece of a paraboloid in $\R^n$, 
\beq\label{S_truncated_paraboloid}
 P^{n-1} = \{ (\xi,|\xi|^2) : \xi \in [0,1]^{n-1}\},
 \eeq
and  the sphere $S^{n-1} \subset \R^n$.\footnote{Here, to be consistent with the presentation in this manuscript, we have stated Theorem \ref{thm_decoupling_paraboloid} in the form of \cite[Thm. 1.1]{BouDem17a}, while in the original paper \cite[Thm. 1.1]{BouDem15} it is stated in terms of a Schwartz function $f$ on $\R^n$ whose Fourier transform is supported in a $\del$-neighborhood of the hypersurface. For example, in the case of $P^{n-1}$, such a function $f$ is  decoupled in terms of  a family $\{f_\theta\}$ constructed by dissecting $\hat{f}$ according to small $\del^{1/2} \times  \cdots \times \del^{1/2} \times \del$ ``slabs'' $\theta$ that form a finitely overlapping cover of the $\del$-neighborhood of $P^{n-1}$.}

As a consequence of Theorem \ref{thm_decoupling_paraboloid}, Bourgain and Demeter \cite[Thm. 1.2]{BouDem15} have also deduced  a sharp $\ell^2$ decoupling result for the truncated cone
\[ C^{n-1} = \{ (\xi_1,\ldots, \xi_{n-1}, (\xi_1^2 + \cdots + \xi_{n-1}^2)^{1/2}), 1 \leq (\sum \xi_j^2)^{1/2}\leq 2 \} \subset \R^n \]
(for which a nice application was anticipated in \cite{PraSee07}).
This is a return to the original setting of $\ell^p$ decoupling in Wolff \cite{Wol00}.
Bourgain and Demeter \cite{BouDem14x} have further proved sharp $\ell^p$ decoupling results (for certain $p$) for compact $C^2$ hypersurfaces with nonvanishing Gaussian curvature; since this includes hyperbolic paraboloids  (such as $\{\xi_1,\xi_2,\xi_3, \xi_1^2 + \xi_2^2 - \xi_3^2\} \subset \R^4$), which contain lines, the result  necessarily holds in a different range than Theorem \ref{thm_decoupling_paraboloid}. See also for example \cite{BouDem15x, DGS16x, Guo17x, BDG16x, BouDem16a, BouDem16b} for decoupling results for other submanifolds in various ambient dimensions.

\subsubsection{Discrete restriction phenomena and beyond}
We have already mentioned  variations on decoupling with number-theoretic consequences throughout \S \ref{sec_mot_cons}; here we mention a few more settings in which decoupling has opened new doors. 
Recall the Tomas-Stein Theorem \ref{thm_SteinTomas} on restriction for a $C^2$ compact hypersurface $S$ with nonvanishing Gaussian curvature. 
One can deduce the following discrete consequence from the statement of Theorem \ref{thm_SteinTomas} for each fixed $p$:
for every $0 < \del \leq 1$, every $\del^{1/2}$-separated set $\Lambda \subset S$, and every sequence $\{a_\xi \}_{\xi \in \Lambda}$ of complex numbers, 
\[ ( \frac{1}{|B|} \int_B | \sum_{\xi \in \Lambda} a_\xi e(\xi \cdot x)|^p dx)^{1/p} \ll_{p,n} \del^{\frac{n}{2p} - \frac{n-1}{4}} ( \sum_{\xi \in \Lambda} |a_\xi|^2)^{1/2},\]
for every ball $B \subset \R^n$ of radius $\del^{-1/2}$. 
We may think of this as capturing cancellation on the spatial scale reciprocal to the separation of the frequencies. With the new decoupling Theorem \ref{thm_decoupling_paraboloid} in hand, Bourgain and Demeter have improved this in \cite[Thm. 2.2]{BouDem15}, saving an additional $\del^{1/2p}$, by working with balls of larger radius $\del^{-1}$:

\begin{theo}\label{thm_discrete_restriction}
Let $S$ be a compact $C^2$ hypersurface in $\R^n$ with positive definite second fundamental form. Let $\Lambda \subset S$ be a $\del^{1/2}$-separated set, and let $R \gg \del^{-1}$. Then for $p \geq \frac{2(n+1)}{n-1}$, for every ball $B_R \subset \R^n$ of radius $R \gg \del^{-1}$ and for every $\ep>0$,
\beq\label{discrete_restriction_conj_ineq}
( \frac{1}{|B_R|} \int_{B_R} | \sum_{\xi \in \Lambda} a_\xi e(\xi \cdot x)|^p)^{1/p} \ll_{p,n,\ep} \del^{\frac{n+1}{2p} - \frac{n-1}{4} - \ep}  ( \sum_{\xi \in \Lambda} |a_\xi|^2)^{1/2}.
\eeq
\end{theo}
As a result, Bourgain and Demeter deduce for example a sharp \emph{discrete restriction estimate} for integer lattice points on  the paraboloid 
\[ P^{n-1}(N) = \{ (\xi_1,\ldots, \xi_{n-1}, \xi_1^2 + \cdots + \xi_{n-1}^2) : \xi_i \in \Z \intersect [-N,N] \},\]
of the form 
\[
  \| \sum_{\xi \in P^{n-1}(N)} a_\xi e(\xi \cdot x) \|_{L^p([0,1]^n)} \ll_{n,\ep} N^\ep ( \sum_{\xi \in  P^{n-1}(N)} |a_\xi|^2)^{1/2},
\]
for $p=\frac{2(n+1)}{n-1}$, for all $n \geq 4$.
This completes the proof of a conjecture in a program of discrete restriction phenomena, initiated in \cite{Bou93}. 
See also work of Bourgain and Demeter on discrete restriction on the sphere \cite{BouDem13a,BouDem15a} and \cite[Thm. 2.7]{BouDem15}. Discrete restriction problems are an interesting bridge between Fourier analytic questions and Diophantine problems, and it is worth noting that Wooley  \cite{Woo17a} has also successfully adapted the efficient congruencing method to the realm of discrete restriction problems, even though they are not strictly translation-dilation invariant.
Discrete restriction problems also relate to Strichartz estimates for Schr\"{o}dinger operators. In the periodic setting, decoupling has contributed to Strichartz estimates for Schr\"{o}dinger operators on both classical and irrational tori \cite[Thm. 2.4]{BouDem15}, on irrational tori \cite{DGG17x} for long time-intervals, and a bilinear Strichartz estimate on irrational tori \cite{FSWW16}, via bilinear decoupling methods.
The utility of decoupling also has been extended to Schr\"{o}dinger maximal estimates in non-periodic settings \cite{DGL16x}.

 A generalization  of Theorem \ref{thm_discrete_restriction} in \cite[Thm. 2.16]{BouDem15} allows one to consider, in effect, not just $\del^{1/2}$-separated points lying directly on a hypersurface $S$ but instead within a small neighborhood of $S$; this is in the vein of Diophantine questions on error terms for the number of integer lattice points within a region such as a circle or sphere. Other discrete applications of decoupling include improved results for certain well-known problems in incidence geometry; see \cite[\S 2.4]{BouDem15}.

\section{Kakeya phenomena}\label{sec_Kakeya}
The methods used to prove $\ell^2$ decoupling rely on ideas initially developed to study two critically important phenomena in harmonic analysis: restriction problems and Kakeya problems. Both of these have multilinear formulations, and a fundamental barrier in these fields is that the resolution of the multilinear conjectures does not fully resolve the original conjectures, which are now termed the ``linear'' cases. This stands in stark contrast to the setting of decoupling: a critical ingredient for success is that multilinear $\ell^2$ decoupling implies linear $\ell^2$ decoupling in full. Moreover, the Kakeya problem and its multilinear formulation, which we now introduce, play a role in the proof of $\ell^2$ decoupling.

\subsection{Kakeya and Besicovitch sets}
A pleasing question of Kakeya \cite{Kak17,FujKak17} asks for the smallest measure of a set in the plane, within which one can rotate a needle in any direction.  Besicovitch \cite{Bes28}  showed that such a set can have arbitrarily small measure. 
Independently of Kakeya, Besicovitch had himself asked for a set in $\R^2$ which contains a unit line segment in every direction (without the requirement that the needle can be rotated); he showed \cite{Bes19} that such a set can have zero measure, via a construction now called Perron trees \cite{Per28}. (There are other constructions as well; for example, a construction of Kahane \cite{Kah69} using the union of lines joining two copies of a Cantor set  can be shown using ideas of \cite{PSS03} to have zero measure.)

More generally,  for any $n \geq 2$ a \emph{Besicovitch set} is a compact set in $\R^n$ which contains a unit line segment in every direction; the Besicovitch construction shows that for every $n \geq 2$ there is a Besicovitch set with zero measure. Thus the sets can be ``small'' in some sense, but in fact they are simultaneously expected to be ``large'' in another sense, that of dimension:

\begin{conj}[Kakeya Conjecture for Hausdorff dimension of Besicovitch sets]\label{conj_Kakeya_H}
For every $n \geq 2$, all Besicovitch sets in $\R^n$ have Hausdorff dimension $n$.
\end{conj}
This conjecture is known to be true in dimension $n=2$, by work of Davies \cite{Dav71}.  (It may be interpreted as trivially true for $n=1$.)
It remains open in all dimensions $n \geq 3$.
Many people have contributed to progress on this conjecture: for example, Wolff \cite{Wol95} proved a landmark lower bound of $(n+2)/2$ for the Hausdorff dimension for $n \geq 3$; for large $n$  Bourgain \cite{Bou99} and then Katz and Tao \cite{KatTao02b} made further improvements. Currently in dimension $n=3$ a new record lower bound of $5/2+\ep_0$ for some small $\ep_0>0$ is due to Katz and Zahl \cite{KatZah17x}.  An analogous (weaker) conjecture is also posed for Minkowski dimension; the record lower bounds for Minkowski dimension are due to Katz, \L aba, Tao \cite{KLT00,LabTao01,KatTao02b}; see the surveys \cite{KatTao02,Lab08}.

As stated, this conjecture about Besicovitch sets might appear at first glance to be an isolated curiosity: find the haystack in the needles. But Besicovitch sets have many striking relationships to other areas. Two examples: first, recall that in studying  the Riemann zeta function one naturally encounters the task of bounding from above certain truncated Dirichlet series. Bourgain \cite{Bou91a} (see also \cite{Wol99}) has illustrated that bounds conjectured by Montgomery  \cite[p. 72-73]{Mon71} for these Dirichlet sums  (in the spirit of discrete restriction estimates) imply not only a result on the density of zeroes of the Riemann zeta function, but also imply Conjecture \ref{conj_Kakeya_H}. 
Second, recall Fefferman's Ball Multiplier Theorem \ref{thm_ball_multiplier} on truncations of Fourier transform integrals: Fefferman's proof of this (negative) result used a Besicovitch set construction, and is closely related to well-known conjectures on Bochner-Riesz means.
This is in keeping  with an important philosophy that estimates for certain oscillatory integrals imply Kakeya-type estimates (see \cite{Fef73b},\cite{Car15x} for an explication of such connections). 
Furthermore, although we cannot do this justice here, Kakeya problems have  significant connections to arithmetic combinatorics (see \cite{Lab08} for an introduction);  there are  important analogous Kakeya-type questions for sets containing circles, spheres, or lower-dimensional discs; and there are discrete analogues (e.g. the joints problem) and a finite field analogue (proposed in \cite{MocTao04} and resolved by Dvir \cite{Dvi09}). But our main interest is  the connection to restriction estimates and decoupling.

\subsection{The Kakeya conjecture for $\del$-tubes}
To see this connection, we must recast Conjecture \ref{conj_Kakeya_H} in terms of $\del$-tubes, that is, cylinders in $\R^n$ of proportions $\del \times \cdots \times \delta \times 1$, with volume $\approx \del^{n-1}$. 
In fact, thinking more quantitatively of Besicovitch sets in terms of collections of $\del$-tubes is very natural, as the computation of the Hausdorff dimension of a set introduces a $\del$-thickening of the set. Moreover, this quantification of the phenomenon  is useful in applications (see e.g. surveys \cite{KatTao02}, \cite{Tao01} and our later application in decoupling).\footnote{Although we do not use this formulation here, modern progress toward Conjecture \ref{conj_Kakeya_H} (e.g. as early as \cite{Bou91,Wol95}) focuses  on maximal functions taking maximal averages of a function over a collection of $\del$-tubes; proving that an appropriate maximal function is bounded from $L^p(\R^n)$ to an appropriate $L^q(\R^n)$ implies that the Hausdorff dimension in $\R^n$ is at least $p$, see e.g. \cite[Ch. 10]{Wol03}.}
Following the presentation of \cite{BCT06}, the relevant conjecture in terms of $\del$-tubes (which implies Conjecture \ref{conj_Kakeya_H}) is:

\begin{conj}[Kakeya conjecture]\label{conj_Kakeya_linear}
For each $\frac{n}{n-1}< q \leq \infty$ there exists a constant $C_q$, such that for every set $\Tcal$ of $\del$-tubes in $\R^n$ whose orientations form a $\del$-separated set of points on $S^{n-1} \subset \R^n$, 
\[  \| \sum_{T \in \Tcal} \onebf_T \|_{L^q(\R^d)} \leq C_q \del^{\frac{n-1}{q}} (\# \Tcal)^{1  - \frac{1}{q(n-1)}},\]
where $\#\Tcal$ denotes the cardinality of $\Tcal$.
\end{conj}

Roughly speaking, we can interpret this conjecture as follows: if many of the tubes $T \in \Tcal$ simultaneously overlapped over a large region, then  $\sum_{T \in \Tcal} \onebf_T$ would be large (up to a maximum $\# \Tcal$) over a large region. The largest the left-hand side could be is $\approx \del^{\frac{n-1}{q}} \# \Tcal$, which would occur if all the tubes overlapped on a set of maximal volume $\del^{n-1}$.  Since the conjectured inequality says the left-hand side cannot get this large, then  tubes that point in different directions must not overlap very much, so that we can intuitively expect the volume of their union (and hence its dimension) to be large.

The  inequality in Conjecture \ref{conj_Kakeya_linear} is trivially true for $q=\infty$ with $C_q=1$.
The analogous bound is known to be false at the endpoint $q = n/(n-1)$ by Besicovitch set constructions (see e.g. \cite{BCT06}),
but at the endpoint $q=n/(n-1)$, the conjecture can be modified to include $\del^{-\ep}$ on the right-hand side (or more precisely, a logarithmic factor).
The conjecture is known to be true for $n=2$ by C\'{o}rdoba \cite{Cor77}.

To see at least intuitively why there is a close connection between restriction phenomena and Kakeya phenomena, consider the truncated paraboloid $S = \{ (x,|x|^2) : x \in [-1,1]^{n-1}\} \subset \R^n$. It is possible to construct, for any small $0<\del<1$, a finitely overlapping cover of $S$ by ``slabs'' of proportions approximately $\del \times  \cdots \times \del \times \del^2$ with normals separated by approximately $\del$. Due to the uncertainty principle of Fourier transforms (which we'll return to in \S \ref{sec_uncertainty}), a smooth bump function supported in one of these slabs would have its Fourier transform essentially concentrated (and essentially constant) on a dual slab, with proportions $\del^{-1} \times \cdots \times \del^{-1} \times \del^{-2}$ (looking more like a tube than a slab); since these all point in different directions, it then seems natural that one encounters Kakeya phenomena.

In fact, it is known that the Restriction Conjecture \ref{conj_restriction_paraboloid} implies the Kakeya Conjecture \ref{conj_Kakeya_linear}. This is stated in Bourgain \cite[Eqn. (0.6)]{Bou91}, building on Fefferman's work on the ball multiplier theorem and \cite{BCSS89}; for a nice exposition see \cite[Prop. 10.5]{Wol03}.
In the other direction, it is thought that the Restriction Conjecture will not follow directly from the Kakeya Conjecture alone. 
Roughly speaking, the  Kakeya Conjecture is a statement about non-negative functions, so the source for the inequality does not come from cancellation between oscillating components, while the likely more subtle Restriction Conjecture is built from objects that are oscillatory by definition.\footnote{However, an argument of Carbery \cite{Car15x} shows that an appropriate square function estimate, if known, would imply the Kakeya conjecture (formulated in terms of maximal functions) as well as the restriction conjecture.}

\subsection{Multilinear phenomena: restriction and Kakeya}\label{sec_mult_Kakeya_intro}

We arrive at an important new setting:  bilinear, and more generally, multilinear analogues of the restriction and Kakeya problems. (Conjectures \ref{conj_restriction_paraboloid} and \ref{conj_Kakeya_linear} will now be considered the ``linear'' cases.)  Following \cite[\S6]{Tao04}, we will describe the spirit of the bilinear approach to proving the (linear) Restriction Conjecture \ref{conj_restriction_paraboloid}. Suppose we try to bound the $L^{p'}$ norm of $(gd\sig)\check{\;}$ for a function $g$ on $S$ when $p'=4$. By Plancherel,
\beq\label{restriction_bilinear_prep}
 \|(gd\sig)\check{\;}\|^4_{L^4(\R^n)} =  \|(gd\sig)\check{\;}(gd\sig)\check{\;}\|^2_{L^2(\R^n)}  = \|(gd\sig) *(gd\sig)\|^2_{L^2(\R^n)}.
 \eeq
Now as a convolution of  two measures supported on $S$, this becomes a problem of geometry rather than Fourier analysis, and is in some cases more amenable. (It is not an accident that this resembles the reasoning of Example 3 in \S \ref{sec_bi_orthogonality}.) More generally, one can replay this argument for $p'$ being any even integer, and reduce to $L^{p'/2}$ estimates for $p'/2$-fold convolutions of  measures. The real strength of this approach comes from considering, instead of just  one function $g$ and one compact region $S$, a bilinear restriction estimate of the form
\beq\label{restriction_bilinear}
\| (g_1d\sig_1)\check{\;}(g_2 d\sig_2)\check{\;} \|_{L^{p'/2}(\R^n)} \leq C_{p,q,S_1,S_2} \|g_1\|_{L^{q'}(S_1,d\sig_1)}\|g_2\|_{L^{q'}(S_2,d\sig_2)}
\eeq
for functions $g_i$ ($i=1,2$) supported on a smooth compact hypersurface $S_i$ with surface measure $d\sig_i$. To obtain a bound of the form (\ref{restriction_bilinear}), one assumes that  $S_1$ and $S_2$ are appropriately \emph{transverse}, in the sense that the set of unit normals to $S_1$ lies in some subset of $S^{n-1}$ that is sufficiently separated from the set of unit normals to $S_2$. Although this assumption of transversality may initially sound extravagant (certainly no transversality was apparent in our motivating example (\ref{restriction_bilinear_prep})), it is still suitable for the problem at hand. Indeed, after excising the (measure zero) diagonal $\Del =\{(\xi,\xi): \xi \in S\}$ from the product manifold $S \times S$, a special type of decomposition  partitions the non-diagonal set $(S \times S) \setminus \Del$ into smaller caps $S_1 \times S_2$ where $S_1,S_2$ are disjoint portions of $S$ that are separated by a distance proportional to their size (this is the characteristic property of a Whitney decomposition), and hence transverse as long as $S$ has curvature. This ultimately allows a restriction estimate for the manifolds $S_1 \times S_2$ to imply a restriction estimate for the manifold $S \times S$, and then finally a restriction estimate for the manifold $S$ (although some strength is lost in this passage).
Bilinear approaches have a long history, and were systematically employed for both restriction and Kakeya in work such as \cite{TVV98} and \cite{Wol01}; see further literature in the survey \cite{Tao04}.

Bennett, Carbery and Tao \cite{BCT06} took the multilinear approach to the limit by considering $n$-linear restriction estimates in $n$-dimensional settings, motivated by the insight that in this extreme case, all assumptions about curvature for the underlying manifold could be replaced by a transversality assumption.\footnote{See e.g. \cite{Bej16x} for the role of curvature in $k$-linear estimates in $\R^n$, where $k < n$.} Here we will record this assumption 
somewhat informally as follows: given for each $j=1,\ldots, n$ a parametrization $\Phi_j(x_j)$ of $S_j$ for $x_j$ in the parameter space $U_j$, then for each value of $x_1,\ldots, x_n$ in the parameter spaces $U_1,\ldots, U_n$, the unit normals $\om_1,\ldots, \om_n$ at the corresponding points on $S_1,\ldots, S_n$ have 
\[ |\det (\om_1,\ldots, \om_n)| \geq c_0>0,\]
for some uniform constant $c_0$ (and in particular they span $\R^n$); see \cite[Conj. 1.3]{BCT06} for the precise assumption. Bennett, Carbery and Tao proposed a conjecture in the following spirit:
\begin{conj}[Multilinear restriction conjecture, adjoint form]\label{conj_restriction_multilinear}
Suppose that $S_1,\ldots, S_n$ are $n$ compact $C^2$ hypersurfaces in $\R^n$ with surface measures $d\sig_1,\ldots, d\sig_n$, which are transverse in the above sense.
Then for $p' \geq \frac{2n}{n-1}$ and $q \leq p'(\frac{n-1}{n})$ there exists a constant $C$ (depending on $p,q,c_0$ and the $C^2$ norms of the parametrizations) such that 
\[ \| (f_1 d \sig_1)\check{\;} \cdots (f_n d\sig_n)\check{\;} \|_{L^{p'/n}(\R^n)} \leq C \| f_1 \|_{L^{q'}(S_1,d\sig_1)} \cdots \| f_n \|_{L^{q'}(S_n,d\sig_n)}, \]
for all $f_1 \in L^{q'}(S_1,d\sig_1), \ldots, f_n \in L^{q'}(S_n,d\sig_n)$.
\end{conj}
In the case where the functions parametrizing the hypersurfaces are linear, this conjecture was already known to be true, in the form of the Loomis-Whitney inequality \cite{LooWhi49}; the case $n=2$ was also known via approaches of C. Fefferman and Sj\"{o}lin.
The novel strategy of Bennett, Carbery and Tao was to access Conjecture \ref{conj_restriction_multilinear} via a multilinear version of the Kakeya conjecture. 

For $j=1,\ldots, n$, let $\Tcal_j$ be a set of $\del$-tubes in $\R^n$. We will say that the family $\Tcal_1,\ldots, \Tcal_n$ is transversal if for each $j=1,\ldots, n$, all the tubes in $\Tcal_j$ point in directions that are within a sufficiently small fixed neighborhood of the $j$-th standard coordinate vector in $S^{n-1}$.

\begin{conj}[Multilinear Kakeya conjecture]\label{conj_Kakeya_multilinear}
For each $n/(n-1) \leq q \leq \infty$ there exists a constant $C$ such that for all $\del>0$ and all transversal families $\Tcal_1,\ldots, \Tcal_n$ of $\del$-tubes in $\R^n$, 
\beq\label{Kakeya_multilinear}
 \| \prod_{j=1}^n ( \sum_{T_j \in \Tcal_j} \onebf_{T_j} )^{1/n} \|_{L^{q}(\R^n)} \leq C\del^{n/q} (\prod_{j=1}^n  \# \Tcal_j)^{1/n}.
 \eeq
\end{conj}

As stated, (\ref{Kakeya_multilinear}) for the value $q=n/(n-1)$ is called the endpoint case; if on the right-hand side of (\ref{Kakeya_multilinear}) we include a factor $\del^{-\ep}$ when $q=n/(n-1)$, this is considered near-optimal at the endpoint case. Note that (\ref{Kakeya_multilinear}) holds trivially for $q=\infty$; thus by interpolation it suffices to prove the stated inequality for the endpoint case $q=n/(n-1)$.  Note also that this multilinear conjecture  differs in flavor from the linear conjecture, both because the endpoint case is included, and because within a fixed set $\Tcal_j$, any number of tubes could be parallel or even identical.\footnote{A version of this last fact is later relevant for decoupling, as it allows analogues of (\ref{Kakeya_multilinear}) to be proved with  $\sum_{T_j \in \Tcal_j} \onebf_{T_j}$ replaced by any function $ \sum_{T_j \in \Tcal_j} c_{T_j}\onebf_{T_j}$ that is merely constant on the tubes $T_j$, with complex coefficients $c_{T_j}$.}

Bennett, Carbery and Tao proved that the  Multilinear Restriction Conjecture is equivalent (in an appropriate sense) to the Multilinear Kakeya Conjecture \cite[Prop. 2.1]{BCT06}, in contrast to relationship between the linear conjectures. Then, remarkably, they proved the Multilinear Kakeya Conjecture for $q>n/(n-1)$ and the near-optimal result (with $\del^{-\ep}$) at the endpoint $q=n/(n-1)$ (using a ``monotonicity'' argument based on properties of heat-flow and ideas of induction on scales). From this they deduced a near-optimal version of the Multilinear Restriction Conjecture \ref{conj_restriction_multilinear}.\footnote{Later, Guth \cite{Gut10}  proved the endpoint case of the Multilinear Kakeya Conjecture via the polynomial method of Dvir. Carbery and Valdimarsson \cite{CarVal13} re-proved this via the Borsuk-Ulam theorem. Finally, Guth \cite{Gut15} gave a 7-page proof of (\ref{Kakeya_multilinear}) (weakened by a factor $\del^{-\ep}$ on the right-hand side), showing the power of induction on scales, after which Bejenaru gave a second, short, proof of multilinear restriction \cite{Bej16xb}.}

\subsection{Multilinear approaches: why  decoupling stands apart}\label{sec_lin_mult_ideas}
The Bennett-Carbery-Tao results revealed the power of working in a multilinear setting, but the question remains what implications this will ultimately have for the linear conjectures. With the Multilinear Restriction Conjecture proved, Bourgain and Guth \cite{BouGut11} set a new record for the \emph{linear} Restriction Conjecture for hypersurfaces with positive definite second fundamental form (since further improved by Guth \cite{Gut16, Gut16xb}). But so far, no one has deduced the linear Restriction Conjecture in full for all dimensions.

In contrast,  \emph{multilinear $\ell^2$ decoupling implies linear $\ell^2$ decoupling in full}. For example, in the context of $\ell^2$ decoupling for hypersurfaces, Bourgain's initial work \cite{Bou13} built on the approach in \cite{BouGut11} of extracting linear restriction estimates from multilinear restrictions estimates, and in particular Bourgain proved the estimate of Theorem \ref{thm_decoupling_paraboloid} in the range $2 \leq p \leq 2n/(n-1)$, later pushed in \cite{BouDem15}  to reach the sharp range $p \leq 2(n+1)/(n-1)$. The multilinear-to-linear strategy also proved critical to the success of the Bourgain-Demeter-Guth work on the moment curve  (see \S \ref{sec_lin_mult_equiv}), and applies to $\ell^p$ decoupling for other $p$.

\subsubsection{Bilinear decoupling implies linear decoupling: a simple example}\label{sec_k2_bilinear_suffices}
The influence of multilinear  decoupling on linear decoupling is worth seeing in action. 
Consider for example our model decoupling inequality (\ref{decoupling_VMVT_model0}) for the moment curve, stated in terms of finite exponential sums.
In the exemplar case of degree $k=2$ and the curve $(t,t^2) \subset \R^2$, $\ell^2$ decoupling for $L^p$ is equivalent to the claim that for all $X \geq 1$,
\beq\label{k2_claim}
 \int_{(0,1]^2} |f_{(0,X]}(x)|^p dx \ll_\ep X^{p/2+\ep} 
 \eeq
for the critical exponent $p=p_2=6$, where for any interval $I  \subset \R$ we let 
\[ f_I(x)=f_I(x_1,x_2) = \sum_{t \in I \cap \Z} e(tx_1 + t^2 x_2).\]
We will show, somewhat informally \`a la Tao \cite{Tao15bblog}, that a bilinear decoupling inequality for these sums implies a corresponding linear decoupling inequality; we do so quite flexibly, with an additional parameter $\Del$:\footnote{It suffices to prove the Main Conjecture for $X$ being any power of 2, so that all constructions considered are dyadic, and iterated dissections of $(0,X]$ are well-defined.}
\begin{lemm}[Bilinear decoupling implies linear decoupling]\label{lemma_k2_bilinear}
Suppose that for some $2 < p \leq 6$ and some $\Del \geq 0$ we know a bilinear estimate of the following form: for every $\ep>0$, there exists $C_\ep = C(\ep,p) \geq 3$ such that 
\beq\label{k2_eta_bilinear}
 \int_{(0,1]^2} ( |f_{I_1} (x)|^{1/2} |f_{I_2}(x)|^{1/2})^p dx \leq C_\ep X^{p/2+ \Del + \ep }
 \eeq
holds for every pair of non-adjacent intervals $I_1,I_2$ in a dissection of $(0,X]$ into intervals of length $X/K$, for every constant $K$ with $C_\ep^{2/\ep} < K \leq 2 C_\ep^{2/\ep}$ and every $X$ a multiple of $K$. Then for the same $p$, 
\beq\label{k2_eta_linear}
 \int_{(0,1]^2} | f_{(0,X]}(x) |^p dx  \ll_{p,\ep} X^{p/2+\Del + \ep},
 \eeq
for all  $X \geq 1$ and every $\ep>0$.
\end{lemm}
The role the parameter $\Del$ would serve in a complete proof of (\ref{k2_claim}) is as follows: certainly (\ref{k2_eta_bilinear}) holds with $\Del =p/2$. If we know that for every $\Del>0$, (\ref{k2_eta_bilinear}) with $\Del$ implies (\ref{k2_eta_linear}) with $\Del$, and if we can show that any time (\ref{k2_eta_bilinear}) holds with $\Del$, it holds with a smaller $\Del' < \Del$, then we can ultimately show that (\ref{k2_eta_linear}) holds with $\Del=0$.

To prove Lemma \ref{lemma_k2_bilinear}, we'll use induction on scales, so we need to understand how an estimate like (\ref{k2_eta_linear}) rescales.
\begin{lemm}[Rescaled decoupling]\label{lemma_k2_rescales_well}
If for some $1 \leq p < \infty$ it is known that
\beq\label{k2_dilation_AX0}
 \int_{(0,1]^2} |f_{(0,X]}(x)|^p dx \ll_\ep X^{p/2+\Del+\ep} 
 \eeq
for all  $X \geq 1$ and every $\ep>0$, then for any $a \in \R$, $0<b<1$, we have that
\beq\label{k2_dilation_AX}
 \int_{(0,1]^2} |f_{(a,a+X^{1-b}]}(x)|^p dx \ll_\ep (X^{1-b})^{p/2+\Del+\ep} 
 \eeq
for all  $X \geq 1$ and every $\ep>0$.
\end{lemm}
We first note that the $L^p$ norm of $f_I$ is invariant under translations of $I$ by any fixed real number $a$, that is,  for  any fixed $X \geq 1$ and $1 \leq p < \infty$,
\beq\label{k2_translation_AX}
 \int_{(0,1]^2} |f_{(0,X]}(x)|^p dx= \int_{(0,1]^2} |f_{(a,a+X]}(x)|^p dx .
 \eeq
 When $p$ is an even integer, this is the statement of the translation-dilation invariance of the Vinogradov system of degree 2. 
To see this for general $p$, we observe that
\[  f_{(a,a+X]} (x_1,x_2) = \sum_{t \in (0,X]} e((t+a) x_1 + (t+a)^2 x_2)
	= e(ax_1 + a^2 x_2) f_{(0,X]} (x_1 + 2a x_2,x_2).\]			
Now upon integrating,
\[\int_{(0,1]^2} |f_{(a,a+X]} (x_1,x_2)|^p dx
	 = \int_{(0,1]^2} |f_{(0,X]} (x_1 + 2ax_2,x_2)|^p dx 
	  = \int_{(0,1]^2} |f_{(0,X]} (x_1 ,x_2)|^p dx, \]
	in which the change of variables $x_1 \mapsto x_1  - 2a x_2$ does not affect the range of integration, due to periodicity. 
Now that (\ref{k2_translation_AX}) is known, to prove (\ref{k2_dilation_AX}) it suffices to observe that for every $Y \geq 1$,
\beq\label{k2_dilation_AX'}
 \int_{(0,1]^2} |f_{(0,Y^{1-b}]}(x)|^p dx \ll_\ep (Y^{1-b})^{p/2+\Del+\ep} 
 \eeq
follows from (\ref{k2_dilation_AX0}) with $X=Y^{1-b}$.

Now to deduce the linear estimate (\ref{k2_eta_linear}) from the bilinear estimate (\ref{k2_eta_bilinear}), it suffices to show that
\beq\label{k2_eta'}
A(X) := X^{-p/2} \int_{(0,1]^2} | f_{(0,X]}(x) |^p dx  \ll_\ep X^{\Del + \ep}
 \eeq
as $X \maps \infty$, for all $\ep>0$. Let $X$ (sufficiently large) be given and let $K \leq X$ be a constant (in the assumed range), and dissect $(0,X]$ into $K$ disjoint intervals $I_1, \ldots, I_K$ of length $X/K$, so that $f_{(0,X]} = \sum_{1 \leq j \leq K} f_{I_j}$.
We apply the triangle inequality and then rearrange, to write
\begin{eqnarray*}
A(X)
	&\leq & X^{-p/2} \int_{(0,1]^2} \left( \sum_{1 \leq j \leq K} |f_{I_j}(x) | \right)^p dx \\
	& =& X^{-p/2} \int_{(0,1]^2} \left( (\sum_{1 \leq j \leq K} |f_{I_j}(x)|)^2 \right)^{p/2} dx\\
	&=& X^{-p/2} \int_{(0,1]^2} \left( \sum_{1 \leq i, j \leq K} |f_{I_i}(x)| |f_{I_j}(x)| \right)^{p/2} dx.
	\end{eqnarray*}
By H\"{o}lder's inequality applied with exponent $p/2 >1$ and its conjugate,
\beq\label{k2_AX_twocases}
 A(X) \leq X^{-p/2} K^{p/2-1}  \int_{(0,1]^2}\sum_{1 \leq i, j \leq K} |f_{I_i}(x)|^{p/2} |f_{I_j}(x)|^{p/2}  d x.
 \eeq

We will consider separately the near-diagonal contribution, from $|i-j| \leq 1$, and the off-diagonal contribution, from $|i-j|>1$.
First, for the strict-diagonal contribution, namely $i=j$, we observe that
\beq\label{k2_neardiag}
 X^{-p/2}K^{p/2-1} \int_{(0,1]^2} \sum_{1\leq i \leq K} |f_{I_i}(x)|^p dx 
	\leq X^{-p/2}K^{p/2-1} K \sup_{I \in \{I_1,\ldots, I_K\}} \int_{(0,1]^2}  |f_{I}(x)|^p dx .
	\eeq
Now we apply Lemma  \ref{lemma_k2_rescales_well} to conclude that this contribution is of size $A(X/K)$. 
Next, for the nearest-neighbor terms (with $|i-j|=1$)  we first apply Cauchy-Schwarz, so that for example in the case of $j=i+1$,
\[ \sum_{1 \leq i \leq K-1} |f_{I_i}|^{p/2} |f_{I_{i+1}}|^{p/2} \leq (\sum_{1 \leq i \leq K-1} |f_{I_i}|^p)^{1/2} 
	 (\sum_{1 \leq i \leq K-1} |f_{I_{i+1}}|^p)^{1/2}  \leq \sum_{1 \leq i \leq K} |f_{I_i}|^p;
\]
from this point, we can apply the analysis of (\ref{k2_neardiag}) and again obtain a contribution of size at most $A(X/K)$.

The off-diagonal contribution to (\ref{k2_AX_twocases}) is bounded by at most
\[
X^{-p/2} K^{p/2-1}  K^2 \sup_{\bstack{1 \leq i,j \leq K}{|i-j| >1}} \int_{(0,1]^2} |f_{I_i}(x)|^{p/2} |f_{I_j}(x)|^{p/2}  d x,
\]
and each of these integrals admits the assumed bilinear estimate  (\ref{k2_eta_bilinear}).
In total, we have shown that for every $\ep$, 
\beq\label{k2_bilinear_control}
 A(X) \leq C_\ep [ A(X/K) + X^{\Del+\ep} K^{p/2+1}],
 \eeq
for all constants $K$ in the assumed range, and $X$ sufficiently large.\footnote{Note the similarity of (\ref{k2_bilinear_control}) to Lemmas \ref{k2_EC_lemma2} and \ref{k3_EC_lemma2} in the efficient congruencing method.}

Now we will iterate this. Let $0<\ep<\Del$ be fixed, and choose $K$  such that
\beq\label{k2_K_large_enough}
C_\ep^{2/\ep} < K \leq 2 C_\ep^{2/\ep}.
\eeq
As a result, $K^{-\Del} \leq K^{-\ep} \leq 1/2$
and $\log C_\ep/ \log K < \ep/2$.
Working recursively we see that for any integer $t \geq 1$, 
\[ A(X) \leq C_\ep^t A(\frac{X}{K^t}) + K^{p/2+1} \sum_{s=0}^{t-1} C_\ep^s(\frac{X}{K^s})^{\Del+\ep}
	 \leq C_\ep^t A(\frac{X}{K^t}) +C_\ep^t X^{\Del+\ep} K^{p/2+1} \sum_{s=0}^{t-1} K^{-s(\Del+\ep)}.
	 \]
The last sum is  bounded by $2$ uniformly in $t$, since $K^{-\Del} \leq 1/2$.
 Thus for all $t \geq 1$,
\[ A(X) \leq 2C_\ep^t \{ A(\frac{X}{K^t}) + X^{\Del+\ep} K^{p/2+1}\} .\]
Now we choose $t$ just  large enough that $XK^{-t} \leq c$ (for any fixed constant $c  \geq 1$, e.g. $c=256$), e.g.  $\log (X/c)/\log K \leq t \leq 2 \log (X/c)/ \log K$;
 we then apply the trivial bound $A(Y)\leq Y^{p/2}$ with $Y=XK^{-t}$. Noting that $C_\ep^t \leq C_{\ep}^{2 \log X/\log K} =  X^{2\log C_\ep/ \log K}$,
 and recalling our constraint (\ref{k2_K_large_enough}), we conclude that 
\[ A(X) \leq 2c X^{2\log C_\ep/\log K}\{ 1 + X^{\Del+\ep } K^{p/2+1} \}
	\leq 2c X^{\ep}\{ 1 + X^{\Del+\ep } K^{p/2+1} \}
	\leq C_\ep' X^{\Del+2\ep},\]
which suffices to prove Lemma \ref{lemma_k2_bilinear}.
For a precise version of Lemma \ref{lemma_k2_bilinear} in the general setting, see Theorem \ref{thm_mult_lin}.

\subsection{An intuitive look at the role of multilinear Kakeya in decoupling}\label{sec_uncertainty}
We have seen a hint that a multilinear perspective is useful for proving $\ell^2$ decoupling results. 
Moreover, a form of multilinear Kakeya  is itself a key tool in proving the sharp multilinear decoupling result at the heart of the Bourgain-Demeter-Guth work. (In fact, their work relies on a hierarchy of multilinear Kakeya theorems, which consider not only tubes of proportions $\del \times \cdots \times \delta \times 1$, i.e. sets that are ``thin'' in $n-1$ dimensions, but also ``plates'' which are thin in $k$ dimensions, for each of $1 \leq k \leq n-1$.)

How does multilinear Kakeya enter the picture? We again consider the simple model case $(t,t^2)$ and see heuristically how bilinear Kakeya plays a role in proving a bilinear estimate such as  (\ref{k2_eta_bilinear}). (Formal reasoning is reserved to \S \ref{sec_ball_inflation_Kakeya}; roughly speaking, in the presentation of  \cite{Tao15bblog}, Kakeya allows the reduction of $\Del$ in (\ref{k2_eta_bilinear}) to some $\Del' < \Del$.) 

Consider a portion of our function of interest,
\beq\label{fJ_sum}
 f_{J}(x) = \sum_{t \in J} e(tx_1  + t^2x_2)=\sum_{n_1,n_2 \in \Z} c_{n_1,n_2} e^{2\pi i (n_1x_1 + n_2x_2)},
 \eeq
for an interval $J \subset (0,X]$, with Fourier coefficients $c_{n_1,n_2} = 1$ if $n_1 \in J$ and $n_2 = n_1^2$, and zero otherwise. If $J$ is an interval $J=(a,a+Y]$, then (\ref{fJ_sum}) gives the Fourier expansion of a function whose Fourier coefficients are supported (i.e. can be nonzero) within a rectangle of proportions $Y \times Y^2$.  Roughly speaking, the \emph{uncertainty principle} then indicates that this restriction on the Fourier side controls the behavior of the original function $f_{J}$ on spatial regions of complementary, or ``dual'' scales. 
This is related to the well-known Heisenberg uncertainty principle (and various other precise inequalities e.g. \cite{FolSit97}); we will develop only a looser, intuitive notion of this phenomenon here.

We first recall the principle of rescaling: if $f \in \Scal(\R^n)$ (so that $\hat{f} \in \Scal(\R^n)$ as well) and we set $f_\del(x) = f(x/\del)$ and $g^\del(x) = \del^n g(\del x)$ for $\del>0$, then 
\beq\label{rescaling_principle}
 (f_\del)\hat{\;}( \xi) = (\hat{f})^\del(\xi) .
\eeq
In particular if we think of $f$ as having compact support, say in the unit ball centered at the origin, then $f_{\del}(x)$ has support in a shrinking ball of radius $\del$ as $\del \maps 0$, while the region in which $(f_\del)\hat{\;}$ is non-negligible becomes increasingly spread out. This is one notion rooted in the uncertainty principle.

Suppose on the other hand that $f$ is an $L^2(\R^n)$ function such that $\hat{f}$ is supported in a ball of radius $R$ centered at the origin. Then we may see that there is a Schwartz function $\phi$ such that 
\beq\label{f_convolution}
f = f* \phi^R.
\eeq
 For indeed, if we select $\phi \in \Scal(\R^n)$ with the property that $\hat{\phi}$ is identically one on the ball of radius $1$ centered at the origin, then analogous to (\ref{rescaling_principle}) we see that $\hat{\phi}(\xi/R) =(\phi^R)\hat{\;}(\xi)$ is identically one on the ball of radius $R$ centered at the origin, so that $( f* \phi^R -f)\hat{\;}$ is identically the zero function, from which (\ref{f_convolution}) follows.
 
From (\ref{f_convolution}), it can be shown that as $x$ ranges over a ball $B^*$ of complementary radius $R^{-1}$, the values of $|f(x)|$ are essentially controlled by the average of $|f|$ over $B^*$, and thus in particular are effectively constant at that scale. In fact, the phenomenon gets even more interesting, and this is critical for its usage in the proof of $\ell^2$ decoupling: suppose instead that $\hat{f}$ is supported in an ellipsoid $E \subset \R^n$ with center $a$, and lengths $r_j$ along orthonormal axes $e_j$ for $1 \leq j \leq n$. Then as $x$ varies over any ellipsoid $E^*$ that is dual to $E$,\footnote{We say an ellipsoid $E^*$ is dual to $E$ if it has lengths $r_j^{-1}$ along the axes $e_j$, and any center.} $|f(x)|$ is essentially controlled by its average over $E^*$, in the sense that for any $x \in E^*$, 
 \beq\label{f_avg_ellipsoid}
  |f(x)| \leq C \frac{1}{|E^*|} \int_{\R^n} |f(x)| \phi_{E^*}(x) dx,
 \eeq
 where $\phi_{E^*}$ is a smooth function that is essentially one inside $E^*$ and decays outside $E^*$, analogous to the function $w_B$ tailored in (\ref{decoupling_weight_dfn}) to a ball. (The constant $C$ in (\ref{f_avg_ellipsoid}) depends on the rate of decay; see e.g. \cite[Ch. 5]{Wol03} for details.)
 An analogous phenomenon holds for a function $f$ with $\hat{f}$ supported in a rectangular region in $\R^n$, or in a parallelogram with a certain orientation in $\R^n$, and these considerations introduce Kakeya phenomena to the proof of decoupling.\footnote{As a technical note, the fact that there is a smooth weight function $\phi_{E^*}$ in (\ref{f_avg_ellipsoid}) and not a sharp cut-off function to $E^*$ indicates that the ``Schwartz tails'' introduced by $\phi_{E^*}$ cannot be avoided. Yet it is natural to use a simplified rubric as a heuristic device; to work rigorously one proceeds via a ``wave packet decomposition'' of $f$, which we will not describe here (see e.g. \cite[Dfn. 3.2]{BouDem15}).}

To see this, we go back and refine our understanding of $f_J$ in (\ref{fJ_sum}), with interval $J=(a,a+Y]$. Its Fourier coefficients are supported in a parallelogram $Q_J$ of slope $Y+2a$ that fits inside a rectangle (with sides parallel to the axes) of  
 proportions $Y \times Y^2$. In particular, the farther the interval $J$ is from the origin (so the larger $a$ is), the steeper the slope of the parallelogram $Q_J$. From our notion of the uncertainty principle, we see that $|f|$ is therefore essentially constant on dual parallelograms $P_J$, of slope $-(Y+2a)^{-1}$ within a rectangle of proportions $Y^{-1} \times Y^{-2}$. This  shows, roughly, that when trying to prove a bilinear decoupling inequality involving pieces $f_J$, it suffices to prove a bilinear decoupling inequality with each $f_J$ replaced by a certain local norm of $f_J$, call it $G_J$, that is essentially constant on  translates of a dual parallelogram $P_J$. (See \S \ref{sec_ball_inflation_Kakeya} for more formal reasoning.) 
 
Proving an appropriate bilinear decoupling inequality akin to (\ref{k2_eta_bilinear}) then ultimately hinges upon proving a relation of the form 
\beq\label{k2_ball_inflation}
\fint_{B} \prod_{j=1}^2 ( \sum_{J_j \subset I_j} G_{J_j}^{p/2}) \ll X^{\ep} \prod_{j=1}^2 ( \fint_B \sum_{J_j \subset I_j} G_{J_j}^{p/2}),
\eeq
where $B$ is an appropriate ball, $I_1$ and $I_2$ are any two  intervals in $(0,X]$ of length $X/K$ that are separated by at least  $X/K$, and for each $j=1,2$, $J_j$ runs over a dissection of $I_j$ into subintervals of smaller length $Y$, for an appropriate scale $Y$. (Here, $\fint_{B} = |B|^{-1} \int_B$.)
By our previous heuristic discussion, the  functions $G_{J_j}$ can be regarded as essentially constant on the dual parallelograms $P_{J_j}$, so for each $j=1,2$ we can think of 
\[ \sum_{J_j \subset I_j} G_{J_j}^{p/2} = \sum_{P \in \Pcal_j} c_P \onebf_P \]
for some coefficients $c_P$, where $\Pcal_j$ is the collection of all translates of the dual parallelograms $P_{J_j}$ that intersect the ball $B$, as $J_j$ varies in $I_j$. With this interpretation, (\ref{k2_ball_inflation}) becomes a type of bilinear Kakeya inequality, and can be verified as long as the two collections $\Pcal_1$ and $\Pcal_2$ are appropriately transversal. Roughly speaking, after book-keeping, (\ref{k2_ball_inflation}) will hold if for all choices $P_1 \in \Pcal_1, P_2 \in \Pcal_2$,
\[ |P_1 \intersect P_2 | \leq \frac{ |P_1||P_2|}{|B|}.\]
 Here, an appropriate choice of $K,Y,B$ provides that e.g. $|P_1|/|B|$ is just a bit smaller than 1, so that effectively one needs to show that for any choice of parallelograms $P_1 \in \Pcal_1$ and $P_2 \in \Pcal_2$, their intersection is appreciably smaller than the area of either of them.
 This transversality will ultimately come from the assumption that $I_1$ and $I_2$ are non-adjacent intervals in $(0,X]$ separated by at least $\gg X/K$, so that if $J_1$ is any subinterval in $I_1$ and $J_2$ is any subinterval in $I_2$, the slope of $P_{J_1} \in \Pcal_1$ differs from the slope of $P_{J_2} \in \Pcal_2$ by at least $\gg X^{-1}$, which turns out to be sufficient.
 This concludes our informal discussion of the role of restriction, Kakeya, and multilinear estimates in the proof of $\ell^2$ decoupling.

\section{Anatomy of the proof of $\ell^2$ decoupling for the moment curve}\label{sec_anatomy}
In this final section, we provide an overview of the Bourgain-Demeter-Guth proof of sharp $\ell^2$ decoupling for the moment curve. Our aim is not to provide a complete proof, but to outline (sometimes in broad terms) the key building blocks and their connections,  and also to illuminate some parallels to efficient congruencing.
While we will reference \cite{BDG16} exclusively, as mentioned before the methods of this paper are closely connected to the earlier Bourgain and Demeter canon.

We define for each $n \geq 2$, $2 \leq p \leq p_n=n(n+1)$ and $0< \del \leq 1$ the \emph{decoupling parameter} $V_{p,n}(\del)$ to be the smallest positive real number such that for each ball $B \subset \R^n$ with radius $\del^{-n}$, and for every   $g:[0,1] \maps \C$, 
 \beq\label{decoupling_statement}
 \| E_{[0,1]} g \|_{L^{p}(w_B)} \leq  V_{p,n}(\del) ( \sum_{\bstack{J\subset [0,1]}{|J| = \del}} \|E_J g\|^2_{L^{p}(w_B)} )^{1/2},
 \eeq
 where $E_{[0,1]} g$ is the extension operator (\ref{extn_Ga_op}) for the moment curve $\Ga$ in $\R^n$.\footnote{While Theorem \ref{thm_decoupling_moment_curve} is stated for all integrable functions $g$, we may at each step of the proof assume for example that $g$ is $C^\infty$; once the statement of Theorem \ref{thm_decoupling_moment_curve} holds for every $C^\infty$ function, then given $g \in L^1([0,1])$ and a sequence of functions $g_j \in C^\infty ([0,1])$ that converges to $g$ in $L^1$ norm, then for every interval $J \subset [0,1]$, $E_Jg_j$ will converge uniformly to $E_Jg$ for all $x \in \R^n$, so that the $L^p(w_B)$ norms do so as well.
 } 
 Theorem \ref{thm_decoupling_moment_curve}, the key $\ell^2$ decoupling result for the moment curve, is the statement that for $p=p_n=n(n+1)$ the critical exponent,
\beq\label{decoupling_parameter_claim}
V_{p_n,n}(\del) \ll_{n,\ep} \del^{-\ep}
\eeq
for every $\ep>0$.
 An advantage of defining the decoupling parameter is that it allows us to consider an inequality of the form (\ref{decoupling_statement}) for any value of $p$, and to compare various results that take us toward decoupling but are weaker than (\ref{decoupling_parameter_claim}). 
 It is helpful to compare (\ref{decoupling_parameter_claim}) to the much larger trivial bound 
\beq\label{decoupling_parameter_claim_trivial}
V_{p,n}(\del)\leq \del^{-1/2}.
\eeq
This follows (for any $1 \leq p  \leq \infty$) from Cauchy-Schwarz (in the spirit of (\ref{Hilbert_L2_ex3})), since 
\[ 
 \| E_{[0,1]} g \|_{L^{p}(w_B)} =  \| \sum_{\bstack{J\subset [0,1]}{|J| = \del}} E_{J} g \|_{L^{p}(w_B)} \leq ( \# J )^{1/2}( \sum_{\bstack{J\subset [0,1]}{|J| = \del}} \|E_J g\|^2_{L^{p}(w_B)} )^{1/2},
\]
and the number of  intervals $J$ in the dissection is $O(\del^{-1})$. (Again, we assume we are working with dyadic intervals, so that dissections into sub-intervals are well-defined.)

 The strategy of Bourgain, Demeter and Guth for proving (\ref{decoupling_parameter_claim}) (also present in earlier works on decoupling by Bourgain and Demeter) is to show that the exponent $\eta$ in a statement $V_{p,n}(\del) \ll \del^{-\eta}$ can be successively lowered to be arbitrarily small; this is achieved via an intricate iterative argument, which shows that $V_{p,n}(\del)$ is bounded above by decoupling parameters appearing in other types of decoupling estimates. 
 
 These other decoupling estimates involve: lower-dimensional settings (i.e. the analogue of (\ref{decoupling_statement}) for the moment curve $(t,t^2,\ldots,t^k)$ in $\R^k$ for $ k  \leq n-1$); estimates for multilinear settings (i.e. estimates like (\ref{decoupling_statement}) but involving products of extension operators associated to a family of separated intervals in $[0,1]$); estimates involving balls $B$  of different radii (recalling that when the ball has a larger radius relative to the scale of decoupling, the estimates are easier to prove); and estimates involving different choices for $p$ (and in particular, the powerful choice $p=2$). The point is that each of these other decoupling problems has its own decoupling parameter analogous to $V_{p,n}(\del)$, and by controlling $V_{p,n}(\del)$ iteratively in terms of these other decoupling parameters, one can eventually pass to an appropriate limiting setting in which all the decoupling parameters are known by some other means to be well-controlled. (In particular, the proof of Theorem \ref{thm_decoupling_moment_curve}  for dimension $n$ assumes that Theorem \ref{thm_decoupling_moment_curve} is already known for dimensions $2 \leq k \leq n-1$; recall that for $k=1$ the result is in some sense trivially true.)
 
  We will now describe  the key ingredients of the Bourgain-Demeter-Guth proof, divided into four stages: first, collecting some initial facts about linear decoupling in various settings; second, showing that linear decoupling may be controlled by multilinear decoupling, which is in turn controlled by averaged multilinear local norms; third, assembling three multilinear tools with a focus on multilinear Kakeya; and fourth, unwinding the iterative argument that ultimately controls the key multilinear estimate.

 \subsection{Key ingredients I:  initial facts about linear decoupling}
We assemble here three ingredients: first, the use of exponents smaller than the critical exponent, second, a rescaling principle for decoupling, and third, the special case of $\ell^2$ decoupling for $L^2$.
 
 \subsubsection{It suffices to prove decoupling for $2 \leq p<p_n$}\label{sec_general_smaller_p}
Recall that  $p_n = n(n+1)$ is the critical exponent. We will use the fact that our desired statement $V_{p_n,n}(\del) \ll_{n,\ep} \del^{-\ep}$ (for every $\ep>0$) follows from showing that  for every $p<p_n$ sufficiently close to $p_n$, $V_{p,n}(\del) \ll_{n,\ep} \del^{-\ep}$ (for every $\ep>0$). (This is effectively \cite[Lemma 9.2]{BDG16}, and may be stated more precisely in terms of $\eta_p$ as defined in (\ref{eta_p_dfn}).)\footnote{This type of argument is in some sense dual to that of \S \ref{sec_EC_k2}, whose analogue here would presumably work directly with the critical exponent, assume that there exists some $\eta_0>0$ for which $V_{p_n,n}(\del) \ll \del^{-\eta_0}$ cannot be improved, and obtain a contradiction.}

To prove this, the first step is to show (see \cite[Cor. 7.2]{BDG16}) that for any interval $J \subseteq [0,1]$ and any ball $B$ of radius $\gg 1$, 
\beq\label{same_ball_bigger_p}
 \| E_{J} g\|_{L^{p_n}(w_B)} \ll  \|E_{J}g\|_{L^p(w_B)},
\eeq
for every $1 \leq p < p_n$ (with an implicit constant independent of $B, g$). 
This is a consequence of certain monotonicity properties of weighted $L^p$ norms, combined with ideas of a  Bernstein inequality, which states that for a frequency-localized function, higher $L^p$ norms are controlled by lower $L^p$ norms.\footnote{In fact the standard Bernstein inequality is a direct consequence of the identity (\ref{f_convolution}) we have already seen: if $f \in L^1 + L^2$ and  $\hat{f}$ is supported in $B(0,R) \subset  \R^n$, then for $1 \leq p \leq q \leq \infty,$ $\|f\|_{L^q(\R^n)} \leq c R^{n(\frac{1}{p} - \frac{1}{q})} \|f\|_{L^p(\R^n)}$; see \cite[Prop. 5.3]{Wol03}.}

The next step is to assume that the decoupling inequality (\ref{decoupling_statement}) holds for $p$, and show that the right-hand side of (\ref{decoupling_statement}) for $p$ may be converted into the right-hand side of (\ref{decoupling_statement}) for $p_n$. To do so we apply H\"{o}lder's inequality with  $q=p_n/p >1$ and its conjugate $q'$ to see that 
\beq\label{desired_Bernstein}
  ( \sum_J \|E_J g\|_{L^p(w_B)}^2 )^{1/2} \ll \| \onebf \|_{L^{pq'}(w_B)} ( \sum_J \| E_J g\|_{L^{p_n}(w_B)}^2)^{1/2},
\eeq
where  $pq' = \frac{pp_n}{p_n-p}$  goes to infinity as $p \maps p_n$. Thus the contribution of $ \| \onebf \|_{L^{pq'}(w_B)}$ is a power of $B$, with the power going to zero as $p \maps p_n$, as desired.

\subsubsection{The linear decoupling inequality rescales well}\label{sec_decoupling_rescaling}

Any portion of the curve $\Gamma_{[0,1]}=\{ (t,t^2,\ldots, t^n) : t \in [0,1]\}$ may be rescaled to cover $\Gamma_{[0,1]}$; this is known as affine invariance (translation-dilation variance by another name). As a consequence, decoupling for $E_{[0,1]}g$ implies a form of decoupling for $E_{I}g$ for shorter intervals $I$. Precisely, we have:
\begin{prop}[Rescaling $\ell^2$ decoupling for $L^p$]\label{prop_decoupling_rescaling}
 For any $0<\del \leq 1$ and any $0<\rho \leq 1$, for every interval $I$ of length $\del^\rho$ and every ball $B \subset \R^n$ of radius $\del^{-n}$, 
\beq\label{decoupling_rescaling}
 \| E_{I} g \|_{L^{p}(w_B)} \leq  V_{p,n}(\del^{1-\rho}) ( \sum_{\bstack{J\subset I}{|J| = \del}} \|E_J g\|^2_{L^{p}(w_B)} )^{1/2}.
 \eeq
 \end{prop}
This is proved in a similar style to our analogous observation in Lemma \ref{lemma_k2_rescales_well}, although this generalization of the parabolic rescaling principle requires a more involved affine change of variables; see \cite[Lemma 7.5]{BDG16}. Let us pause to appreciate how useful this rescaling principle is: suppose that we first apply (\ref{decoupling_statement}) to decouple down to $\del_1$ and then apply (\ref{decoupling_rescaling}) with $|I| = \del_1$ to decouple down to $\del_2< \del_1$. Then we have 
 \begin{eqnarray*}
  \| E_{[0,1]} g \|_{L^{p}(w_B)}& \leq & V_{p,n}(\del_1) ( \sum_{\bstack{I\subset [0,1]}{|I| = \del_1}} \|E_I g\|^2_{L^{p}(w_B)} )^{1/2}\\
  	& \leq & 	V_{p,n}(\del_1) V_{p,n}(\del_2/\del_1)( \sum_{\bstack{I\subset [0,1]}{|I| = \del_1}} \sum_{\bstack{J\subset I}{|J| = \del_2}} \|E_J g\|^2_{L^{p}(w_B)})^{1/2}	,
	\end{eqnarray*}
which tells us that $V_{p,n}(\del_2) \leq V_{p,n}(\del_1)V_{p,n}(\del_2/\del_1)$; more generally the decoupling parameters satisfy a type of multiplicativity. This makes decoupling estimates prime candidates for methods that involve many scales simultaneously. 

 \subsubsection{Linear $\ell^2$ decoupling for $L^2$ is simple to prove}\label{sec_ell2_L2_linear}
It is natural that something special should occur when we consider $\ell^2$ decoupling for $L^2$: on the one hand, we recall e.g. from (\ref{Hilbert_L2_ex4}) that $\ell^2$ and $L^2$ interact nicely with each other, and in addition, $L^2$ is the nicest space in which to apply principles of orthogonality (or almost orthogonality) for a family of functions, via some knowledge of their Fourier supports. 
\begin{prop}[Linear $\ell^2$ decoupling for $L^2$]\label{prop_ell2_L2_linear}
Fix $n \geq 1$. For every $0<\del \leq 1$, for any interval $I$ (of length a multiple of $\del$), and for  any  ball $B$  of radius $\del^{-1}$ in $\R^n$, we have 
\beq\label{ell2_L2_ineq}  
\| E_I g \|_{L^2 (w_B)} \ll ( \sum_{\bstack{J \subset I}{|J| = \del}} \| E_J g\|_{L^2(w_B)}^2)^{1/2},
\eeq
for a dissection of $I$ into subintervals $J$ of length $\del$. 
\end{prop}
(Note that the statement of Theorem \ref{thm_decoupling_moment_curve} for dimension $n=1$ is an immediate consequence of Proposition \ref{prop_ell2_L2_linear}, since in this case the critical exponent is $p_1 = 2$.) 

This decoupling inequality has a very special flavor: this is claiming that decoupling down to the scale $\del$ can be detected over spatial balls of radius $\del^{-1}$; in contrast, our main  Theorem \ref{thm_decoupling_moment_curve} for $\ell^2$ decoupling for $L^p$ only claims that decoupling down to the scale $\del$ can be detected over spatial balls of the much bigger radius $\del^{-n}$.  Thus while the typical $\ell^2$ decoupling for $L^p$ result is restricted by the relationship 
$\mathrm{radius} = \mathrm{scale}^{-n}$, $\ell^2$ decoupling for $L^2$ allows the much better relationship $\mathrm{radius} = \mathrm{scale}^{-1}$. 

The proof of (\ref{ell2_L2_ineq}) rests on the following ideas. Invoking  a monotonicity property of weights \cite[Lemma 7.1]{BDG16} shows that (\ref{ell2_L2_ineq}) will follow if we can prove  that the sharply truncated $L^2(B)$ norm with weight $w=\onebf_B$ satisfies the property 
\beq\label{ell2_L2_ineq_sharp}
\| E_I g \|_{L^2(B)} \ll ( \sum_{\bstack{J \subset I}{|J| = \del}} \| E_J g\|_{L^2(\eta_{B})}^2)^{1/2},
\eeq
 for all balls $B$ with radius $\del^{-1}$, for a particular smooth weight $\eta_{B}$ we will choose advantageously, as in \cite[Lemma 8.1]{BDG16}. Fix $\eta$ to be a positive Schwartz function that is $\geq 1$ on the unit ball centered at the origin in $\R^n$, and such that the Fourier transform of $\sqrt{\eta}$ is supported, as a function of $\xi$, in a small neighborhood of the origin in $\R^n$, say $|\xi | \leq c_0 <1$. For any ball $B \subset \R^n$ with center $c$ and radius $R$, define 
\[ \eta_B(x) = \eta ( \frac{x-c}{R}).\]
Now the key point is that 
\[ \|E_I g\|^2_{L^2(B)} \ll \|E_I g\|^2_{L^2 (\eta_{B})}  = \| \sqrt{\eta_{B}} E_I g\|_{L^2(\R^n)}^2
	=\|  \sum_{\bstack{J \subset I }{|J| = \del}} \sqrt{\eta_{B}} E_J g  \|_{L^2(\R^n)}^2, \]
where in the first inequality we used that $\eta \gg 1$ on the unit ball.
If we can show that the  functions $\{\sqrt{\eta_{B}} E_J g)\}_J$ are orthogonal as long as $J,J'$ are distinct, non-adjacent intervals, then (\ref{ell2_L2_ineq_sharp}) follows immediately by an almost orthogonality argument in the style of Example 1 \S \ref{sec_almost_orthogonal_example}. By Plancherel's theorem, we can prove this orthogonality property by showing that  the collection of Fourier transforms $\{(\sqrt{\eta_{B}} E_J g) \hat{\;}\}_J$ have disjoint supports as long as $J,J'$ are distinct, non-adjacent intervals. 

We assume for simplicity that $B$ is centered at the origin.\footnote{The case of $B$ centered at another point just modulates $(\sqrt{\eta_{B}})\hat{\;}$ by an exponential of norm 1, since in general for any function $F(x)$, $F( \cdot - c)\hat{\;}(\xi) = e^{-2\pi i c \cdot \xi} \hat{F}(\xi).$}
Recall that for two functions $F, G$, 
\[ (FG)\hat{\;}  (\xi) =( \hat{F} * \hat{G} )(\xi) = \int_{\R^n} \hat{F}(\om)\hat{G}(\xi - \om) d\om.\]
Thus, if for example $\hat{F}$ is supported in a ball of radius $\mu$ then $\hat{F} * \hat{G}$ ``blurs'' the support of $\hat{G}$ by enlarging it by an  $O(\mu)$-neighborhood. In our case, we take
\[\hat{F}(\xi) = (\sqrt{\eta_{B}})\hat{\;} (\xi)=(\sqrt{\eta}(\frac{\cdot}{\del^{-1}}))\hat{\;}  (\xi) = \del^{-n} (\sqrt{\eta})\hat{\;}  (  \del^{-1} \xi),\]
 which by construction is supported in a small neighborhood of the origin, where $|\xi| \leq c_0 \del<\del$. 
On the other hand, we take $G=E_J g=(gd\sig)\check{\;}$ for $d \sig$ the surface measure of $\Gamma_J = \{(t,t^2,\ldots, t^n): t \in J\}$, so that $\hat{G} (\xi) = (E_J g)\hat{\;}(\xi)$ is supported on this arc $\Gamma_J$. 
Now by our previous observation, $\hat{F} * \hat{G} = (\sqrt{\eta_{B}} E_J g) \hat{\;}$ is supported on $\Gamma_J$ thickened by a $c_0\del$-neighborhood; since $J,J'$ are of length $\del$, the only way two of these thickened supports can intersect is if $J,J'$ are adjacent (or identical). This provides the almost orthogonality that we needed.

\subsection{Key ingredients II: how linear decoupling is controlled by multilinear objects}
 We now show how the ``linear'' decoupling statement (\ref{decoupling_statement}) follows from an analogous statement in a multilinear setting; this is a key means of gaining traction on the problem.

\subsubsection{The multilinear decoupling parameter}
We define a multilinear decoupling parameter $V_{p,n}(\del,K)$ as follows. Our multilinear objects will be $M$-multilinear for some large $M$ ($M=n!$ suffices), relative to intervals of length $1/K$, where $K$ is assumed to be sufficiently large, and certainly $K>M$.
 Given $n \geq 2$, $2 \leq p < p_n$, $0<\del \leq 1$ and a sufficiently large integer $M=M(n)$, define $V_{p,n}(\del,K)$ to be the smallest positive real number such that for every collection $I_1, \ldots, I_M$ of intervals of the form $[\frac{i}{K},\frac{i+1}{K}]$ that are pairwise non-adjacent, for each ball $B \subset \R^n$ with radius $\del^{-n}$, and for every $g:[0,1] \maps \C$  we have
\beq\label{decoupling_statement_mult}
 \|  ( \prod_{j=1}^M E_{I_j} g )^{1/M}\|_{L^{p}(w_B)} \leq  V_{p,n}(\del,K) (\prod_{j=1}^M( \sum_{\bstack{J\subset I_j}{|J| = \del}} \|E_J g\|^2_{L^{p}(w_B)} )^{1/2})^{1/M}.
 \eeq
 Roughly speaking,  we ultimately aim to show that for appropriate $M,K$, we have
$V_{p,n}(\del,K) \ll_{n,\ep} \del^{-\ep}$, for every $\ep>0$.
We will see that it is important that the multilinear inequality assumes that we only consider (distinct) non-adjacent intervals.

\subsubsection{Multilinear decoupling is equivalent to linear decoupling}\label{sec_lin_mult_equiv}
We first observe that linear decoupling implies multilinear decoupling, that is, $V_{p,n}(\del,K) \leq V_{p,n}(\del)$ for every integer $K \geq 1$. Formally, given an appropriate collection of pairwise non-adjacent intervals $I_1, \ldots, I_M$ of length $1/K$ and a function $g$ for which to prove (\ref{decoupling_statement_mult}),  by H\"{o}lder's inequality, 
\beq\label{mult_lin_step1}
  \|  ( \prod_{j=1}^M E_{I_j} g )^{1/M}\|_{L^{p}(w_B)}
	\leq  \prod_{j=1}^M ( \|  E_{I_j} g \|_{L^{p}(w_B)})^{1/M}.
	\eeq
The right-hand side stays unchanged if we replace $g$ by $\sum_{j=1}^M g_j$ where each $g_j=g$ on $I_j$ and vanishes outside $I_j$, so that for each $j$,  $E_{I_j}g = E_{[0,1]}g_j$.
We apply (\ref{decoupling_statement}) to each $j$-th factor on the right-hand side of (\ref{mult_lin_step1}); in total this shows that (\ref{mult_lin_step1}) may be bounded by 
\[
V_{p,n}(\del) \prod_{j=1}^M  ( \sum_{\bstack{J\subset I_j}{|J| = \del}} \|E_J g\|^2_{L^{p}(w_B)} )^{1/2M},
\]
which suffices for the multilinear claim.

But the real interest is in the other direction. In contrast to  restriction and Kakeya phenomena, multilinear decoupling implies linear decoupling  \cite[Thm. 7.6]{BDG16}:\footnote{We recall that this is true not just of decoupling for the moment curve, but for  other  decoupling settings, c.f.  \cite[Thm. 5.3]{BouDem15}.} 
\begin{theo}[Multilinear decoupling implies linear decoupling]\label{thm_mult_lin}
For every $2 \leq p \leq p_n$ and every integer $K \geq 1$, there exists a constant $C_{K,p}$ and $\ep_p(K)>0$, with $\lim_{K \maps \infty} \ep_p(K)=0$,  such that for every $0<\del \leq 1$,
\beq\label{mult_implies_lin}
V_{p,n}(\del) \leq C_{K,p} \del^{-\ep_p(K)} \sup_{\del \leq \del' <1}V_{p,n}(\del',K).
\eeq
\end{theo}
The argument to prove (\ref{mult_implies_lin}) is by induction on scales, in the style of  Bourgain-Guth  \cite{BouGut11} and  similar to the bilinear argument of \S \ref{sec_k2_bilinear_suffices}.  We sketch only the main points.
For a partition of $[0,1] = \union I_j$ into $K$ intervals of length $1/K$, we first gather the intervals into $K/M$ collections $\Ical$, each comprised of $M$ pairwise non-adjacent intervals $I_j$. Then trivially $|E_{[0,1]}g| \leq \sum_{\Ical} | \sum_{I_j \in \Ical }E_{I_j}g|.$ 
For each collection $\Ical$, the sum over $I_j \in \Ical$ is either dominated by one term, or the terms are all comparable and the sum is dominated by the geometric mean over all $M$ entries. The $L^p(w_B)$ norm of the first type of term is dominated by $\ll_{n,K} \sup_{|I| = K^{-1}} \|E_I g\|_{L^p(w_B)}$. After summing over the collections $\Ical$, the relevant supremum is trivially dominated by 
\beq\label{multi1}
\ll_{n,K} ( \sum_{\bstack{I \subset [0,1]}{|I| = K^{-1}}} \| E_I g\|_{L^p(w_B)}^2 )^{1/2},
 \eeq
so that on this portion we have effectively performed a linear decoupling down to scale $K^{-1}$ (without picking up a decoupling parameter). The second type of term we can bound in $L^p(w_B)$ norm via the definition of the multilinear decoupling parameter, that is, by
\[ \| \prod_{I_j \in \Ical} (E_{I_j} g)^{1/M} \|_{L^p(w_B)} \leq V_{p,n}(\del,K) ( \prod_{j=1}^M ( \sum_{\bstack{J \subset I_j}{|J| = \del}} \| E_J g\|_{L^p(w_B)}^2 )^{1/2})^{1/M}.  \]
Thus after summing over the $K/M$ collections $\Ical$, this contribution is in total
\beq\label{multi2}
\leq V_{p,n}(\del,K) \sum_{\Ical} (\prod_{j=1}^M ( \sum_{\bstack{J \subset I_j}{|J| = \del}} \| E_J g\|_{L^p(w_B)}^2 )^{1/2})^{1/M}
	\ll_{K,M} V_{p,n}(\del,K) ( \sum_{\bstack{J \subset [0,1]}{|J| = \del}} \| E_J g\|_{L^p(w_B)}^2)^{1/2} ,
	\eeq
	where we have employed the arithmetic-geometric mean inequality. 
At this stage, we have bounded $\|E_{[0,1]} g\|_{L^p(w_B)}$ by the sum of (\ref{multi1}) and (\ref{multi2}).

We now apply this same procedure not to $E_{[0,1]}$ but to $E_{I}$ for each $|I| = K^{-1}$ that appears in (\ref{multi1}); this effectively bounds each $\|E_{I} g\|_{L^p(w_B)}$ by a sum of terms like (\ref{multi1}) and (\ref{multi2}), but with (\ref{multi1}) replaced by an analogous term decoupling down to intervals of length $K^{-2}$, and with (\ref{multi2}) exhibiting a factor $V_{p,n}(\del K, K)$ instead of $V_{p,n}(\del,K)$ (by the rescaling principle, cf. Proposition \ref{prop_ell2_L2_linear}). Iterating this $\ell$ times, where $K^{-\ell} \approx \del$, we have brought the scale in (\ref{multi1}) down to $\del$, while picking up from (\ref{multi2}) a finite linear combination of terms $V_{p,n}(\del',K)$ for $\del \leq \del' \leq 1$. Ultimately (\ref{mult_implies_lin}) can then be deduced.

\subsubsection{Introducing the key multilinear players}\label{sec_mult_players}

It will be convenient to work now not with the $L^p(w_B)$ norm but with the normalized $L^p_\#(w_B)$ norm, which is defined as 
\[ \|F\|_{L^p_\# (w_B) }= ( \frac{1}{|B|} \int |F|^p w_B)^{1/p} \]
for $1 \leq p < \infty$.\footnote{As mentioned before, all steps in the rigorous induction are proved for $w_B$ with exponent of decay $E \geq 100n$ arbitrarily large. For simplicity, we omit this from our notation in this presentation.}
From now on, we also adopt the convention that for any $u >0$, $B^u$ denotes any ball of radius $\del^{-u}$ in $\R^n$.

Unwinding the definitions shows that for any $p$, the key multilinear decoupling inequality (\ref{decoupling_statement_mult}) is equivalent to the corresponding statement 
\beq\label{decoupling_statement_mult_sharp}
 \|  ( \prod_{j=1}^M E_{I_j} g )^{1/M}\|_{L^{p}_\#(w_{B^n})} \leq V_{p,n}(\del,K) (\prod_{j=1}^M( \sum_{\bstack{J\subset I_j}{|J| = \del}} \|E_J g\|^2_{L^{p}_\#(w_{B^n})} )^{1/2})^{1/M},
\eeq
for $I_1,\ldots, I_M$ any non-adjacent intervals of length $1/K$ (with $0 < \del \leq 1/K$).
We are now motivated to define  (for $1 \leq t < \infty$ and  $q,r>0$) the quantity
\[ D_t(q,B^r,g) =  ( \prod_{j=1}^M (\sum_{\bstack{J_j \subset I_j}{|J_j| = \del^q}} \|E_{J_j} g\|_{L^{t}_\#(w_{B^r})}^2 )^{1/2} )^{1/M};\]
this is the type of quantity we would  see on the right-hand side of a multilinear decoupling inequality (\ref{decoupling_statement_mult_sharp}) if we were able to prove an $\ell^2$ decoupling inequality for $L^t$, which detected decoupling at scales $\del^q$ over spatial balls $B^r$ of radius $\del^{-r}$. In particular, it is important to note that $D_p(1,B^n,g)$ appears on the right-hand side of (\ref{decoupling_statement_mult_sharp}), so that we may re-write for each $2 \leq p \leq p_n$ the definition of $V_{p,n}(\del,K)$ as the least positive real number such that 
\beq\label{decoupling_statement_mult_sharp_D}
 \|  ( \prod_{j=1}^M E_{I_j} g )^{1/M}\|_{L^{p}_\#(w_{B^n})} \leq V_{p,n}(\del,K)D_p(1,B^n,g),
\eeq
for all balls $B^n$ and functions $g$. 

A simple argument (employing e.g. the Cauchy-Schwarz and Minkowski inequalities \cite[\S9]{BDG16}) shows that for any small $0<u<n$ (to be specified later),
our object of interest on the left-hand side of (\ref{decoupling_statement_mult_sharp_D}) is bounded by 
\beq\label{mult_sumD}
  \|  ( \prod_{j=1}^M E_{I_j} g )^{1/M}\|_{L^{p}_\#(w_{B^n})} 
	\ll \del^{-u/2}  ( \frac{1}{| \Bcal_u(B^n)| } \sum_{B^u \in \Bcal_u(B^n)} D_p(u,B^u,g)^p)^{1/p},
	\eeq
where $\Bcal_u(B^n)$ is a finitely overlapping cover of $B^n$ by balls $B^u$ of smaller radius $\del^{-u}$. 
Then, since we only consider $p \geq 2,$ a Bernstein-type property as in (\ref{same_ball_bigger_p}) shows that the right-hand side of (\ref{mult_sumD}) is dominated by 
\beq\label{decoupling_D_average}
\del^{-u/2}  ( \frac{1}{| \Bcal_u(B^n)| } \sum_{B^u \in \Bcal_u(B^n)} D_2(u,B^u,g)^p)^{1/p},
\eeq
in which the $L^2$ norm plays  a key role.

This now motivates the definition of the  averaged multilinear local norm at the heart of the Bourgain-Demeter-Guth strategy:
for any ball $B^r$ of radius $\del^{-r}$, given a finitely overlapping cover $\Bcal_s(B^r)$ of $B^r$ by balls $B^s$ of smaller radius $\del^{-s}$, we define 
\[ A_p(q,B^r,s, g) = ( \frac{1}{| \Bcal_s(B^r)| } \sum_{B^s \in \Bcal_s(B^r)} D_2(q,B^s,g)^p)^{1/p}.\]
We can think of this as an $\ell^p$ average, over a cover of a ball of radius $\del^{-r}$ by balls of smaller radius $\del^{-s}$, of the quantity $D_2(q,B^s,g)$ which governs multilinear $\ell^2$ decoupling for $L^2(\R^n)$ down to scale $\del^q$, measured with respect to the relevant spatial ball of radius $\del^{-s}$. 
Recalling the strength of $\ell^2$ decoupling on $L^2$,  it is in particular reasonable to study $A_p(u,B^r,u,g)$, which picks out $\ell^2$ decoupling for $L^2(\R^n)$ down to scale $\del^u$, with spatial ball of the complementary radius $\del^{-u}$.

In particular,  using  (\ref{mult_sumD}) and (\ref{decoupling_D_average}) (and a careful argument passing from sharp cut-offs to smooth weights, \cite[\S 9]{BDG16}), it can be shown that 
\beq\label{mult_sharp_D}
\|( \prod_{j=1}^M E_{I_j} g )^{1/M}\|_{L^{p}_\#(w_{B^n})}   
	\ll \del^{-u/2} A_p(u,B^n,u,g).
\eeq
To prove multilinear $\ell^2$ decoupling, our central goal is now to control $A_p(u,B^n,u,g)$. 

\subsubsection{The principal result for $A_p(u,B^n,u,g)$}
In order to assemble the next important bound in precise terms, it helps to define, for each $2 \leq p \leq p_n$, the parameter $\eta_p \geq 0$ to be the unique real number (depending on $p,n$) such that 
\beq\label{eta_p_dfn}
 \lim_{\del \maps 0} V_{p,n}(\del) \del^{\eta_p + \sig} =0, \qquad  \text{and} \qquad \limsup_{\del \maps 0} V_{p,n}(\del) \del^{\eta_p - \sig} =\infty,
 \eeq
for every $\sig>0$. (Our goal is to show that $\eta_p=0$.)

The main technical iteration \cite[Thm. 8.3]{BDG16} of the Bourgain-Demeter-Guth proof, as repackaged in \cite[Thm. 9.1]{BDG16}, results in the following statement:
\begin{theo}\label{thm_BDG_Ap}
Fix $n \geq 3$ and let $2 \leq p < p_n$ be sufficiently close to $p_n = n(n+1)$. Suppose that Theorem \ref{thm_decoupling_moment_curve} holds for all dimensions $k \leq n -1$ (including $k=2$). Then for every positive number $W>0$ and for every sufficiently small $u>0$, we have for every $g: [0,1] \maps \C^n$, every $0<\del \leq 1$, and every ball $B^n \subset \R^n$ of radius $\del^{-n}$, 
\beq\label{Ap_bound}
A_p(u,B^n,u,g) \ll_{\sig,\ep,K,W}\del^{-\ep} \del^{-(\eta_p+\sig)(1 - uW)} D_p(1,B^n, g),
\eeq
for every $\ep,\sig>0$.
\end{theo}
If we compare (\ref{decoupling_statement_mult_sharp_D}) and (\ref{mult_sharp_D}), we can see right away that this will be useful in controlling the multilinear decoupling parameter $V_{p,n}(\del,K)$. 

\subsubsection{The endgame of $\ell^2$ decoupling for the moment curve}
With  the definition (\ref{decoupling_statement_mult_sharp_D}) of the multilinear decoupling parameter $V_{p,n}(\del,K)$ in hand, we apply Theorem \ref{thm_BDG_Ap} in (\ref{mult_sharp_D}) and take the supremum over all functions $g$, collections of pairwise non-adjacent intervals $I_j$ of length $1/K$ and balls $B^n$,  to see that 
\[ V_{p,n}(\del,K) \ll_{\sig,\ep,K,W} \del^{-u/2} \del^{-\ep} \del^{-(\eta_p + \sig)(1-uW)} ,\]
for every $\ep,\sig>0$.
Upon also recalling the key fact (\ref{mult_implies_lin}) that multilinear decoupling controls linear decoupling, we conclude that 
\[ 
\del^{-(\eta_p - \sig)}  \ll_{\sig,\ep,K,W} \del^{-\ep_p(K)} \cdot \del^{-u/2} \del^{-\ep} \del^{-(\eta_p + \sig)(1-uW)} ,
\]
for every $\ep,\sig$,
for a sequence of $\del \maps 0$. After a simple rearrangement, we deduce from the relationship of the exponents that 
\[ \eta_p \leq \frac{1}{2W} + \frac{\ep + \ep_p(K) + \sig (2 - uW)}{uW}.\]
The key property of this relation is that it holds for every $\ep,\sig>0$ arbitrarily small and every $K$ arbitrarily large (with $\ep_p(K) \maps 0$ as $K \maps \infty$), so we can conclude that
\beq\label{eta_compare_W}
\eta_p \leq \frac{1}{2W}.
\eeq
 Theorem \ref{thm_BDG_Ap} now allows us to take $W$ arbitrarily large, so that necessarily $\eta_p=0$. 
This proves $V_{p,n}(\del) \ll_{n,\ep} \del^{-\ep}$ for every $\ep>0$, and this provides the final step in the proof of Theorem \ref{thm_decoupling_moment_curve}.

 \subsection{Key ingredients III:  three tools to control $A_p(u,B^n,u,g)$}
The heart of the matter is thus to control the averaged multilinear local norms $A_p(u,B^n,u,g)$ as in Theorem \ref{thm_BDG_Ap}. This relies on an iteration using three main tools: (1) multilinear $\ell^2$ decoupling for $L^2$, (2) lower-dimensional decoupling, and (3) \emph{ball inflation} via a hierarchy of multilinear Kakeya estimates. We will now briefly describe these tools, and then turn to the philosophy of the iteration method.

\subsubsection{Multilinear $\ell^2$ decoupling for $L^2$}
We stated the linear case of $\ell^2$ decoupling for $L^2$ in Proposition \ref{prop_ell2_L2_linear}; the multilinear version is \cite[Lemma 8.1]{BDG16}:
\begin{prop}[Multilinear $\ell^2$ decoupling for $L^2$]\label{prop_ell2_L2_mult}
For every $0<\del \leq 1$, let $\Ical_1, \ldots , \Ical_M$ be collections of intervals of length a multiple of $\del$, such that each $\Ical_j$ is comprised of pairwise disjoint intervals. Then for any ball $B$ of radius $\del^{-1}$, we have 
\beq\label{ell2_L2_ineq_mult}  
( \prod_{j=1}^M ( \sum_{I \in \Ical_j} \| E_I g \|^2_{L^2_\#(w_B)} )^{1/2} )^{1/M} \ll ( \prod_{j=1}^M ( \sum_{\bstack{J \subset I, \; \mathrm{some}\; I \in \Ical_j}{|J| = \del}} \| E_J g\|_{L^2_\#(w_B)}^2)^{1/2})^{1/M}.
\eeq
\end{prop}
To prove this, we apply the linear case of $\ell^2$ decoupling for $L^2$ to each term $E_I g$.
(In fact, this method of proof works for any $M \geq 1$.)
As in Proposition \ref{prop_ell2_L2_linear}, this is a very strong kind of decoupling, which decouples down to scale $\del$ by detecting cancellation on spatial balls in $\R^n$ of radius $\del^{-1}$, for all dimensions $n \geq 1$. (Note that we made no assumption on the intervals being non-adjacent, as  transversality did not play a role.)

\subsubsection{Lower-dimensional decoupling}
The second tool used in the iteration scheme to prove Theorem \ref{thm_BDG_Ap} arises in an induction on dimension.
Roughly speaking, the idea is that if one views a particular portion of the curve $\Gamma = \{ (t,t^2,\ldots, t^n): t \in [0,1]\}$ at an appropriate scale, then it can be well-approximated by a curve in a lower-dimensional space. This can be advantageous for decoupling: we would only expect, if we look at spatial balls in $\R^n$ of radius $\del^{-n}$, to decouple down to scale $\del$ (that is, $\mathrm{radius}^{-1/n}$) on the curve in $\R^n$; but if we are instead able to translate this into decoupling on a curve in dimension $k<n$, we would expect to be able to decouple down to the smaller scale $\del^{n/k}$ (that is, $\mathrm{radius}^{-1/k}$).  For example, if in the case $n=3$ we look at the arc $\{ (t,t^2,t^3) : t \in [0,\del]\}$ for some small $\del<1$ (so that $\del^3 < \del$), then this is contained in an $O(\del)$-thickening of the arc $\{ (t,t^2,0) : t \in [0,\del]\}$ in the 2-dimensional plane, and this enables a decoupling result on balls of radius $\del^{-3}$ to decouple all the way down to scale $\del^{3/2}$, not just to scale $\del$.
Precisely, the result is \cite[Lemma 8.2]{BDG16}:
\begin{prop}[Lower-dimensional decoupling]\label{prop_lower_dim_decoupling}
Let $n \geq 3$ be fixed. For every $0<\del \leq 1$ and every $3 \leq k \leq n$, for any interval $I \subset [0,1]$ of length a multiple of $\del^{\frac{n}{k-1}}$, for every ball $B$ in $\R^n$ of radius $\del^{-n}$ and every $2\leq p \leq p_n$,
\[ \|E_I g\|_{L^p_\sharp(w_B)} \ll V_{p,k-1}(\del^{\frac{n}{k-1}}) ( \sum_{\bstack{J \subset I}{|J| = \del^{\frac{n}{k-1}}}} \| E_J g\|_{L^p_\sharp(w_B)}^2)^{1/2}.\]
Here for each $k$, $V_{p,k}(\cdot)$ denotes the decoupling parameter for $\ell^2$ decoupling in $L^p(\R^k)$, analogous to the definition (\ref{decoupling_statement}).
\end{prop}
Here we see vividly why in proving Theorem \ref{thm_decoupling_moment_curve} for dimensions $n \geq 3$, it is useful to assume that it is already known for  dimensions $k \leq n-1$. We observed in \S \ref{sec_ell2_L2_linear} why the theorem is true for $n=1$, and the truth of the theorem for $n=2$ is known by \cite{BouDem15}; hence we focus now on $n \geq 3$.

\subsubsection{The key ball inflation statement}
Finally, the most important tool is the mechanism of ball inflation, which is enabled by a sophisticated application of multilinear Kakeya-type results. Recall that if we are considering larger spatial balls, we expect to be able to decouple down to a smaller scale.  The mechanism of ball inflation, roughly speaking, allows us to control a statement for smaller spatial balls by a statement for larger spatial balls while preserving the scale of decoupling, thus opening up the door to decoupling down to a smaller scale. 
\begin{prop}[Ball inflation]\label{prop_ball_inflation}
Fix $ 1 \leq k \leq n-1$, $2n \leq p \leq p_n$, and take $M=n!$. For any ball $B^{k+1}$ in $\R^n$ with radius $\del^{-(k+1)}$ and a cover $\Bcal_k(B^{k+1})$ of $B^{k+1}$ by balls $B^k$ of radius $\del^{-k}$, then for every $g:[0,1] \maps \C$, 
\begin{multline}\label{ball_inflation_ineq}
\frac{1}{|\Bcal_k(B^{k+1})|} \sum_{B^k \in \Bcal_k(B^{k+1})} ( \prod_{j=1}^M (\sum_{\bstack{J_j \subset I_j}{|J_j| = \del}} \|E_{J_j} g\|_{L^{\frac{pk}{n}}_\#(w_{B^k})}^2 )^{1/2} )^{p/M}
\\
\ll_{\ep,K} \del^{-\ep} ( \prod_{j=1}^M (\sum_{\bstack{J_j \subset I_j}{|J_j| = \del}} \|E_{J_j} g\|_{L^{\frac{pk}{n}}_\#(w_{B^{k+1}})}^2 )^{1/2} )^{p/M}
\end{multline}
with an implicit constant independent of $g,\del$ and the ball $B^{k+1}$.
\end{prop}
A similar result holds for $k=n$ \cite[Eqn. 17]{BDG16}.
Note that (\ref{ball_inflation_ineq}) is not itself a decoupling inequality, since the intervals $J_j$ remain the same on both sides:  instead the point is that if initially we were using spatial balls of radius $\del^{-k}$, we have now transformed the multilinear quantity into being measured with respect to larger spatial balls of radius $\del^{-(k+1)}$, so that when decoupling is performed, we are allowed to go down to a smaller scale. 
(In terms of the notation of the previous section, this is allowing us, for each $1 \leq k <n$, to bound a certain average $A_{pk/n}(1, B^{k+1},k,g)$ by $\del^{-\ep} D_{pk/n}(1,B^{k+1},g)$.) 
We note also that a version of one layer of the ball inflation hierarchy appeared in earlier works of Bourgain and Demeter; see for example the case $k=n-1$ in \cite[Thm. 9.2]{BouDem17a}.

\subsubsection{How multilinear Kakeya for plates enables ball inflation}\label{sec_ball_inflation_Kakeya}
We consider how to prove the ball inflation Proposition \ref{prop_ball_inflation} for a fixed $1 \leq k <n$. 
By applying dyadic pigeonholing, we may assume we have reduced our attention to proving (\ref{ball_inflation_ineq}) where for each $j=1,\ldots,M$ we have restricted the smaller intervals $J_j \subset I_j$ to a set (of cardinality $N_j$, say)  for which 
\beq\label{E_J_const}
C \leq \|E_{J_j} g\|_{L^{\frac{pk}{n}}_\#(w_{B^{k+1}})} \leq 2C
\eeq
for some constant $C$. 
(Since this norm is trivially bounded above by $O(\del)$, we need only consider $O(\log (\del^{-1}))$ values for $C$, so this introduces allowable logarithmic losses $\ll_\ep \del^{-\ep}$.)
From now on we will only consider such $J_j$, and aim to prove (\ref{ball_inflation_ineq}) for such terms.

By applying H\"{o}lder's inequality to the sum over $J_j$ on left-hand side of (\ref{ball_inflation_ineq}) (with exponent $q=pk/(2n)\geq 1$ and its conjugate, recalling $p \geq 2n$) we may bound it by
\[ ( \prod_{j=1}^M N_j^{\frac{1}{2} - \frac{n}{pk} })^{p/M} \frac{1}{|\Bcal_k(B^{k+1})|} \sum_{B^k \in \Bcal_k(B^{k+1})} ( \prod_{j=1}^M \sum_{J_j \subset I_j} \|E_{J_j} g\|_{L^{\frac{pk}{n}}_\#(w_{B^k})}^{pk/n})^{n/kM},\]
which unifies the $pk/n$-moment with the $L^{pk/n}$ norm.
We observe that the ball inflation inequality (\ref{ball_inflation_ineq}) for this restricted collection of $J_j$ will follow if we can show that 
\beq\label{ball_inflation_step1}
  \frac{1}{|\Bcal_k(B^{k+1})|} \sum_{B^k \in \Bcal_k(B^{k+1})} ( \prod_{j=1}^M \sum_{J_j \subset I_j} \|E_{J_j} g\|_{L^{\frac{pk}{n}}_\#(w_{B^k})}^{pk/n} )^{n/kM}
\ll \del^{-\ep}  ( \prod_{j=1}^M \sum_{J_j \subset I_j} \|E_{J_j} g\|_{L^{\frac{pk}{n}}_\#(w_{B^{k+1}})}^{pk/n}  )^{n/kM}. 
\eeq
This is because of the dyadic pigeonholing assumption, since (\ref{E_J_const}) allows us to verify that 
\[ ( \prod_{j=1}^M N_j^{\frac{1}{2} - \frac{n}{pk} })^{p/M}  ( \prod_{j=1}^M \sum_{J_j \subset I_j} \|E_{J_j} g\|_{L^{\frac{pk}{n}}_\#(w_{B^{k+1}})}^{pk/n} )^{n/kM}
	\approx  ( \prod_{j=1}^M (\sum_{J_j \subset I_j} \|E_{J_j} g\|_{L^{\frac{pk}{n}}_\#(w_{B^{k+1}})}^2 )^{1/2} )^{p/M},\]
	which yields the desired right-hand side in (\ref{ball_inflation_ineq}).

The key passage taking place in (\ref{ball_inflation_step1}) is from a cover of $B^{k+1}$ by smaller balls $B^{k}$ on the left, to the larger ball $B^{k+1}$ on the right. We will describe this passage in broad terms. The idea is that for each fixed interval $J$ encountered on the left-hand side, we will cover $\union B^k$ by a collection of ``plates'' that are tailored to $J$ in the following way. We dissect $\union B^k$ into a collection of disjoint plates such that in $k$ dimensions they have length $\del^{-k}$ and in the remaining $n-k$ dimensions they have longer length $\del^{-(k+1)}$; we can do this so that all these plates fit within a dilation of $B^{k+1}$ by at most 4, and such that any smaller ball $B^k$ in the union fits within a dilation of one of the plates by a factor of 2. 

Finally, we specify the orientation of the plates, determined by the interval $J$ under consideration. We specify the parametrization 
\[ \Gamma = \{\Phi(t) = (t,t^2,\ldots, t^n): t \in [0,1]\}\]
of the moment curve. Supposing that $J = [t_J - \del/2,t_J +\del/2]$, we specify that the first $k$ axes of the plates tailored to $J$ span the $k$-dimensional linear subspace 
\beq\label{plate_orientation}
 \langle \Phi'(t_J),\Phi''(t_J),\ldots, \Phi^{(k)}(t_J)\rangle 
 \eeq
and the remaining $n-k$ axes span the orthogonal complement of this subspace in $\R^n$. 
Since the plates are disjoint, for any $x \in B^{k+1}$ we may specify the unique plate $P(x)$ that contains $x$.

Now our goal is to modify the quantities we consider in (\ref{ball_inflation_step1}) so that we can consider functions that are constant on each of the plates $P$ associated to $J$ (thus taking a step toward a Kakeya setting). 
It suffices to define a function $G_J$ as follows: for any $x$ in the union of the plates associated to $J$ and covering $\union B^k$, we set 
\[ G_J (x) = \sup_{y \in 2P(x)} \| E_J g\|_{L^{\frac{pk}{n}}_\# (w_{B(y,\del^{-k})})}.\]
This jitters the center of the ball over which the norm is taken (while preserving its radius), and then takes a sup over the possible centers. Note that this function is constant as $x$ varies over any single fixed plate $P$ associated to $J$. Moreover, if $x$ lies in a particular ball $B^k$ of radius $\del^{-k}$, recalling that this ball lies in $2P(x)$ we see that 
\[ \| E_J g\|_{L^{\frac{pk}{n}}_\# (w_{B^k})} \leq G_J(x).\]
A far more subtle argument (via Fourier multipliers and monotonicity of weights, see \cite[Thm. 6.6]{BDG16}) shows that in the other direction,
\[  \| G_J \|_{L^{\frac{pk}{n}}_\# (w_{4B^{k+1}})} \leq  \| E_J g\|_{L^{\frac{pk}{n}}_\# (w_{B^{k+1}})}.\]
As a result of these two comparisons, we may ultimately conclude that the ball inflation inequality (\ref{ball_inflation_step1}) will hold if we can verify that 
\beq\label{ball_inflation_step2}
 \fint_{4B^{k+1}} \prod_{j=1}^M (\sum_{J_j \subset I_j} G_{J_j}^r)^s
 	\ll_{\ep,K}\del^{-\ep} \prod_{j=1}^M (  \fint_{4B^{k+1}} \sum_{J_j \subset I_j} G_{J_j}^r)^s,
\eeq
where $r=pk/n$ and $s = n/kM$, and we let $\fint_{\Omega} = |\Omega|^{-1} \int_\Omega$ for any compact region $\Omega$.

We now return to the geometry at hand. Fix an index $1 \leq j \leq M$. For this fixed $j$, for each $J_j \subset I_j$ we have a collection of  plates  within $4B^{k+1}$ with proportions $(\del^{-k})^k \times (\del^{-(k+1)})^{n-k}$ and with orientation along the moment curve specified by (\ref{plate_orientation}) according to the location of the interval $J_j$ within $I_j$. We can assemble all of the plates corresponding to all the $J_j$ we consider within the fixed interval $I_j$, and call this collection $\Pcal_j$; this now possibly includes many parallel plates and many overlapping plates. But importantly, if we define for each fixed $1 \leq j \leq M$ the function 
\[ F_j = \sum_{J_j \subset I_j} G_{J_j}^r, \]
since $G_{J_j}$ is constant on each plate associated with $J_j$, we see that $F_j$ behaves like a linear combination 
\beq\label{FPMj}
 F_j =  \sum_{P \in \Pcal_j} c_P \onebf_P.
 \eeq
Now (\ref{ball_inflation_step2}) is equivalent to the claim (recalling $s=n/kM$) that 
\[  \fint_{4B^{k+1}} (\prod_{j=1}^M F_j )^s
 	\ll_{\ep,K}\del^{-\ep} \prod_{j=1}^M (  \fint_{4B^{k+1}}F_j)^s.\]
This is in fact the statement of a strong, uniform multilinear Kakeya inequality for plates that are appropriately thin in $k$ directions and thick in $n-k$ directions, and it will hold as long as the families $\Pcal_1,\ldots, \Pcal_M$ appearing in (\ref{FPMj}) are sufficiently transverse to each other, as families. (Note that we are not assuming transversality within any particular family $\Pcal_j$, in keeping with the setting described for multilinear Kakeya for tubes in \S \ref{sec_mult_Kakeya_intro}.)
This transversality of the families  is assured by the fact that for each $j$, the orientations of the plates in $\Pcal_j$ are determined by points within the interval $I_j$, and in the multilinear setting we have guaranteed that $I_1,\ldots, I_M$ are non-adjacent intervals in $[0,1]$ of length $1/K$. Thus for any $j \neq j'$, the points controlling the orientation of plates in $\Pcal_j$ are separated from those controlling the orientation of plates in $\Pcal_{j'}$ by at least $1/K$.
With this key insight in hand, the veracity of the relevant multilinear Kakeya result  is proved in \cite[\S 6]{BDG16}, building on the methods of \cite{Gut15}. (The utility of the choice $M=n!$ may be seen in the proof of \cite[Lemma 6.3]{BDG16}, where it allows the verification of the nondegeneracy condition \cite[Eqn. (8)]{BDG16}, required for the application of a multilinear H\"{o}lder-Brascamp-Lieb inequality of \cite{BCCT10}.)

\subsection{Key ingredients IV: the iteration}
\subsubsection{Iteration: a philosophy}
We have now seen three tools that propel us toward being able to decouple down to a finer scale: $\ell^2$ decoupling for $L^2$, lower-dimensional decoupling, and ball inflation.
How do we assemble these to prove Theorem \ref{thm_BDG_Ap}?  
We will not be able to do justice to the full structure of the iterative argument, but we can simulate its philosophy. (Tao \cite{Tao15ablog} points out the iteration can alternatively be seen from the perspective of Bellman functions, surveyed in \cite{NazTre96}.)

To gain some intuition, let us recall what was so important about the output of Theorem \ref{thm_BDG_Ap}: taking for the moment a very simplistic view, it allowed us  to see that for arbitrarily large $W$ (and comparably  small $u$) we could write
\[  \del^{-\eta_p} \approx V_{p,n}(\del)\approx V_{p,n}(\del,K) \ll \del^{-\eta_p(1-uW) -u/2 },\]
 so that we could cancel the factors $\del^{-\eta_p}$ and conclude that $\eta_p \leq 1/2W$. 

How could we achieve such a majorization of $V_{p,n}(\del,K)$? First, suppose we knew that for any positive number $\be$ (and $u$ sufficiently small that $u\be<1$), we could prove the majorization\footnote{In this discussion, to simplify notation, we have stopped notating the dependence of various quantities on the function $g$.}
\beq\label{Ap_Dp_relation}
 A_p(\be u, B^n, \be u) \ll V_{p,n}(\del^{1-u \be}) D_p(1,B^n).
 \eeq
 This is desirable, as by the definition of $\eta_p$ in (\ref{eta_p_dfn}) we can see $V_{p,n}(\del^{1-u \be}) \approx \del^{-\eta_p(1-u\be)}$, and furthermore we know that $D_p(1,B^n)$ is intimately related to the multilinear decoupling parameter $V_{p,n}(\del,K)$  via (\ref{decoupling_statement_mult_sharp_D}).
 
Second, suppose hypothetically that we could show that the quantity $A_p(u,B^n,u)$ we consider in Theorem \ref{thm_BDG_Ap} (always for $u$ sufficiently small) could be majorized, for any large $L$ of our choice, by a product of $L$ such quantities,
\beq\label{Ap_fake_product}
A_p(u,B^n,u) \ll \prod_{\ell=1}^L A_p(\be_\ell u , B^n, \be_\ell u),
\eeq
with each $\be_\ell \geq 1$, so that as $L$ increases, $\sum_{\ell \leq L } \be_\ell$ may be made arbitrarily large. (This would be a significant achievement, since this would be saying that $A_p(u,B^n,u)$, which decouples down to scale $\del^{u}$, could be controlled by quantities that decouple down to the  smaller scale $\del^{\be_\ell u}$.) With such a result in hand, we would apply (\ref{Ap_Dp_relation}) to each of the $L$ factors in (\ref{Ap_fake_product}), and see  by the definition (\ref{eta_p_dfn}) of $\eta_p$ that
\beq\label{Ap_fake_product'}
 A_p(u,B^n,u)\ll \del^{-\eta_p (L - u  \sum_{\ell \leq L} \be_\ell)} D_p(1,B^n)^{L}.
 \eeq
Now if we compare this to the claim (\ref{Ap_bound}) of Theorem \ref{thm_BDG_Ap}, it looks somewhat hopeful: the exponents $L$ should instead be $1$ as in (\ref{Ap_bound}) (and this is important), but we do at least see that by taking $L$ large we could make the sum $\sum_{\ell \leq L} \be_\ell$ in the exponent arbitrarily large, and in particular as large as $W$, for any $W$.

In fact, the relation (\ref{Ap_Dp_relation}) is true and relatively simple to prove. (H\"{o}lder's inequality converts the left-hand side into $D_p(\be u, B^n)$ and then rescaling via Proposition \ref{prop_decoupling_rescaling} leads to the right-hand side, see \cite[Eqns. 74--75]{BDG16}.) 
The hypothetical relation (\ref{Ap_fake_product}) is not correct, but it provides some inspiration for a  result that is true. 

The more complicated version of (\ref{Ap_fake_product}), which is true,  takes the shape
\beq\label{Ap_prep}
 A_p(u,B^n,u) \ll \del^{-\ep} V_{p,n}(\del)^{1 - \sum_{0}^r \ga_j} \times \prod_{j=0}^r A_p(b_ju, B^n, b_ju)^{\ga_j} \times D_p(1,B^n)^{1 - \sum_0^r \ga_j}
\eeq
where we may take $r$ as large as we like, but $p<p_n$ must be sufficiently close to the critical exponent $p_n$ and $u$ must be sufficiently small. Here $\ga_j$ and $b_j$ are certain positive values (with $b_j>1$) such that $\sum_{j=0}^\infty \ga_j b_j$ may be computed precisely, but here we note only that it is a positive constant that is strictly bigger than 1 (although not necessarily large); this is the content of \cite[Thm. 10.1]{BDG16}.

We could at this point apply the (correct) relation (\ref{Ap_Dp_relation}) to each of the $A_p$ factors on the right-hand side of (\ref{Ap_prep}); this would show (after cancelling certain factors that align nicely) that (\ref{Ap_prep}) becomes
\beq\label{Ap_weak}
 A_p(u,B^n,u) \ll  \del^{-\ep} \del^{\eta_p( 1 - u \sum_0^r b_j \ga_j)} D_p(1,B^n).
\eeq
Comparing this to the claim (\ref{Ap_bound}) of Theorem \ref{thm_BDG_Ap}, we see we have made some progress over our hypothetical (incorrect) relation (\ref{Ap_fake_product'}): instead of the troublesome large powers (the exponents $L$ in (\ref{Ap_fake_product'})) we have precisely one factor of $D_p(1,B^n)$, as we hoped. But (\ref{Ap_weak}) is not yet what we need: recall that if we take $r$ sufficiently large then $\sum_0^r b_j \ga_j >1$, but it is not necessarily large (and we need it to be as large as an arbitrary $W$ of our choice). 

The key insight is now to do one final iteration by applying the expansion (\ref{Ap_prep}) in turn to each of the factors $A_p(b_ju, B^n, b_ju)$ in (\ref{Ap_prep}); this results in a larger expansion containing factors $A_p(b_jb_{j'}u, B^n, b_jb_{j'}u)$ to which we again apply (\ref{Ap_prep}), and so on, iterating $L$ times, say. (Applying (\ref{Ap_prep}) is, in effect, decoupling a multilinear quantity to smaller and smaller scales.) After $L$ iterations, this results in an expansion of the form 
\beq\label{Ap_strong}
 A_p(u,B^n,u) \ll \del^{-\ep} V_{p,n}(\del)^{1 - (\sum_{\jbf} \ga_\jbf)} \times \prod_{j_1=0}^r  \cdots \prod_{j_L=0}^r A_p(\be_{\jbf} u, B^n,\beta_\jbf u)^{\ga_\jbf} \times D_p(1,B^n)^{1 - (\sum_{\jbf} \ga_\jbf)}
\eeq
where for $\jbf = (j_1,\ldots, j_L) \in [0,r]^L$, we have defined  $\be_\jbf = b_{j_1} \cdots b_{j_L}$ and $\ga_\jbf  = \ga_{j_1} \cdots \ga_{j_L}$. (The expansion (\ref{Ap_strong}) is valid for $p$ sufficiently close to the critical exponent $p_n$, and for $u$ sufficiently small relative to the $\be_\jbf$.)
This turns out to be exactly what we need: now we apply (\ref{Ap_Dp_relation}) to each of the $A_p$ factors, cancel certain terms that match advantageously, and see that 
(\ref{Ap_strong}) implies that 
\[  A_p(u,B^n,u) \ll \del^{-\ep} \del^{\eta_p(1 - u \sum_{\jbf} \be_\jbf \ga_\jbf)} D_p(1,B^n).\]
While we couldn't force a single sum $\sum_{j=0}^r b_j \ga_j$ to be large (although we said it could be made $>1$ by taking $r$ sufficiently large), if we consider instead the sum
\[ \sum_{\jbf} \be_\jbf \ga_\jbf  = ( \sum_{j=0}^r b_j \ga_j)^L,
\]
by taking $r$ and then $L$ sufficiently large we may make this as large as we like, and in particular larger than any $W$ in the hypothesis of Theorem \ref{thm_BDG_Ap}. Success! 

To summarize, in order to prove Theorem \ref{thm_BDG_Ap}, the key is to prove the relation (\ref{Ap_prep}). 
The proof of this relation is the technical heart of the Bourgain-Demeter-Guth argument.

\subsubsection{Iteration: a sketch of the protocol}\label{sec_iteration_sketch}
To sketch very briefly the protocol for proving the key relation (\ref{Ap_prep}), it is helpful to reverse perspectives, as in \cite[\S 8]{BDG16}: instead of thinking of a fixed ball size and passing to ever smaller scales of decoupling (as (\ref{Ap_prep}) currently states), we will think of a fixed decoupling scale while passing to ever larger ball sizes; thus we may aim to bound $A_p(1,\Bcal,1)$ by a relation akin to (\ref{Ap_prep}) for balls $\Bcal$ of iteratively larger radii; once this has been achieved, we can re-name the variables and reverse perspectives, to deduce (\ref{Ap_prep}).  

Let us consider for example the case of $n=3$. Recall that the ball inflation Proposition \ref{prop_ball_inflation} is good at passing from an average over smaller balls on the left to an expression for larger balls on the right, and moreover in the case of $n=3$, this works for $L^{pk/n}$ for $k=1,2$. After an initial step that is a bit different (using $k=1$ and inflation the ball from size $\del^{-1}$ to size $\del^{-2}$), the standard iteration step applies Proposition \ref{prop_ball_inflation} with $k=2$, in this fashion: if we already know how to bound $A_p(1,\Bcal,1)$ for a ball $\Bcal$ of radius $\del^{-2(\frac{3}{2})^r}$, by bounding it as as a product of $A_p$, $D_p$, and $D_{2p/3}$ factors related to balls of size $\del^{-2(\frac{3}{2})^r}$, we then prepare to apply Proposition \ref{prop_ball_inflation} again by averaging over such balls in order to cover a larger ball $\Bcal'$ of radius $\del^{-2(\frac{3}{2})^{r+1}}$. We then apply Proposition \ref{prop_ball_inflation}, and produce (after additional applications of H\"{o}lder's inequality, $\ell^2$ decoupling for $L^2$, and lower-dimensional decoupling) an upper bound for $A_p(1,\Bcal',1)$  in terms of a product of $A_p$, $D_p$, and $D_{2p/3}$ factors related to balls of size $\del^{-2(\frac{3}{2})^{r+1}}$. 

After each of these steps, the $A_p,D_p$ factors are already ``good,'' that is, in the form we want to see them appear in (\ref{Ap_prep}), while the $D_{2p/3}$ factor is ripe for the next ball-inflation step. (If we regard this iteration as a tree, the $A_p,D_p$ factors are the final ``leaves'' while $D_{2p/3}$ factors will be split again.) If we terminate this process at step $r$, the last $D_{2p/3}$ factor we have produced can be bounded  in terms of a $D_p$ factor (via H\"{o}lder, see \cite[p. 664]{BDG16}), and we arrive at an expression of the form (\ref{Ap_prep}) (after, as mentioned before, renaming variables so that the size of the ball is fixed, while the decoupling scale is commensurably smaller). 
Note that in this description of the $n=3$ case, the multiplicative factors $b_j$ in (\ref{Ap_prep}) that  arise from this process are the factors $2$ (once) or $3/2$ ($r$ times); the exponents $\ga_j$ are accrued from various interpolations and applications of H\"{o}lder's inequality.

In general, the argument for dimension $n$ includes $n-1$ types of ball inflations, increasing the exponents by multiplicative factors of $2, 3/2, 4/3, \cdots, n/(n-1)$, and these multiplicative factors form the factors $b_j$ in (\ref{Ap_prep}), while again the exponents $\ga_j$ arise from interpolation and H\"{o}lder's inequality.
This concludes our sketch of the method for obtaining the crucial relation (\ref{Ap_prep}), which enabled the iteration method to prove Theorem \ref{thm_BDG_Ap}, the final piece in the $\ell^2$ decoupling theorem for the moment curve.

\subsection{Glimmers of parallel structure}\label{sec_parallels}
Having taken a look at both the methods of efficient congruencing and $\ell^2$ decoupling, it is worth noting  parallels between the arguments.

\emph{Rescaling:}
In efficient congruencing, we called it translation-dilation invariance and used it to pass back and forth between counting solutions $(x_1,\ldots, x_{2s})$ to the Vinogradov system with $x_i \leq X$ and $x_i$ satisfying some congruence restriction modulo $q$, to solutions $(y_1,\ldots, y_{2s})$ to the Vinogradov system with $y_i \leq X/q$ and no congruence restriction.
In  decoupling, we called it affine invariance, and used it to relate decoupling at various scales (\S \ref{sec_decoupling_rescaling}). 

\emph{Multilinear estimates:}
The role of multilinear estimates in the decoupling approach is clear, but now we may notice them also in the arithmetic argument. For example, the definition (\ref{EC_bilinear}) is a bilinear expression, which we use in Lemma \ref{k2_EC_lemma2} to dominate a linear quantity, analogous to the relation (\ref{k2_bilinear_control}). More sophisticated applications of efficient congruencing use multilinear ``well-conditioned'' $m$-tuple products, e.g. \cite[Eqn 3.3.]{Woo12a}, \cite[Eqn 2.18--2.20]{Woo17b}. In each setting, the rescaling principle was critical for showing that the linear case followed from a multilinear case.

\emph{Iteration:} In both settings, we set up an iterative system that allowed us to carry out induction on scales.
In the arithmetic setting, this meant repeatedly passing from congruences modulo $p^b$ to congruences modulo $p^{\kappa^\ell b}$ for certain values of $\kappa,\ell$ depending on the degree of the Vinogradov system; we may see this as playing a parallel role to the ball inflation in \S \ref{sec_iteration_sketch}. Each method may be viewed from either of two dual perspectives: in efficient congruencing, we can think of passing to congruences modulo higher powers of $p$ as  ``inflation'' or equivalently as $p$-adic ``concentration.'' In decoupling, we can fix a scale of decoupling and inflate the ball, or by changing parameters, fix the ball size and decouple down to a smaller scale.

\emph{The role of transversality:} This leads to the question of whether a Kakeya phenomenon may be found within the structure of the efficient congruencing approach. A hint appears in our sketch in \S \ref{sec_k3_EC} of efficient congruencing for degree 3 Vinogradov systems, where we noted already the encroachment of singular solutions (that is, solutions for which certain gradient vectors are proportional modulo $p$). Wooley removes singular solutions via a ``conditioning'' step, and it is possible to interpret this as imposing a form of ($p$-adic) transversality.

The methods of efficient congruencing and decoupling reach deeply into a web of interconnected questions with long histories, and one imagines, long futures. Unifying these two methods within one over-arching framework will likely lead to new understanding of applications as well.

\subsection*{Acknowledgements}
 I express my thanks to C. Demeter, S. Guo, K. Hughes and T. Wooley  for helpful conversations, and to D.R. Heath-Brown, E. Kowalski, D. Oliveira e Silva, and P.-L. Yung for remarks on this manuscript.  In preparing these notes,  I have benefited from the expository work of many authors, in lectures, lecture notes and surveys, including online resources (not formally published) by D.R. Heath-Brown, J. Hickman and M. Vitturi, I. {\L}aba, and T. Tao, for which I express my great appreciation.

The author is partially supported by NSF grant DMS-1402121 and CAREER grant DMS-1652173, and was hosted at MSRI during a portion of the spring 2017 semester, under NSF Grant No. 1440140.
\bibliographystyle{alpha}
\bibliography{NoThBibliography}

%
%
\end{document}